%% file: unseen_13Mar2022-arxiv.tex
\documentclass[bj,preprint]{imsart}

\RequirePackage{amsthm,amsmath,amsfonts,amssymb}
\RequirePackage[numbers]{natbib}
\RequirePackage[colorlinks,citecolor=blue,urlcolor=blue,breaklinks=true]{hyperref}
\RequirePackage{graphicx}

\usepackage{breakcites}

\usepackage[pdftex,dvipsnames]{xcolor}
\usepackage{mathrsfs}

\usepackage{stmaryrd}
\usepackage{dsfont}
\usepackage{mathtools}
\usepackage{bbm}
\usepackage{bm}

\usepackage{booktabs}
\usepackage{pgf}

\startlocaldefs

\newcommand{\M}{\mathbf{M}}

\newcommand{\ddr}{\mathrm{d}}

\def\simiid{\stackrel{\mbox{\scriptsize{\rm iid}}}{\sim}}

\theoremstyle{plain}
\newtheorem{thm}{\textsc{Theorem}}
\newtheorem{definition}{\textsc{Definition}}
\newtheorem{lem}{\textsc{Lemma}}

\newtheorem{rmk}{Remark}
\newtheorem{prp}{\textsc{Proposition}}

\global\long\def\Reals{\mathbb{R}}
\global\long\def\Nats{\mathbb{N}}
\global\long\def\Ints{\mathbb{Z}}
\global\long\def\NNInts{\Ints_{+}}
\global\long\def\NNReals{\Reals_{+}}

\global\long\def\intd{\mathrm{d}}

\newcommand{\1}{\boldsymbol{1}}

\providecommand\given{} 
\newcommand\SetSymbol[1][]{
  \nonscript\,#1:\nonscript\,\mathopen{}\allowbreak}
\DeclarePairedDelimiterX\Set[1]{\lbrace}{\rbrace}%
{ \renewcommand\given{\SetSymbol[]} #1 }

\newcommand{\EE}{\mathbb{E}}
\newcommand{\var}{\mathrm{var}}
\newcommand{\PP}{\mathbb{P}}
\newcommand{\QQ}{\mathbb{Q}}

\newcommand{\KL}{\mathsf{KL}}

\newcommand{\equalDist}{\overset{d}{=}}

\newcommand{\gammaDist}{\mathrm{Gamma}}

\newcommand{\PK}{\mathrm{PK}}

\endlocaldefs

\begin{document}

\begin{frontmatter}
\title{Near-optimal estimation of the unseen under regularly varying tail populations}
\runtitle{Near-optimal estimation of the unseen under regularly varying tail populations}
\runauthor{Favaro and Naulet}

\begin{aug}
  \author[A]{\fnms{Stefano} \snm{Favaro}\thanksref{m1}\ead[label=e1]{stefano.favaro@unito.it}}
  \and
  \author[B]{\fnms{Zacharie} \snm{Naulet}\ead[label=e2]{zacharie.naulet@universite-paris-saclay.fr}}
  
  \thankstext{m1}{Also affiliated to IMATI-CNR ``Enrico  Magenes" (Milan, Italy).}

\address[A]{Department of Economics and Statistics,
  University of Torino and Collegio Carlo Alberto, Torino, 10134, Italy.\\
\printead{e1}}

\address[B]{Laboratoire de math\'ematiques d'Orsay, CNRS,
Universit\'e Paris-Saclay, 91405 Orsay, France.\\
\printead{e2}}


\end{aug}

\begin{abstract}
Given $n$ samples from a population of individuals belonging to different species, what is the number $U$ of hitherto unseen species that would be observed if $\lambda n$ new samples were collected? This is the celebrated unseen-species problem, of interest in many disciplines, and it has been the subject of recent breakthrough studies introducing non-parametric estimators of $U$ that are minimax near-optimal and consistent all the way up to $\lambda\asymp\log n$. These works do not rely on any assumption on the underlying unknown distribution $p$ of the population, and therefore, while providing a theory in its greatest generality, worst-case distributions may hamper the estimation of $U$ in concrete applications. In this paper, we strengthen the non-parametric framework for estimating $U$, making use of suitable assumptions on the underlying distribution $p$. In particular, inspired by the estimation of rare probabilities in extreme value theory, and motivated by the ubiquitous power-law type distributions in many natural and social phenomena, we make use of a semi-parametric assumption of regular variation of index $\alpha\in(0,1)$ for the tail behaviour of $p$. Under this assumption, we introduce an estimator of $U$ that is simple, linear in the sampling information, computationally efficient, and scalable to massive datasets. Then, uniformly over our class of regularly varying tail distributions, we show that the proposed estimator has provable guarantees: i) it is minimax near-optimal, up to a power of $\log n$ factor; ii) it is consistent all of the way up to $\log \lambda\asymp n^{\alpha/2}/\sqrt{\log n}$, and this range is the best possible. This work presents the first study on the estimation of the unseen under regularly varying tail distributions. Our results rely on a novel approach, of independent interest, which combines the renowned method of the two fuzzy hypotheses for minimax estimation of discrete functionals, with Bayesian arguments under Poisson-Kingman priors for $p$. A numerical illustration of our methodology is presented for synthetic and real data.
\end{abstract}

\begin{keyword}[class=MSC]
\kwd[Primary ]{62G05 }
\kwd[; Secondary ]{62C20}
\end{keyword}

\begin{keyword}
\kwd{Multinomial model} \kwd{optimal minimax estimation} \kwd{Poisson-Kingman prior} \kwd{power-law data} \kwd{regularly varying tails} \kwd{tail-index} \kwd{useen-species problem}
\end{keyword}

\end{frontmatter}


\section{Introduction}

Estimating the number of unseen species is an important problem in many disciplines. It first appeared in ecology (\citet{Fis(43),Goo(56),Cha(92),Bun(93)}), and its importance has grown in recent years driven by applications in biological sciences (\citet{Kro(99),Gao(07),Ion(09),Dal(13)}). Consider a population of individuals, with each individual being endowed with a value from a set of species labels $\mathbb{S}=\{s_{j}\}_{j\geq1}$, and for $n,m\geq1$ consider $(n+m)$ samples from the population to be modeled as a random sample $\mathbf{X}_{n+m}=(X_{1},\ldots,X_{n},X_{n+1},\ldots,X_{n+m})$ from an unknown distribution $p=\sum_{j\geq1}p_{j}\delta_{s_{j}}$ on $\mathbb{S}$. Assuming that only the first $n$ elements of $\mathbf{X}_{n+m}$ are observable, the unseen-species problem calls for estimating
\begin{displaymath}
U_{n,m}\coloneqq U_{n,m}(\mathbf{X}_{n+m})=|\{X_{n+1},\ldots,X_{n+m}\}\setminus \{X_{1},\ldots,X_{n}\}|,
\end{displaymath}
namely the number of hitherto unseen species that would be observed in $m$ additional (unobservable) sample. For $\lambda=m/n>0$, the Good-Toulmin estimator (\citet{Goo(56)}) and the smoothed Efron-Thisted estimator (\citet{Efr(76)}) are arguably the most popular estimators of $U_{n,\lambda n}$. Motivated by the increasing interest in the range $\lambda>1$, especially in biological applications, the smoothed Efron-Thisted estimator has been the subject of recent breakthrough studies (\citet{Orl(17),Wu(16), Wu(19),Pol(20)}). In particular, under the Multinomial model for the observable samples, \citet{Orl(17)} showed that a smoothed Efron-Thisted estimator is such that: i) its worst-case normalized mean-square error (NMSE) is minimax near-optimal for all $\lambda>1$; ii) it consistently estimates $U_{n,m}$ all of the way up to $\lambda\asymp\log n$, and this range is the best possible. These results do not rely on any assumption on the underlying distribution $p$, and they hold uniformly over the class $\mathcal{P}$ of distributions on $\mathbb{S}$, thus providing a grounded theory in its greatest generality. Under this non-parametric framework, however, worst-case distributions may hamper both the estimation of $U_{n,\lambda n}$ and the study of theoretical guarantees, leading to unreliable results in numerous concrete applications.

\subsection{Our contributions}

In this paper, we strengthen the non-parametric framework for estimating $U_{n,\lambda n}$, making use of suitable assumptions on the underlying distribution $p$. For $p\in\mathcal{P}$, we define the function $\bar{F}_p : (0,1) \to \NNInts$ such that
\begin{equation*}
  \bar{F}_p(x) \coloneqq \sum_{j\geq 1}\1_{\{p_j > x\} },\qquad x\in (0,1),
\end{equation*}
which characterizes the distribution of $U_{n,\lambda n}$, and we argue that reasonable assumptions on the distribution $p$ can be made by prescribing the behaviour of $\bar{F}_p(x)$ near $x\approx 0$. We denote by $K_{n}$ the number of distinct species in a random sample $(X_{1},\ldots,X_{n})$ from $p$, such that $U_{n,m}=K_{n+m}-K_{n}$. Since $\EE_p[K_{n+m}] = (n+m)\int_{(0,1)}\bar{F}_p(1-x)x^{n+m-1}\intd x$, that is $\EE_p[K_{n+m}]/(n+m)$ is the Mellin transform of $y \mapsto \bar{F}_p(1-y)$, the estimation of $K_{n+m}$ is dual to the estimation of $\bar{F}_p(x)$ for small values of $x$. Therefore, the problems of estimating $U_{n,m}$, $K_{n+m}$ and the Mellin transform of $y \mapsto \bar{F}_p(1-y)$ at large arguments, or characterizing $\bar{F}_p(z)$ for small $z$, are essentially equivalent problems. This observation allows to draw a parallel with extreme value theory (EVT) (\citet{Bin(89),Deh(06)}), where interest is in estimating probabilities of rare events; that is, given a real-valued random sample $(Z_1,\dots,Z_n)$, interest is in the estimation of the survival function $z \mapsto \PP(Z_1 > z )$ for large values of $z$, especially for $z$ larger than $\max\{Z_{1},\ldots,Z_{n}\}$. We see that $\bar{F}_p$ in the unseen-species problem plays an analogous role to that of the survival function in EVT, and therefore it is not  surprising that the estimation of $U_{n,m}$ suffers from the same issues encountered in the estimation of probabilities of rare events. In both cases the interest is in estimating $\bar{F}_p$ outside the range allowed by available data, which, unless $\lambda$ is small or a weak loss functions is considered, is impossible because of the arbitrary behaviour of $\bar{F}_p$ nearly zero. This precludes generic non-parametric assumptions in many realistic scenarios.

Inspired by the well-known problem of estimating probabilities of rare events in the context of EVT, we consider a semi-parametric assumption on $\bar{F}_p$ prescribing the behaviour of $\bar{F}_p(x)$ near $x \approx 0$. In particular, for constants $L,\,L^{\prime}>0$ and $\alpha\in(0,1)$, we consider the class of regularly varying distributions
  \begin{equation*}
\mathcal{P}(\alpha,L,L^{\prime})\coloneqq\left\{p\in\mathcal{P}\text{ : }\sup_{x\in(0,1)}\frac{|x^{\alpha}\bar{F}_{p}(x)-L|}{\sqrt{x^{\alpha}\log\left(e/x\right)}}\leq L^{\prime}\right\},
\end{equation*}
i.e. distributions with regularly varying tails of index $\alpha$ (\citet{Fel(71),Bin(89)}). We note that not all distributions with regularly
varying tails of index $\alpha$ permit a reliable estimation of the unseen for large $\lambda$. This is the reason to impose, in the definition of $\mathcal{P}(\alpha,L,L^{\prime})$, that the error rate of approximating $\bar{F}_p(x)$ by $Lx^{-\alpha}$ must decay rapidly enough when $x \to 0$. Such a constraint is common in EVT, for instance to quantify the rate of convergence of Hill's estimator (\citet{Hil(75),Hal(84),Hal(85),Dre(98)}). A distribution $p$ satisfying $\lim_{x\to 0}x^{\alpha}\bar{F}_p(x) \sim L$ for
some $(\alpha,L)\in (0,1)\times (0,\infty)$ is referred to as a power-law
distributions, and data from such a distribution are typically referred to as power-law
data. As in EVT, for power-law type distributions the value of $\EE_p[K_{n+m}]$
can be deduced from $\alpha$ via Tauberian arguments (\citet{Fel(71)}). \citet{Gne(07)} show that
$\lim_{n\rightarrow+\infty}n^{-\alpha}K_{n}=L\Gamma(1-\alpha)$ almost surely,
with $\Gamma(\cdot)$ being the Gamma function. The power-law behaviour is a
common assumption in EVT, because it entails that the law of
$\max\{Z_{1},\ldots,Z_{n}\}$ properly centered and rescaled converge to a Fr\'echet
distribution. In the unseen-species problem, interest in the power-law behaviour
is two-fold: first it enables for tractable analysis and accurate estimation of
$U_{n,\lambda n}$, and second, power-law data have readily identifiable
signatures which permits to assess the plausibility of the assumption on real
data. See Appendix~\ref{sec:elem-diagn-check} for details on checking the power-law assumption.

\begin{rmk}
  \label{rmk:1}
  One could argue that since the behaviour of $U_{n,m}$ for large $m$ is
  determined by $\bar{F}_p(x)$ for small $x$, the class
  $\mathcal{P}(\alpha,L,L')$ is too restrictive, and we may be tempted to use
  instead the weaker class
  \begin{equation*}
    \tilde{\mathcal{P}}(\alpha,L,L',x_0) \coloneqq \Set*{p \in \mathcal{P} \given \sup_{x\in (0,x_0)}\frac{|x^{\alpha} \bar{F}_p(x )-L|}{\sqrt{x^{\alpha}\log(e/x)} } \leq L'}.
  \end{equation*}
For any $x_0 \in (0,1)$ we have $\mathcal{P}(\alpha,L,L') \subset \mathcal{P}(\alpha,L,L',x_0)$ and for any 
  $p \in \mathcal{P}$, we also have for $x \geq x_0$ that
  $\bar{F}_p(x) \leq \bar{F}_p(x_0) = \sum_{j\geq 1}\1_{p_j > x_0} \leq x_0^{-1}\sum_{j\geq 1}p_j = x_0^{-1} < \infty$.
  This shows that for any $p \in \mathcal{P}$ and any $(\alpha,L)$, it must be
  that
  $\sup_{x\in (x_0,1)} |x^{\alpha} \bar{F}_p(x )-L| (x^{\alpha}\log(e/x))^{-1/2} < \infty$. That is, for all $(\alpha,L,L',x_0)$ there is a $L'' > 0$ such that
  $\mathcal{P}(\alpha,L,L') \subset \mathcal{P}(\alpha,L,L',x_0) \subset \mathcal{P}(\alpha,L,L'')$.
  Using the weaker class $\mathcal{P}(\alpha,L,L',x_0)$ in place of the class $\mathcal{P}(\alpha,L,L')$ is
  interesting only if the dependency in $(L',x_0)$ can be made explicit in the
  rates of convergence of the estimators. Unfortunately, this is known to be a tough problem, and then we make
  the choice to stick with $\mathcal{P}(\alpha,L,L')$ in order to facilitate the analysis.
\end{rmk}

The estimation of $U_{n,\lambda n}$ over $\mathcal{P}(\alpha,L,L^{\prime})$ is directly connected to the problem of estimating the tail-index $\alpha\in(0,1)$, and hence we construct estimators for both $\alpha$ and $U_{n,\lambda n}$. We assume the Multinomial model for the observable samples and, as a practical tool, we rely on Poisson-Kingman priors (\citet{Pit(03),Pit(06)}) for the underlying unknown $p\in\mathcal{P}(\alpha,L,L^{\prime})$ to obtain estimators of $\alpha$ and $U_{n,\lambda n}$ and to establish minimax lower bounds under suitable loss functions. The proposed estimators $\hat{\alpha}_{n}$ and $\hat{U}_{n,\lambda n}(\hat{\alpha}_{n})$ of $\alpha$ and $U_{n,\lambda n}$, respectively, are simple, linear in the sampling information, computationally efficient and scalable to massive datasets. In particular, uniformly over $\mathcal{P}(\alpha,L,L^{\prime})$, we show that: i) under a quadratic loss, $\hat{\alpha}_{n}$ consistently estimates $\alpha$ at rate $n^{-\alpha/2}\sqrt{\log n}$; ii) under a quadratic loss normalized by $(n\lambda)^{\alpha}$, the estimator $\hat{U}_{n,\lambda n}(\hat{\alpha}_{n})$ consistently estimates $U_{n,\lambda n}$ all of the way up to $\log \lambda\asymp n^{\alpha/2}/\sqrt{\log n}$. We also establish lower bounds for the minimax risks of estimating $\alpha$ and $U_{n,\lambda n}$. In particular, uniformly over $\mathcal{P}(\alpha,L,L^{\prime})$, we show that: i) the estimators $\hat{\alpha}_{n}$ and $\hat{U}_{n,\lambda n}(\hat{\alpha}_{n})$ are near-optimal, in the sense of matching minimax lower bounds, up to a power of $\log n$ factor; ii) the range $\log \lambda\asymp n^{\alpha/2}/\sqrt{\log n}$ is the best possible for consistently estimating $U_{n,\lambda n}$. This work presents the first study on the estimation of the unseen under regularly varying tail distributions. Our results rely on a
novel approach, of independent interest, that combines the renowned method of the two fuzzy hypotheses for minimax estimation of discrete functionals (\citet{tsybakov2009,Wu(16),Wu(19)}), with Bayesian arguments under Poisson-Kingman priors for the underlying $p\in\mathcal{P}(\alpha,L,L')$; interestingly, our approach is critical both for deriving an estimator of the tail-index $\alpha$ and for establishing minimax lower bounds. An illustration of our methodology is presented for synthetic data, text data arising from Gutenberg books and Wikipedia pages, and data from humans' electronic activities in email communications and Twitter posts. All these data have power-law behaviours, and our empirical analysis emphasizes the substantial gain in estimation, over smoothed Efron-Thisted estimators, allowed by leveraging such a power-law behaviour. This is of course expected, and we make clear that such a results do not mean that our estimator is preferable in general. The good empirical performance of our estimator on these datasets is due to their power-law behaviour, on which our estimator is designed to work at best.

Power-law data occur in many situations of scientific interest, and nowadays they have significant consequences for the understanding of natural and social phenomena. See \citet{Cla(09)}, and references therein, for a detailed account on power-law data. Besides the context of text data in natural languages (\citet{Zip(49),Can(03),Har(01)}),  
power-law phenomena have emerged for data arising from humans' electronic activities, e.g. patterns of website visits, email messages, relations and interactions on social networks, password innovation, tags in annotation systems and edits of webpages (\citet{Hub(99),Bar(05),Ryb(09),Muc(13)}). In these contexts, the estimation of $U_{n,m}$ is critical in decisions 
concerning with, e.g., collective attention monitoring, resources managements,
novelty and popularity triggering, and language learning. For power-law data,
our study shows a remarkable gain in the estimation accuracy of
$U_{n,\lambda n}$ in the large $\lambda$ regime if we leverage this prior
information on the structure of $\bar{F}_{p}$. Of course, as in the context of
EVT, estimates have to be considered with precautions, since they rely on a
modeling hypothesis that can be easily rejected if not satisfied, but it is
impossible to accept with absolute certainty on the sole basis of the observed
data (at best we can be satisfied with plausible evidences). Thus, validating
the result of the estimation shall always be made based upon external 
information about the data generating mechanism, such as biological or physical
insights, that can confirm it has an actual power-law behaviour
(\citet{Cla(09)}).

\begin{rmk}
In general, the problem of estimating the unseen in the large $\lambda$ regime do require to make modeling assumptions on the underlying $p$. In this paper, we claim that assuming $\bar{F}_p(x)$ behaves as a polynomial near zero is a convenient assumption, often realistic in many problems of practical interest, and thus is worth investigating. Other forms of modeling assumptions can be made, such as constraining the shape of the distribution. In this direction, \citet{Che(16)}, \citet{Gig(18)} and \citet{Bal(20)} investigate the same problem under the modeling assumption that $p$ is $k$-monotone. Although not directly related to the unseen-species problem, the approaches developed in \citet{Ane(17)} and \citet{Jan(20)} may be adapted to construct estimators for the unseen.
\end{rmk}

\subsection{Organization of the paper}

The paper is structured as follows. In Section \ref{sec1} we recall the Multinomial and the Poisson-Kingman partition models, and we highlight the major challenges in deriving lower bounds for the minimax risk of estimating $U_{n,\lambda n}$. Section \ref{sec2} contains our results: i) we introduce an estimator of $\alpha$ and a plug-in estimator of $U_{n,\lambda n}$, and we establish their convergence rates; ii) we study optimality of the proposed estimators by establishing lower bounds for the minimax risks of estimating $\alpha$ and $U_{n,\lambda n}$. In Section \ref{sec3} we present a numerical illustrations of our method on synthetic and real data, and Section \ref{sec4} contains a discussion of our work, including related estimation problems, and open challenges. Proofs are deferred to appendices.


\section{Partition models for $U_{n,\lambda n}$ and preliminaries}\label{sec1}

The Multinomial model for the observable samples is the most common model for estimating $U_{n,\lambda n}$ (\citet{Goo(56),Efr(76),Orl(17)}). Let $p=\sum_{j\geq1}p_{j}\delta_{s_{j}}$ be an unknown distribution on the set of
species labels $\mathbb{S}$, that is $p\in\mathcal{P}$. Under the Multinomial model it is assumed that $n\geq1$ observable samples are modeled as a random sample $\mathbf{X}_{n}=(X_{1},\ldots,X_{n})$ from  $p$, i.e.,
\begin{equation*}
X_{i}\stackrel{\text{iid}}{\sim}p\qquad i=1,\ldots,n.
\end{equation*}
That is, the random variable $X_{i}$ takes value $s_{j}\in\mathbb{S}$ with unknown probability $p_{j}$, for $j\geq1$. Species labels $s_{j}$'s are not relevant in our context, and therefore they are assumed to be fixed. Because of the discreteness of $p$, the random sample $\mathbf{X}_{n}$ from $p$ induces a random partition $\Pi_{n}$ of the set $\mathbb{N}_{n}=\{1,\ldots,n\}$ whose blocks are the (equivalence) classes induced by the equivalence relations $i\sim j\iff X_{i}=X_{j}$ almost surely. Similarly, for $m\geq1$ we denote by $\Pi_{n+m}$ the random partition of the set $\mathbb{N}_{n+m}=\{1,\ldots,n,n+1,\ldots,n+m\}$ induced by the random sample $\mathbf{X}_{n+m}=(X_{1},\ldots,X_{n},X_{n+1},\ldots,X_{n+m})$ from $p$. According to its definition, $U_{n,\lambda n}$ is uniquely determined by the random partition $\Pi_{n+m}$. We assume that the random sample $\mathbf{X}_{n}$ can be ideally extended to a sequence $\mathbf{X}=(X_{i})_{i\geq1}$, of which the first $n$ elements $\mathbf{X}_{n}$ are observable, and we denote by $\Pi=(\Pi_{n})_{n\geq1}$ the random partition of $\mathbb{N}$ induced by the sequence $\mathbf{X}$. The conditional distribution of $\Pi$ given $\mathbf{X}_{n}$ coincides with the conditional distribution of $\Pi$ given $\Pi_{n}$, and hence it is sufficient to consider a model for $\Pi_{n}$ rather than for $\mathbf{X}_{n}$. The next definition summarizes the model for the random partition $\Pi_{n}$ under $\mathcal{P}$.

\begin{definition}[Multinomial partition model]
  \label{def1}
  Let $\mathbf{X}_{n}$ be a random sample from $p\in\mathcal{P}$, and let
  $\mathbb{P}_{p}^{n}$ be the distribution of the random partition $\Pi_{n}$ of
  $\mathbb{N}_{n}$ induced by $\mathbf{X}_{n}$. Then the Multinomial partition model
  for $\Pi_{n}$ is $\{\PP_{p}^{n}\text{ : }p\in\mathcal{P}\}$, and we
  denote by $\EE_{p}^{n}$ the expectation under $\PP_{p}^{n}$. Similarly,
  we denote by $\bar{\PP}_{p}$ and $\bar{\EE}_p$ the
  distribution of the random partition $\Pi$ and the corresponding expectation under
  $\bar{\PP}_p$, respectively.
\end{definition}


Under the model of Definition \ref{def1}, the derivation of sharp lower bounds for the minimax risk of an estimator of $U_{n,\lambda n}$ is challenging, the main difficulty being that $\mathbf{X}_{n}$ contains only few information of $U_{n,\lambda n}$
(\citet{Orl(17),Wu(19)}); this is a common feature of discrete functionals of $\Pi_{n+m}$, of which the unseen $U_{n,\lambda n}$ is an instance. Although the estimation of $U_{n,\lambda n}$ will be carried under sub-models of the Multinomial partition model, i.e. over $p\in\mathcal{P}(\alpha,L,L^{\prime})$, we will also make use of a partition
model where the distribution $p$ is itself random. Such a model will serve as a practical tool to obtain an estimator of $U_{n,\lambda n}$, and then to establish a lower bound for its minimax risk. In particular, we assume
that $p$ is a discrete (almost surely) random probability measure in the broad class of Poisson-Kingman models (\citet{Pit(03)}), namely discrete random probability measures obtained by a suitable normalization of Poisson point processes. To define a Poisson-Kingman model, let $(J_{j})_{j\geq1}$ be the decreasing ordered jumps of an inhomogeneous Poisson point process on $\NNReals$ with L\'evy measure $\rho$ such that
\begin{equation}
  \label{levycond}
\int_{\mathbb{R}^{+}}\min(1,x)\rho(\intd x)<+\infty.
\end{equation}
Under the assumption \eqref{levycond} it holds that $T \coloneqq \sum_{j\geq1}J_{j}<+\infty$, i.e. finite total mass, almost surely
(\citet{Kin(93)}). A Poisson-Kingman distribution with parameter $\rho$ is defined
as the law of the random probability masses $(J_{j}/T)_{j\geq1}$, or equivalently the law of the discrete (almost surely) random probability measure
\begin{displaymath}
P \coloneqq \sum_{j\geq1}\frac{J_{j}}{T}\delta_{s_{j}}.
\end{displaymath}
We denote by $\mathcal{L}$ the class of L\'evy measures satisfying the assumption \eqref{levycond}, and we shall write $P \sim \PK(\rho)$ for abbreviating the Poisson-Kingman distribution with L\'evy measure $\rho \in \mathcal{L}$ (\citet[Chapter 3]{Pit(06)})

According to the celebrated de Finetti's representation theorem for exchangeable random variables, a random sample from a Poisson-Kingman model $P$ is part of an exchangeable sequence $\mathbf{X}=(X_{i})_{i\geq1}$ with directing (de Finetti) measure $\PK(\rho)$, that is $\lim_{n\rightarrow+\infty}n^{-1}\sum_{1\leq i\leq n}\delta_{X_{i}}=P$ almost surely (\citet[Chapter 2 and Chapter 4]{Pit(06)}). For $n\geq1$, let $\mathbf{X}_{n}$ be a random sample from $P$, that is we write 
\begin{align}\label{exch}
X_i\,|\,P & \quad\simiid\quad P \qquad i=1,\ldots,n,\\[0.2cm]
\notag P & \quad\sim\quad \PK(\rho).
\end{align}
Under the sampling model \eqref{exch}, the distribution $\PK(\rho)$ takes on the
natural interpretation as a (nonparametric) prior distribution for the unknown
$P \in\mathcal{P}$. Hence, the distribution $\PK(\rho)$ is also referred to as
Poisson-Kingman prior (\citet{Lij(10)}). Because of the (almost sure) discreteness of $P$, the
random sample $\mathbf{X}_{n}$ from $P$ induces a random partition $\Pi_{n}$ of
the set $\mathbb{N}_{n}$ whose blocks are the (equivalence) classes induced by
the equivalence relations $i\sim j\iff X_{i}=X_{j}$ almost surely. The random
partition $\Pi_{n}$ is exchangeable, namely the distribution of $\Pi_{n}$ is a
symmetric function of the sizes of its blocks, and it is referred to as
Poisson-Kingman random partition with parameter $\rho$. Moreover, the sequence
of exchangeable random partitions $\Pi=(\Pi_{n})_{n\geq1}$ is consistent in the
sense that $\Pi_{r}$ is the restriction of $\Pi_{n}$ to the first $r$ elements,
for all $r<n$ (\citet[Chapter 2 and Chapter 3]{Pit(06)}). The assumption of consistency of $\Pi$ implies that
$\Pi$ defines an exchangeable random partition of $\mathbb{N}$, where
exchangeability of $\Pi$ means that the distribution of $\Pi$ is invariant under
finite permutations of its elements.

\begin{definition}[Poisson-Kingman partition model]
  \label{def2}
  Let $\mathbf{X}_{n}$ be a random sample from $P$, where $P\sim \PK(\rho)$, and let $\PP_{\rho}^{n}$ be the distribution of the random partition $\Pi_{n}$
  of $\mathbb{N}_{n}$ induced by $\mathbf{X}_{n}$. Then the Poisson-Kingman partition
  model for $\Pi_{n}$ is $\{\mathbb{P}_{\rho}^{n}\text{ :
  }\rho\in\mathcal{L}\}$, and we denote by $\EE_{\rho}^{n}$ the expectation
  under $\mathbb{P}_{\rho}^{n}$.
\end{definition}


\section{Near-optimal estimation of $U_{n,\lambda n}$}\label{sec2}

Under the Multinomial model for the observable samples, with $p\in\mathcal{P}(\alpha,L,L^{\prime})$, we consider the problem of estimating the unseen $U_{n,\lambda n}$ for $\lambda>1$. We start by introducing an estimator of the tail-index $\alpha$ and a plug-in estimator of $U_{n,\lambda n}$, and by establishing their convergence rates under a suitable class of loss functions. We make use of the Poisson-Kingman partition model of Definition \ref{def2}, for a suitable specification of $\rho\in\mathcal{L}$. In particular, we consider the class of L\'evy measures $\mathcal{L}^{\prime}=\{\rho_{\alpha}\text{ : }\alpha\in(0,1)\}$ where
\begin{displaymath}
\rho_{\alpha}(\ddr x)=\frac{\alpha}{\Gamma(1-\alpha)}x^{-\alpha-1}\ddr x.
\end{displaymath}
$\mathcal{L}^{\prime}$ is referred to as the class of $\alpha$-stable L\'evy measures. Assuming that $P \sim \PK(\rho_{\alpha})$ (\citet{Kin(75)}), then
\begin{equation}\label{randsur}
\lim_{x\rightarrow0}x^{\alpha}\bar{F}_{P}(x)=\frac{T^{-\alpha}}{\Gamma(1-\alpha)}\quad\text{almost surely}.
\end{equation}
See \citet{Kin(75)} and \citet{Gne(07)} for a proof of \eqref{randsur}. According to \eqref{randsur}, the parameter $\alpha\in(0,1)$ of the $\alpha$-stable Poisson-Kingman model is precisely the tail-index of the random function $\bar{F}_{P}$ almost-surely. This observation suggests to make use of the Poisson-Kingman partition model $\{\mathbb{P}_{\rho_{\alpha}}^{n}\text{ : }\alpha\in(0,1)\}$ to obtain an estimator of the tail-index $\alpha$. We remark that the model $\{\mathbb{P}_{\rho_{\alpha}}^{n}\text{ : }\alpha\in(0,1)\}$ is used uniquely to construct an estimator of $\alpha$, whereas the study of such an estimator will be carried under the Multinomial partition model. The next theorem characterizes the maximum likelihood estimator $\hat{\alpha}_{n}$ of the parameter $\alpha$ under $\{\mathbb{P}_{\rho_{\alpha}}^{n}\text{ : }\alpha\in(0,1)\}$, and it guarantees its existence and uniqueness under weak conditions.

\begin{thm} \label{thm:1} Let $\Pi_{n}$ be a partition of $\Nats_n$, and denote by $C_{n,\ell}$ the number of blocks of $\Pi_n$ with at least $\ell$
  elements, for $\ell=1,\dots,n$. Then, under the model
  $\{\mathbb{P}_{\rho_{\alpha}}^{n}\text{ : }\alpha\in(0,1)\}$ the maximum
  likelihood estimator $\hat{\alpha}_{n}$ of the parameter $\alpha$ based on the
  observation of $\Pi_n$ satisfies
  $\varphi(\hat{\alpha}_n) = 0$, where $\varphi$ is such that
\begin{equation*}
\varphi_n(\alpha)\coloneqq
\sum_{k=1}^{n-1}\frac{\alpha}{k-\alpha} \cdot C_{n,k+1}- (C_{n,1}-1).
\end{equation*}
The equation $\varphi_n(\alpha)x = 0$ has a unique solution in $(0,1)$
whenever $C_{n,1}\neq n$ and $C_{n,1}\neq 1$.
\end{thm}

Theorem \ref{thm:1} shows that the maximum likelihood estimator of the parameter $\alpha\in(0,1)$ exists uniquely whenever $\Pi_n$ is not the partition consisting on $n$ singletons or the partition consisting of a single block. Therefore, without loss of generality, we assume that the estimator $\hat{\alpha}_n$ always exists; alternatively, we may assume the convention that $\hat{\alpha}_n = 0$ if $\Pi_{n}$ is the partition that consists of a single block and $\alpha=1$ if $\Pi_{n}$ is the partition that consists of $n$ blocks. Under the Multinomial partition model with $p\in\mathcal{P}(\alpha,L,L^{\prime})$, the next theorem shows that $\hat{\alpha}_{n}$ is a consistent estimator of the tail-index $\alpha$, as $n\rightarrow+\infty$, with respect to a quadratic loss function. Moreover, it holds that the convergence rate of $\hat{\alpha}_{n}$ is $n^{-\alpha/2}\sqrt{\log n}$.

\begin{thm} \label{thm:3}
For $\overline{L}> 0$ and $0 < \underline{\alpha} < \overline{\alpha} < 1$ let
$\mathbb{A}= \{(\alpha,L,L')\text{ :
} L,L' \in (0,\overline{L}),\ \alpha \in (\underline{\alpha},\overline{\alpha})\}$.
For every $\overline{L} > 0$, $0 < \underline{\alpha} < \overline{\alpha} < 1$, and $\varepsilon > 0$ there exists a constant $D > 0$ such that
\begin{equation*}
\limsup_{n\to \infty} \sup_{\substack{(\alpha,L,L')\in \mathbb{A} \\ p\in \mathcal{P}(\alpha,L,L^{\prime})}}
\bar{\mathbb{P}}_{p}\left( |\hat{\alpha}_n - \alpha|^2 > \frac{D \log(n)}{n^{\alpha}} \right) \leq \varepsilon.
\end{equation*}
\end{thm}

According to Theorem \ref{thm:1} and Theorem \ref{thm:3}, we make use of  $\hat{\alpha}_{n}$ to introduce an estimator of the unseen $U_{n,\lambda n}$ for $\lambda>1$. Because $\lambda >1$ and $p\in\mathcal{P}(\alpha,L,L^{\prime})$, it is well-known that $U_{n,\lambda n} \asymp(\lambda n)^{\alpha}$ almost-surely (\citet{Gne(07)}). This suggests to measure the performance of an estimator of
$U_{n,\lambda n}$ through the loss function
\begin{equation}\label{loss}
\ell_{\alpha}(u,v)=\frac{(u-v)^{2}}{(\lambda n)^{2\alpha}}.
\end{equation}
Observe that, under the loss functions \eqref{loss}, the null estimator $\hat{U}=0$ is such that $\ell_{\alpha}(\hat{U},U_{n,\lambda n})=O_p(1)$. Thus, we look for estimators of $U_{n,\lambda n}$ such that the loss function diminishes, in probability, as $n\rightarrow+\infty$. It is known that the number $K_{n}$ of distinct species in a random sample $\mathbf{X}_{n}$ from $p\in\mathcal{P}(\alpha,L,L^{\prime})$ is such that $\lim_{n\rightarrow+\infty}n^{-\alpha}K_{n}=L\Gamma(1 - \alpha)$ almost-surely (\citet{Gne(07)}). Accordingly, by assuming that the tail-index $\alpha$ is known, it is natural to consider as an estimator of $U_{n,\lambda n}$ the following quantity 
\begin{equation}\label{est_u}
K_{n}((1+\lambda)^{\alpha}-1).
\end{equation}
In general, the tail-index $\alpha$ is unknown, and hence it must be estimated. This
leads to introduce an estimator of $U_{n,\lambda n}$ by combining \eqref{est_u} with $\hat{\alpha}_{n}$. Here we consider the plug-in (threshold) estimator of the form
\begin{displaymath}
\hat{U}_{n,\lambda n} \coloneqq  K_{n}((1+\lambda)^{\hat{\alpha}_{n}}-1)\1_{\Set{ \log(\lambda) \leq C \sqrt{n^{\hat{\alpha}_n}/\log(n)}   }},
\end{displaymath}
where $C > 0$ is any constant. Under the Multinomial partition model of Definition \ref{def1}, with $p\in\mathcal{P}(\alpha,L,L^{\prime})$, and with respect to the class of loss
functions \eqref{loss}, the next theorem shows that
$\hat{U}_{n,\lambda n}$ consistently estimates $U_{n,\lambda n}$ all of the way
up to $\log( \lambda )\asymp \sqrt{n^{\alpha}/\log n}$. As a direct consequence of this result, we obtain a range or threshold for $\lambda>1$ in order to consistently estimate $U_{n,\lambda n}$ over the class $\mathcal{P}(\alpha,L,L^{\prime})$.

\begin{thm}\label{thm:6}
  For $\overline{L}> 0$ and $0 < \underline{\alpha} < \overline{\alpha} < 1$ let
  $\mathbb{A}= \{(\alpha,L,L')\text{ :
  } L,L' \in (0,\overline{L}),\ \alpha \in (\underline{\alpha},\overline{\alpha})\}$.
  Then for every $\bar{L} > 0$, $0 < \underline{\alpha} < \bar{\alpha} < 1$, and
  $\varepsilon > 0$ there exists a constant $D > 0$ such that
\begin{equation*}
\limsup_{n\to \infty} \sup_{\substack{(\alpha,L,L')\in \mathbb{A} \\ p\in \mathcal{P}(\alpha,L,L^{\prime})}}
\bar{\mathbb{P}}_{p}\left( \ell_{\alpha}(\hat{U}_{n,\lambda n},U_{n,\lambda n}) >
  D\min\Big\{1,\, \frac{\log(n)\log(\lambda)^2}{n^{\alpha}} \Big\} \right) \leq \varepsilon.
\end{equation*}
\end{thm}

\subsection{Optimality of the estimators $\hat{\alpha}_{n}$ and
  $\hat{U}_{n,\lambda n}$}\label{sec12}

We establish lower bounds for the minimax risks of estimating the tail-index $\alpha$ and the unseen $U_{n,\lambda n}$. Under the Multinomial model for the observable samples, with $p\in\mathcal{P}(\alpha,L,L^{\prime})$, we start by determining a lower bound for the minimax risks of estimating $\alpha$, with respect to a quadratic loss, and then we determine a lower bound for the minimax risks of estimating $U_{n,\lambda n}$, with respect to the loss function \eqref{loss}. These results, in combination with Theorem~\ref{thm:3} and Theorem~\ref{thm:5}, allows for an assessment on the optimality of the estimators $\hat{\alpha}_{n}$ and $\hat{U}_{n,\lambda n}$, where optimality is in the sense of matching minimax lower bounds. To construct minimax lower bounds we make use of the Poisson-Kingman partition model of Definition \ref{def2}. Precisely, minimax lower bounds are obtained by means of a suitable variation of the method the two fuzzy hypotheses (\citet[Chapter 2]{tsybakov2009}), which exploits the Poisson-Kingman distribution as a fuzzy prior. Our problem requires a minor adaptation over \citet[Theorem 2.14]{tsybakov2009}, mostly to deal with the problem that we aim at deriving minimax lower bounds over a class $\mathcal{P}^{\prime}\subseteq\mathcal{P}$ but our priors may have a larger support. The adaption of the method is straightforward as long as the prior charges enough $\mathcal{P}^{\prime}$. For the sake of completeness, the next theorem gives the desired extension of \citet[Theorem 2.14]{tsybakov2009}, and it will be critical for the derivation
of our minimax lower bounds.

\begin{thm}\label{thm:2}
Let $\nu_0$ and $\nu_1$ be prior distributions over $\mathcal{P}$, let $\mathcal{P}' \subseteq \mathcal{P}$, and let $\delta(\Pi,p) \in \mathbb{R}$ be a (possibly random) functional to be estimated. Write $F_0$ and $F_1$ for the joint distributions of $(\Pi,P)$, respectively under the prior $\nu_0$ and $\nu_1$. If there exist $c \in \mathbb{R}$ and $\kappa,\xi \in (0,1)$ such that $\min\{F_0(\delta(\Pi,P) \leq c),\, F_1(\delta(\Pi,P) \geq c + 2\varepsilon) \} \geq 1 - \kappa$ and $\min\{\nu_0(\mathcal{P}'),\nu_1(\mathcal{P}')\} \geq 1 - \xi$, then it holds
\begin{equation*}
\inf_{\hat{\delta}}\sup_{p \in \mathcal{P}'}\bar{\mathbb{P}}_{p}(|\hat{\delta}(\Pi_n) - \delta(\Pi,p)| > \varepsilon)\geq \frac{1}{2}(1 - \|\nu_0(\mathbb{P}_{P}^n) - \nu_1(\mathbb{P}_{P}^n)\|_{\mathrm{TV}}) - \kappa - \xi,
\end{equation*}
where $\nu_j(\mathbb{P}_{P}^n)(\cdot)$ are the mixture distributions $\int_{\mathcal{P}}\mathbb{P}_{P}^n(\cdot)\nu(\ddr P)$, $j=0,1$, and where the infimum is taken over all measurable maps $\hat{\delta}$ that are solely function of $\Pi_n$.
\end{thm}

The idea consists in the application of Theorem \ref{thm:2} to the case where $\nu_0$ and $\nu_1$ are Poisson-Kingman distributions with L\'evy measures $\rho_0$ and $\rho_1$, respectively. In this case, $\nu_j(\mathbb{P}_{P}^n)$ is equal to the distribution $\PP_{\rho_j}^n$ of the Poisson-Kingman partition model of Definition~\ref{def2}. The main challenge then reduces to find an upper bound on $\|\PP_{\rho_0}^n - \PP_{\rho_1}^n\|_{\mathrm{TV}}$ for general L\'evy measures $\rho_0$ and $\rho_1$. In particular, with the help of Pinsker's inequality, it is enough to upper bound $\KL(\PP_{\rho_0}^n;\PP_{\rho_1}^n)$ which turns out to be more convenient. In the next proposition, we derive an upper bound on $\KL(\PP_{\rho_0}^n;\PP_{\rho_1}^n)$ depending only on the moments
\begin{equation*}
\mu_j(u,\ell) \coloneqq%
\int_0^{\infty} x^{\ell}e^{-ux}\rho_j(\ddr x),\qquad j=0,1,\quad \ell=1,\dots,n,
\end{equation*}
and the statistic $\M_{n} \coloneqq (M_{n,1},\dots,M_{n,n})$, where $M_{n,\ell}$ denotes the number of blocks of size $\ell$ in the partition $\Pi_n$.
\begin{prp} \label{thm:4} Let $G_n \sim \gammaDist(n,1)$ be independent of the
  Poisson process $(J_j)_{j\geq 1}$ used to construct the random probability
  measure $P = \sum_{j\geq 1}\frac{J_j}{T}\delta_{s_j}$, and let
  $U_n \coloneqq G_n/T$. Then,
\begin{align*}
    \KL(\PP_{\rho_0}^n;\PP_{\rho_1}^n)%
    &\leq%
      \sum_{\ell=1}^n%
      \mathbb{E}_{\rho_0}^n\Big[ M_{n,\ell} \Big\{%
      -\log \frac{\mu_1(U_n,\ell)}{\mu_0(U_n,\ell)} +
      \frac{\mu_1(U_n,\ell)}{\mu_0(U_n,\ell)} -1%
      \Big\} \Big].
  \end{align*}
\end{prp}

We apply Theorem \ref{thm:2} and Proposition \ref{thm:4} to determine a lower bound for the minimax risk of estimating the tail-index $\alpha$. In general, both Theorem~\ref{thm:2} and Proposition~\ref{thm:4} are of independent interest, and they may be applied in the context of the estimation of other discrete functionals under the Multinomial model for the observable samples, with $p\in\mathcal{P}(\alpha,L,L^{\prime})$. We refer to Section \ref{sec4} for a discussion. Under the Multinomial partition model for the observables samples, with $p\in\mathcal{P}(\alpha,L,L^{\prime})$, and with respect to a quadratic loss function, the next theorem shows that the estimator $\hat{\alpha}_{n}$ is near-optimal up to a power of $\log n$ factor.

\begin{thm}\label{thm:5}
For $\overline{L}> 0$ and $0 < \underline{\alpha} < \overline{\alpha} < 1$ let
$\mathbb{A}= \{(\alpha,L,L')\text{ :
} L,L' \in (0,\overline{L}),\ \alpha \in (\underline{\alpha},\overline{\alpha})\}$.
For every $0 < \underline{\alpha} < \overline{\alpha} < 1$ there exists a constant $D > 0$ such that
  \begin{equation*}
    \liminf_{n\to\infty}
    \inf_{\hat{\alpha}_n} \sup_{\substack{(\alpha,L,L')\in \mathbb{A} \\ p\in \mathcal{P}(\alpha,L,L^{\prime})}}%
    \bar{\mathbb{P}}_{p}\left( |\hat{\alpha}_{n} - \alpha|^2 \geq \frac{D}{n^{\alpha}\log(n)^2} \right) \geq \frac{1}{4},
  \end{equation*}
  where the infimum is with respect to all measurable maps $\hat{\alpha}_n$ depending only on $\Pi_n$.
\end{thm}

Under the Multinomial partition model with $p\in\mathcal{P}(\alpha,L,L^{\prime})$, and with respect to the loss function \eqref{loss}, the next theorem shows that the estimator $\hat{U}_{n,\lambda n}(\hat{\alpha}_{n})$ is near-optimal up to a power of $\log n$ factor. The proof of the theorem consists in establishing that the estimation of $U_{n,\lambda n}$ is a harder problem than the estimation of $\alpha$ if $\lambda$ is sufficiently large, thus reducing the problem to the estimation of the tail-index. 

\begin{thm}\label{thm:7}
For $\overline{L}> 0$ and $0 < \underline{\alpha} < \overline{\alpha} < 1$ let
$\mathbb{A}= \{(\alpha,L,L')\text{ :
} L,L' \in (0,\overline{L}),\ \alpha \in (\underline{\alpha},\overline{\alpha})\}$,
and let $b > 0$ be such that it holds true $\log\lambda\geq b\log n$. Then for
every $0 < \underline{\alpha} < \bar{\alpha} < 1$ there exists constant $D > 0$ such that
\begin{equation*}
\liminf_{n\to\infty} \inf_{\hat{U}} \sup_{\substack{(\alpha,L,L')\in \mathbb{A} \\ p\in \mathcal{P}(\alpha,L,L^{\prime})}}%
    \bar{\mathbb{P}}_{p}\left( \ell_{\alpha}(\hat{U},U_{n,\lambda n}) \geq D\min\left\{1,\, \frac{\log(\lambda)^2}{n^{\alpha}\log(n)^2} \right\} \right) \geq \frac{1}{8},
  \end{equation*}
where the infimum is with respect to all measurable maps $\hat{U}$ depending only on $\Pi_n$.
\end{thm}

According to Theorem \ref{thm:7}, no estimator can estimate $U_{n,\lambda n}$ in the range $\log(\lambda) \geq n^{\alpha/2}\log(n)$; moreover, according to Theorem \ref{thm:6}, this range or threshold (of estimability) is almost attained by the estimator $\hat{U}_{n,\lambda n}$. As a corollary of Theorem \ref{thm:7} and Theorem \ref{thm:6} it holds that, up to polylog factors, the range $\log \lambda\asymp n^{\alpha/2} $ is the best possible range for consistently estimating $U_{n,\lambda n}$, and such a range is achieved by the estimator $\hat{U}_{n,\lambda n}$. In the range $\log(\lambda) \gtrsim n^{\alpha/2}\log(n)$, the trivial estimator $\hat{U} = 0$ is a minimax estimator, with no extra $\log(n)$ factor. That is, in the range $\log(\lambda) \gtrsim n^{\alpha/2}\log(n)$ the best we can do to estimate $U_{n,\lambda n}$ is to do nothing. We observe that Theorem \ref{thm:7} holds true for large values of $\lambda$, typically $\lambda \gtrsim n^{b}$. This is because of the arguments used in the proof, which establish that the problem of estimating $U_{n,\lambda n}$ is at least as hard as estimating the tail-index $\alpha$; these two problems are essentially equivalent, provided that $\lambda$ is not too small. For small values of $\lambda$, however, we do not believe that the two problems are equivalent. In particular, if $\lambda < 1$ we expect that smoothed Efron-Thisted estimator (\citet{Goo(56),Efr(76),Orl(17)}) may have a stronger theoretical and experimental performance than our estimator. Furthermore, we believe the two problems may be equivalent at a much lower range of $\lambda$ than that provided by Theorem \ref{thm:6}. Understanding this range requires to obtain an exact matching for the upper and lower bounds of the risk in the estimation of $\alpha$, which remains a challenging problem. This constitutes an interesting problem to investigate in future research.


\section{Numerical illustrations}\label{sec3}

Power-law type distributions are arguably the most common examples of
distributions in the class of regularly varying tail distributions
$\mathcal{P}(\alpha,L,L^{\prime})$. Power-law type distributions occur in many
situations of scientific interest, and nowadays have significant consequences
for the understanding of natural and social phenomena (\citet{Cla(09)}). In this section, we present a numerical
illustration of our methodology for synthetic data from Zipf distributions (a
favorable scenario), modified Zipf distributions with two tail indices (a less
favorable scenario), text data from Gutenberg books and Wikipedia pages, and
data from humans' electronic activities in email communication and Twitter
posts. We compare our estimator to smoothed Efron-Thisted estimators, which may
be considered as the state-of-the-art estimators for the unseen-species problem
(\citet{Goo(56),Efr(76),Orl(17)}).

\subsection{Synthetic data : favorable scenario}
\label{sec:simulated-data}

\begin{figure}[!htb]
    \label{fig:simulated-tailindex}

    \centering
    \resizebox{0.8\linewidth}{!}{\input{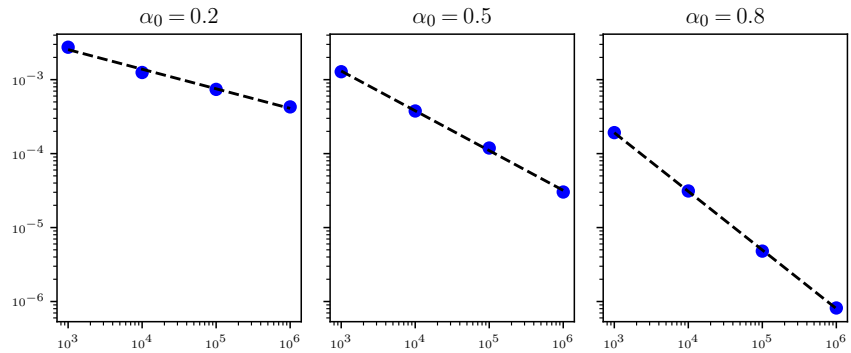}}%
    \caption{Monte Carlo estimates of $\EE[(\hat{\alpha}_n - \alpha_0)^2]$ based
      on 1000 replicates of $n$ data points generated (i.i.d) from
      $\mathrm{Zipf}(1/\alpha_0)$.The blue points represent the
      estimated values of $\EE[(\hat{\alpha}_n - \alpha_0)^2]$ for
      $n \in \{1000,10000,100000,10000\}$. The dashed line is a fit
      corresponding to $\EE[(\hat{\alpha}_n - \alpha_0)^2] = e^{-4.14}n^{-0.27}$
      (left), $\EE[(\hat{\alpha}_n - \alpha_0)^2] = e^{-2.92}n^{-0.54}$
      (center), and $\EE[(\hat{\alpha}_n - \alpha_0)^2] = e^{-3.08}n^{-0.79}$
      (right).}
  \end{figure}

\begin{figure}[!htb]
    \label{fig:tailindex-hists}

    \centering
    \resizebox{0.8\linewidth}{!}{\input{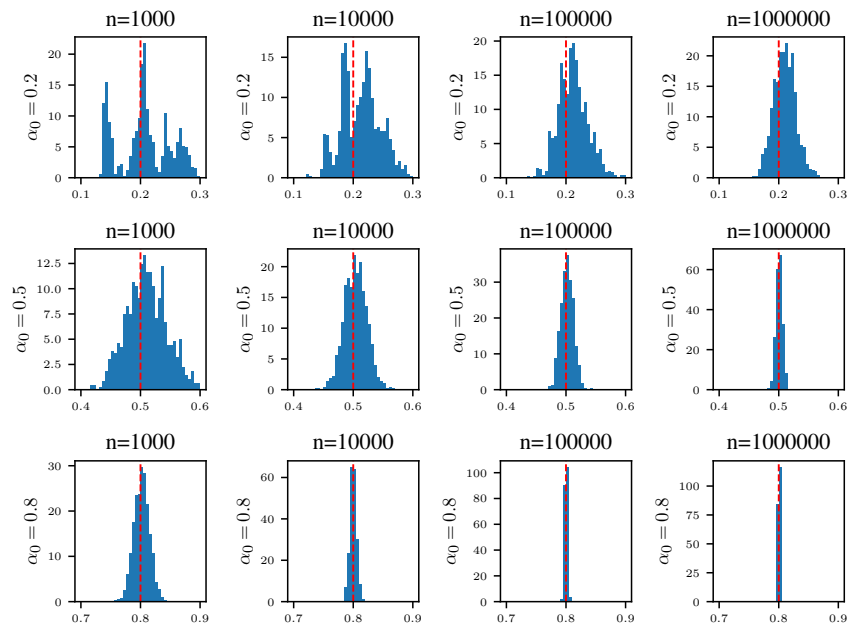}}%
    \caption{Histograms of $\hat{\alpha}_n$ based
      on 1000 replicates of $n$ data points generated (i.i.d) from $\mathrm{Zipf}(1/\alpha_0)$.}
  \end{figure}

The most favorable scenario for our estimators is undoubtedly when the model is
well-specified, i.e. $\Pi_n$ has a Poisson-Kingman distribution with Lévy
measure $\rho_{\alpha}$ for some $\alpha \in (0,1)$. As we are mainly interested
in data generated according to the Multinomial model, we leave the ideal case
apart and consider instead a situation where the model
$\Set{\PP_{\rho_{\alpha}} \given \alpha \in (0,1)}$ is misspecified, but in a
benign fashion. A prototypical example of such benign misspecification is when
$(X_1,\dots,X_n)$ generated (i.i.d) from a Zipf distribution (\citet{Cla(09)}).
For $s > 1$, the $\mathrm{Zipf}(s)$ distribution on $\NNInts$ with parameter $s$
has probability mass function $P(j)\propto j^{-s}$ for all $j \in \NNInts$. The
parameter $s$ controls the tail behaviour of Zipf distribution: the smaller $s$
the heavier is the tail of the distribution, i.e., the smaller $s$ the larger
the fraction of species with low frequency. In particular, it follows that the
$\mathrm{Zipf}(s)$ distribution has a tail-index $\alpha = 1/s$ and belongs to
$\mathcal{P}(\alpha,L,L')$ for some $L,L' > 0$. Therefore, we expect synthetic
Zipf data to be a favorable situation for our estimators. We start by an
empirical validation of the maximum likelihood estimator $\hat{\alpha}_n$ of the
tail-index. We compute Monte Carlo (MC) estimates of
$\EE[(\hat{\alpha}_n - \alpha_0)^2]$, where the expectation is understood with
respect to data generated from a $\mathrm{Zipf}(1/\alpha_0)$ distribution. The
expectation is approximated using $1000$ MC samples. Simulations were run for
$n \in \{1000, 10000, 100000, 1000000\}$ and $\alpha_0 \in \{0.2,0.5,0.8\}$.
Estimated values of $\EE[(\hat{\alpha}_n - \alpha_0)^2]$ are reported in
Table~\ref{tab:simulated-tailindex} and plotted in
Figure~\ref{fig:simulated-tailindex}. About the plots, we have adjusted a curve
$\log \EE[(\hat{\alpha}_n - \alpha_0)^2] = -r \log(n) + C$ using least squares
estimation. We found estimated values of $r$ that are $0.27$, $0.54$, and
$0.79$, respectively for $\alpha_0 = 0.2$, $\alpha_0= 0.5$, and
$\alpha_0 = 0.8$, which are coherent with our theoretical finding that
$\hat{\alpha}_n$ converges at rate $n^{-\alpha_0}$.
Figure~\ref{fig:tailindex-hists} contains the histograms of $\hat{\alpha}_n$,
showing that the distribution of $\hat{\alpha}_n$ does concentrate on
$\alpha_0$, and this concentrations happens faster as $\alpha_0$ gets larger.

  \begin{table}[!htb]
    \centering
    \caption{Monte Carlo estimates of $\EE[(\hat{\alpha}_n - \alpha_0)^2]$ based on 1000 replicates of $n$ data points generated (i.i.d) from
    $\mathrm{Zipf}(1/\alpha_0)$}
    \label{tab:simulated-tailindex}
    \begin{tabular}{c|c|c|c}
      \toprule
      $n$%
      &$\alpha_0 = 0.2$%
      &$\alpha_0 = 0.5$%
      &$\alpha_0 = 0.8$\\
      \midrule
      $1000$%
      & $2.75\cdot 10^{-3}$%
      & $1.28\cdot 10^{-3}$%
      & $1.92\cdot 10^{-4}$\\
      $10000$%
      & $1.3\cdot 10^{-3}$%
      & $3.77\cdot 10^{-4}$%
      & $3.12\cdot 10^{-5}$\\
      $100000$%
      & $7.39 \cdot 10^{-4}$%
      & $1.19\cdot 10^{-4}$%
      & $4.80\cdot 10^{-6}$\\
      $1000000$%
      & $4.28\cdot 10^{-4}$%
      & $3.02\cdot 10^{-5}$%
      & $8.18 \cdot 10^{-7}$\\
      \bottomrule
    \end{tabular}
  \end{table}

\begin{figure}[!htb]
    \label{fig:unseen-largelambda}

    \centering
    \resizebox{7cm}{!}{\input{unseen_large.pgf}}%
    \caption{Monte Carlo estimates of $\tilde{R}_{n,\lambda n}$ for our
      estimator $\hat{U}_{n,\lambda n}$ based on 100 replicates of $n$ data points generated (i.i.d) from
      $\mathrm{Zipf}(1/\alpha_0)$, with $\alpha_0 = 0.8$. The
      growth of $\tilde{R}_{n,\lambda n}$ in $\lambda$ is logarithmic, which is
      consistent with the theory.}
  \end{figure}
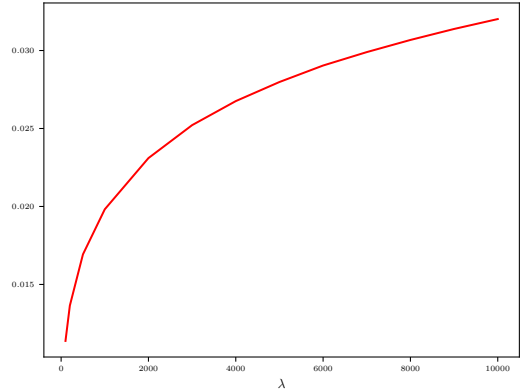

We consider the same MC setting to make an empirical validation of $\hat{U}_{n,\lambda n}$.  We compute MC estimates of the risk
\begin{displaymath}
\tilde{R}_{n,\lambda n} \coloneqq \frac{\EE[(\hat{U}_{n,\lambda n} - U_{n,\lambda n})^2]}{\EE[U_{n,\lambda n}]^2}.
\end{displaymath}
In particular, we fix $n=1000$, and increasing values of $\lambda$ from $\lambda=1.1$ to $\lambda= 40$. We compare our estimator with the state-of-the-art estimators of \citet{Orl(17)},
which are smoothed Efron-Thisted estimator based on Poisson smoothing or Binomial smoothing distributions. Estimated values of $\tilde{R}_{n,\lambda n}$ are reported in Table~\ref{tab:simulated-tailindex} and plotted in Figure~\ref{fig:simulated-tailindex}. In Table~\ref{tab:simulated-tailindex}, we only report the best value for smoothed Efron-Thisted estimators, which is not necessarily attained using the same smoothing in the various scenarios. Our results show that smoothed Efron-Thisted estimators can estimate $U_{n,\lambda n}$ with good accuracy only when $\lambda$ is small, which is expected, whereas by estimating the tail-index the estimation of the unseen remains accurate even for large values of $\lambda$. In particular, $\hat{U}_{n,\lambda n}$ can estimate the unseen with more accuracy than smoothed Efron-Thisted estimators in the situations where the data come from a distribution with a tail-index. It does not make clear, however, how large $\lambda$ can be taken. To illustrate that the unseen can eventually be estimated with good accuracy even for large $\lambda$, we ran the same simulation for $\hat{U}_{n,\lambda n}$, $\alpha_0 = 0.8$, $n=1000$, but with $\lambda$ ranging from $100$ to $10000$. From Figure~\ref{fig:unseen-largelambda} it is clear that  $\tilde{R}_{n,\lambda n}$ grows logarithmically with $\lambda$, which agrees with our theoretical findings. Hence, in this example, even for values of $\lambda$ as large as $\lambda = 10000$ the accuracy of estimating the unseen remains satisfactory.

  \begin{table}[!htb]
    \centering
    \caption{Monte Carlo estimates of $\tilde{R}_{n,\lambda n}$ for various
      estimators of the unseen based on 1000 replicates of $n$ data points generated (i.i.d) from $\mathrm{Zipf}(1/\alpha_0)$. Simulations are made on the basis of $n=1000$ observed
      samples. The column $\EE[U_{n,\lambda n}]$ is a MC estimate of the expected
      number of unseen, the column FN gives a MC estimate of
      $\tilde{R}_{n,\lambda n}$ for $\hat{U}_{n,\lambda n}$, and
      the column O gives the best value of $\tilde{R}_{n,\lambda n}$ obtained
      using the three estimators proposed in \citet{Orl(17)}.}
    \label{tab:simulated-unseen}
    \resizebox{0.8\linewidth}{!}{
      \begin{tabular}{c|c|c|c|c|c|c|c|c|c}
        \toprule
        & \multicolumn{3}{c|}{$\alpha_0=0.2$}%
        & \multicolumn{3}{c|}{$\alpha_0=0.5$}%
        & \multicolumn{3}{c}{$\alpha_0=0.8$}\\
        \midrule
        $\lambda$%
        & $\EE[U_{n,\lambda n}]$%
        & FN%
        & O%
        & $\EE[U_{n,\lambda n}]$%
        & FN%
        & O%
        & $\EE[U_{n,\lambda n}]$%
        & FN%
        & O\\%
        \midrule
        1.1%
        & 1.20%
        & 0.99%
        & 2.35
        & 62.53
        & $2.56\cdot 10^{-2}$
        & $3.97\cdot 10^{-2}$
        & 1741
        & $9.66\cdot 10^{-4}$
        & $1.35 \cdot 10^{-3}$
        \\
        1.5
        & 1.54
        & 0.80
        & 1.75
        & 80.77
        & $2.24\cdot 10^{-2}$
        & $4.89 \cdot 10^{-2}$
        & 2323
        & $8.75\cdot 10^{-4}$
        & $8.48 \cdot 10^{-3}$
        \\
        2
        & 1.86
        & 0.69
        & 1.64
        & 101.63
        &  $2.13\cdot 10^{-2}$
        &  $7.28\cdot 10^{-2}$
        & 3025
        & $8.30\cdot 10^{-4}$
        & $4.12 \cdot 10^{-2}$
        \\
        3
        & 2.41
        & 0.56
        & 1.41
        & 138.39
        &  $1.92\cdot 10^{-2}$
        &  $1.51\cdot 10^{-1}$
        & 4364
        & $7.65\cdot 10^{-4}$
        & $1.66 \cdot 10^{-1}$
        \\
        4
        & 2.86
        & 0.49
        & 1.21
        & 170.81
        &  $1.88\cdot 10^{-2}$
        &  $2.46\cdot 10^{-1}$
        & 5638
        & $7.61\cdot 10^{-4}$
        & $2.87 \cdot 10^{-1}$
        \\
        5
        & 3.24
        & 0.45
        & 1.40
        & 200.43
        &  $1.82\cdot 10^{-2}$
        &  $2.20\cdot 10^{-1}$
        & 6861
        & $7.65\cdot 10^{-4}$
        & $2.61 \cdot 10^{-1}$
        \\
        6
        & 3.55
        & 0.43
        & 1.98
        & 227.71
        &  $1.84\cdot 10^{-2}$
        &  $2.86\cdot 10^{-1}$
        & 8044
        & $7.79\cdot 10^{-4}$
        & $3.43 \cdot 10^{-1}$
        \\
        8
        & 4.06
        & 0.40
        & 1.64
        & 276.05
        &  $1.83\cdot 10^{-2}$
        &  $4.04\cdot 10^{-1}$
        & 10312
        & $8.05\cdot 10^{-4}$
        & $4.80 \cdot 10^{-1}$
        \\
        10
        & 4.54
        & 0.38
        & 1.37
        & 319.49
        &  $1.84\cdot 10^{-2}$
        &  $5.00\cdot 10^{-1}$
        & 12483
        & $8.30\cdot 10^{-4}$
        & $5.83 \cdot 10^{-1}$
        \\
        15
        & 5.42
        & 0.35
        & 1.09
        & 414.50
        &  $1.85\cdot 10^{-2}$
        &  $6.61\cdot 10^{-1}$
        & 17598
        & $8.87\cdot 10^{-4}$
        & $7.42 \cdot 10^{-1}$
        \\
        20
        & 6.15
        & 0.33
        & 0.99
        &494.98
        &  $1.92\cdot 10^{-2}$
        &  $7.55\cdot 10^{-1}$
        & 22399
        & $9.52\cdot 10^{-4}$
        & $8.25 \cdot 10^{-1}$
        \\
        30
        & 7.26
        & 0.32
        & 0.95
        & 631.13
        &  $2.00\cdot 10^{-2}$
        &  $8.52\cdot 10^{-1}$
        & 31377
        & $1.03\cdot 10^{-3}$
        & $9.04 \cdot 10^{-1}$
        \\
        40
        & 8.12
        & 0.31
        & 0.95
        & 746.16
        &  $2.12\cdot 10^{-2}$
        &  $8.99\cdot 10^{-1}$
        & 39783
        & $1.10\cdot 10^{-3}$
        & $9.39 \cdot 10^{-1}$
        \\
        \bottomrule
      \end{tabular}
    }
  \end{table}

  \begin{figure}[!htb]
    \centering
    \resizebox{0.8\linewidth}{!}{\input{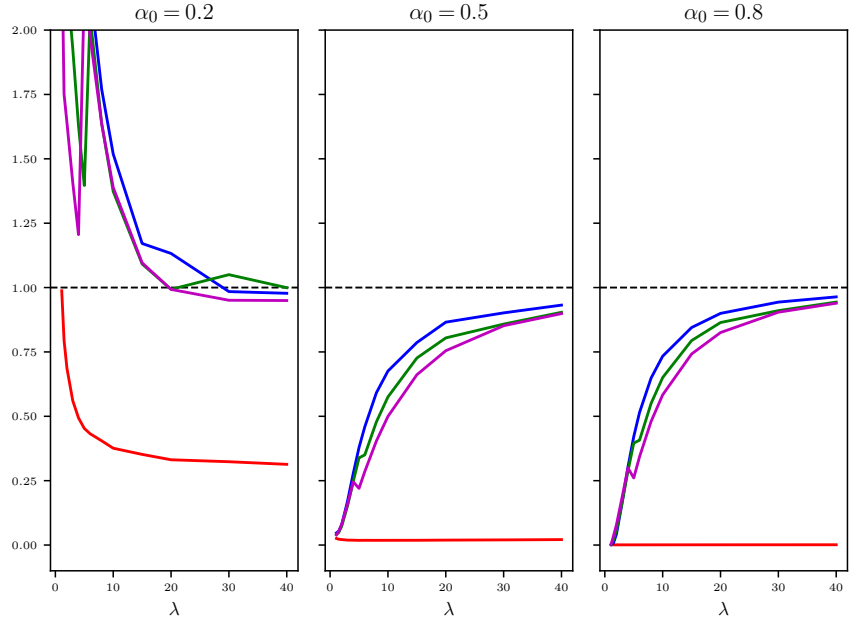}}%
    \caption{Monte Carlo estimates of $\tilde{R}_{n,\lambda n}$ for various
      estimators of the unseen based on 1000 replicates of $n$ data points generated (i.i.d) from $\mathrm{Zipf}(1/\alpha_0)$. Simulations are made on the basis of $n=1000$ observed
      samples. The red curve is the estimator
      $\hat{U}_{n,\lambda n}$ and the three other curves
      correspond to the three estimators proposed in \citet{Orl(17)},
      respectively for a Poisson smoothing (blue) and two variants of a Binomial
      smoothing (green and magenta) of the Good-Toulmin estimator.}
  \end{figure}

\subsection{Synthetic data: less favorable scenario}\label{sec:synthetic-data-new}

We consider a collection of $n$ data points modeled as i.i.d. random samples $(X_1,\dots,X_n)$ coming from a distribution $p$ such that $p_j \propto j^{-1/\alpha}$ for $j \in \{1,\dots,J\}$ and $p_j \propto j^{-1/\beta}$ for $j > J$; we shall refer to such distribution as ``double Zipf distribution''. We assume that $0 < \beta < \alpha < 1$. Such a double Zipf distribution is less favorable for our estimator of the unseen, because it is constructed in such a way that it has a true tail-index of $\beta$, while if $J$ is too large it will look like data are
generated from a $\mathrm{Zipf}(1/\alpha)$ distribution. Indeed, the
Lemma~\ref{lem:double-zipf} below, shows that the double Zipf distribution will
always be in $\mathcal{P}(\beta,L,L')$ for some $L,L' > 0$, but the larger $J$
is taken, the larger $L'$ will be. Performances of our estimators at finite
sample size $n$, however, are strongly dependent on the value $L'$ (the smaller
the better).

\begin{lem}
  \label{lem:double-zipf}
  For any $1 > \alpha > \beta > 0$ and $J \in \Nats$, let
  $p \propto \sum_{j=1}^Jj^{-1/\alpha}\delta_j + \sum_{j>J}j^{-1/\beta}\delta_j$.
  There exist constants $L,C,J_0 > 0$ such that for all $J \geq J_0$,
  \begin{equation*}
    \sup_{x \in (0,1)} \frac{|x^{\beta}\bar{F}_p(x) - L|}{\sqrt{x^{\beta}\log(e/x)}}%
    \leq C \cdot (\alpha - \beta) \sqrt{\frac{J^{2\alpha/\beta} \log(J)}{J}}.
  \end{equation*}
\end{lem}

\begin{figure}[!htb]
  \label{fig:double-zipf}

  \centering
  \resizebox{6cm}{4.5cm}{\input{mc_estimate_risk_alpha_double_zipf.pgf}}%
  \resizebox{6cm}{4.5cm}{\input{mc_estimate_risk_unseen_double_zipf.pgf}}
  \caption{On the left: Monte Carlo estimates of
    $\EE[(\hat{\alpha}_n - \beta)^2]$ as a function of $J$ based on 1000
    replicates of $n$ data points generated (i.i.d) from a double Zipf
    distributions with parameters $(1/0.5,1/0.4,J)$. On the right: Monte Carlo
    estimates of $\tilde{R}_{n,\lambda n}$ for our estimator
    $\hat{U}_{n,\lambda n}(\hat{\alpha}_n)$ as a function of $J$ based on 100
    replicates of $n$ data points generated (i.i.d) from a double Zipf
    distributions with parameters $(1/0.5,1/0.4,J)$.}
\end{figure}
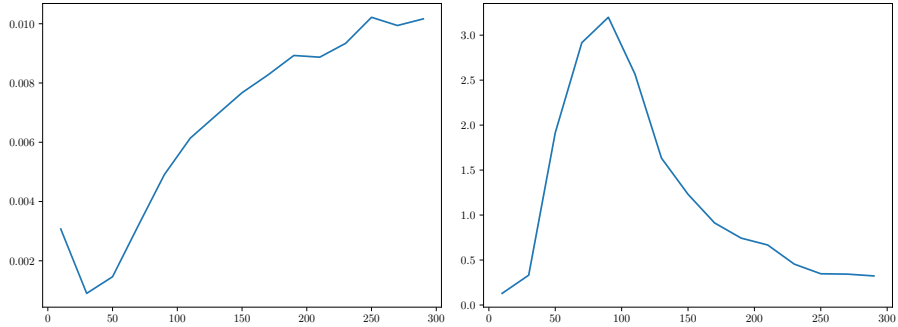

The simulations displayed in Figure~\ref{fig:double-zipf} illustrate how the risks for
$\alpha$ and for the unseen deteriorate as a function of $J$ for fixed $n$. We
take $n= 1000$, $\alpha= 0.5$, $\beta= 0.4$, $\lambda= 50$ and $J$ ranging from
$10$ to $290$. For both $\hat{\alpha}_n$ and $U_{n,\lambda n}$ we compute MC
estimates of the risks (as defined in the previous section) using $1000$ samples
for the risk of estimating $\beta$, and $100$ samples for the unseen [to make
the simulation run in short time]. We see that the risk
$\bar{\EE}_p[(\hat{\alpha}_n - \beta)^2]$ behaves as expected; that is is small
when $J$ is reasonably small and increases with $J$, until the point where it
gets of the order of $(\alpha - \beta)^2 = 0.01$. The plot for
$\tilde{R}_{n,\lambda n}$ might look more surprising at first, but it is indeed
what is expected. We see three different regions for the risk. First, if $J$ is
sufficiently small, then $\hat{\alpha}_n$ estimates $\beta$ with good accuracy
and the number of unseen is of order $(\lambda n)^{\beta}$; everything goes well
in this region. Then, when $J$ gets larger it becomes impossible to estimate
$\beta$, and we rather have $\hat{\alpha}_n \approx \alpha$, but the number of
unseen is still of order $(\lambda n)^{\beta}$; whence a very high risk.
Finally, when $J$ gets very large, we again have
$\hat{\alpha}_n \approx \alpha$, but the number of unseen is also of oder
$(n\lambda)^{\alpha}$. Basically everything happens as if $p$ was a
$\mathrm{Zipf}(1/\alpha)$ distribution; thus everything goes well again here.
We emphasize that the reason the risk gets apparently better for large $J$ is
just an artifact of the simulation that uses a fixed value of $n$ and $\lambda$.

\subsection{Real data}\label{sec:real-data}

We consider text data arising from Gutenberg books and Wikipedia pages, and data from humans' electronic activities in email communication and Twitter posts. The email dataset (eMail) consists of a collection of data from the email communication activity of a Department at Universit\`a degli Studi di Padova in the years 2012 and 2013 (\citet{For(14)}); the dataset is in the form of a table, where each row contains information of (sender, receiver, timestamp); for our analysis we set the sender identity to be a species, and the number of emails sent from a sender to be the frequency of the species. The Twitter dataset (Twitter) consists of a collection data from the communication activity of Twitter (\citet{Tri(14),Mon(17)}); the dataset is in the form of a table, where each row contains information of (timestamp, hashtag, user); for our analysis we set an hashtag to be a species, and the number of tweets containing an hashtag to be the frequency of the species. The Wikipedia (Wiki) and Gutenberg (Gut) datasets consist of a collection of data from Wikipedia pages and Gutenberg books (\citet{Tri(14),Mon(17)}); for our analysis we set a word to be a species, and the number of occurrence of a word to be the frequency of the species. Summaries of these datasets are reported in Table~\ref{tab:summary-realdata}.

\begin{figure}[!htb]
  \label{fig:unseen-realdata}
  \centering%
  \resizebox{0.8\linewidth}{!}{\input{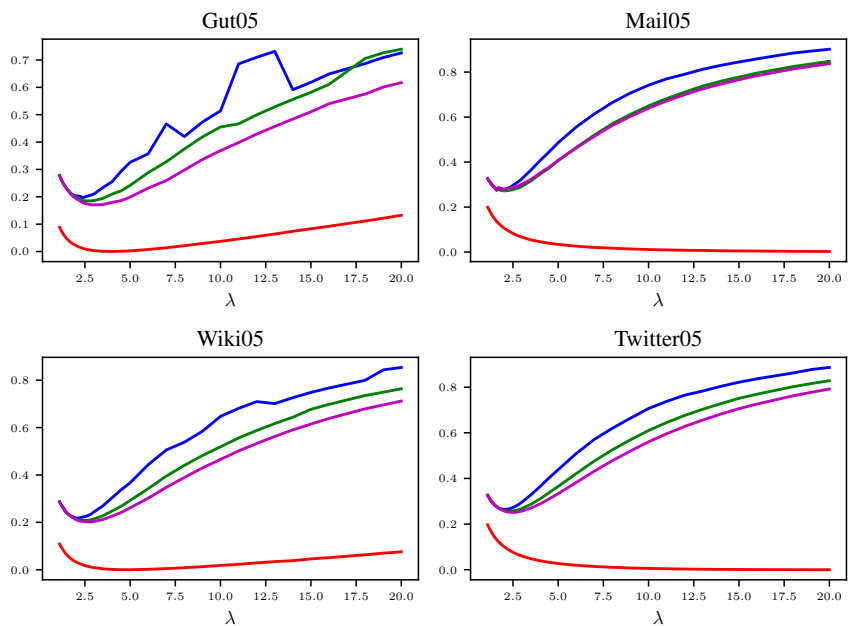}}%
  \caption{$\widehat{MSE}$ for various estimators of the unseen on real datasets
    (exact procedure detailed in Section~\ref{sec:real-data}). The red curve is the
    estimator $\hat{U}_{n,\lambda n}$ and the three other curves
    correspond to the three estimators proposed in \citet{Orl(17)}, respectively
    for a Poisson smoothing (blue) and two variants of a Binomial smoothing
    (green and magenta) of the Good-Toulmin estimator.}
\end{figure}

For each dataset, we make an empirical evaluation of $\hat{U}_{n,\lambda n}$ as follows. Given that the dataset has size $N$, see Table~\ref{tab:summary-realdata}, for each value of $\lambda$ we obtain a sample of size $n = \lfloor N/(1+\lambda) \rfloor$ from the dataset by sampling individuals without replacement. Then, from this sample, we estimate $U_{n,\lambda n}$ in the dataset and we report the normalized error $E = (\hat{U}_{n, \lambda n} - U_{n,\lambda n})/U_{n,\lambda n}$. We repeat this process $100$ times in order to obtain $100$ copies $E_1,\dots,E_{100}$ and we compute the average squared error $\widehat{MSE} = \frac{1}{100}\sum_{1\leq i\leq 100}E_i^2$. We also apply the same procedure to the three smoothed Efron-Thisted estimators proposed in \citet{Orl(17)}. Results of our experiments are available in Figure~\ref{fig:unseen-realdata}, from which we see that for all datasets considered there is an interest in leveraging the power-law behaviour and estimating the tail-index to predict the unseen, as opposed to a fully agnostic approach \textit{à la} Efron-Thisted.

\begin{table}[!htb]
  \centering
  \caption{Summary statistics of the real datasets.}
  \label{tab:summary-realdata}
  \begin{tabular}{l|c|c|c|c}
    \toprule
    & Gut%
    & eMail%
    & Wiki%
    & Twitter\\
    \midrule
    Size $N$ of the dataset
    & $6314483$%
    & $345469$%
    & $1480277$
    & $1734849$\\
    Number of species%
    & $147801$%
    & $69462$%
    & $137030$%
    & $608644$\\
    \bottomrule
  \end{tabular}
\end{table}

\section{Discussion}\label{sec4}

The unseen-species problem dates back to the seminal works of \citet{Fis(43)}, \citet{Goo(56)} and \citet{Efr(76)}, and it has been the subject of recent breakthrough studies by \citet{Orl(17)}, \citet{Wu(16),Wu(19)} and \citet{Pol(20)}. In this paper, we offered new advances on these recent studies, presenting the first work on the estimation of $U_{n,\lambda n}$ under the assumption that the underlying unknown distribution $p$ has regularly varying tails of index $\alpha\in(0,1)$. This is motivated by the ubiquitous power-law type distributions, which nowadays occur in many natural and social phenomena (\citet{Cla(09)}). Under a semi-parametric assumption of  regular variation for the tail behaviour of $p$, we introduced an estimator of $U_{n,\lambda n}$ and we showed that it has the following theoretical guarantees: i) it is minimax near-optimal up to a power of $\log n$ factor; ii) it is consistent all of the way up to $\log \lambda\asymp n^{\alpha/2}/\sqrt{\log n}$, and this range is the best possible. The proposed estimator is simple, linear in the sampling information, computationally efficient and scalable to massive datasets, and its provable guarantees hold uniformly over the class of regularly varying tail distributions. From an empirical perspective, it is shown that our estimator outperforms existing ones on several synthetic and real datasets. To the best of our knowledge, this is the first work on the estimation of $U_{n,\lambda n}$ under assumptions on the underlying unknown distribution $p$. Our results rely on a novel approach, of independent interest, which combines the method of the two fuzzy hypotheses for minimax estimation of discrete functionals (\citet{tsybakov2009,Wu(16),Wu(19)}) with Bayesian arguments under Poisson-Kingman priors for $p$.

Our semi-parametric assumption entails that $|\bar{F}_p(x) - Lx^{-\alpha}| \lesssim x^{-\alpha/2}\log(e/x)^{1/2}$, which
guarantees that $\bar{F}_p(x)$ approaches $Lx^{-\alpha}$ fast enough when
$x \to 0$. An interesting problem is to study the case
where $x^{\alpha}\bar{F}_p(x) \to L$ is true but convergence happens at a slower
rate than $x^{-\alpha/2}\log(e/x)^{1/2}$; say at rate $x^{-\beta}$ for
$\beta \in (0,\alpha/2)$. In general, we can adapt the proof of
Theorem~\ref{thm:3} to obtain the following weaker Proposition.

\begin{prp}
  \label{pro:mle-weaker-bias}
  Suppose there exist $\alpha \in (0,1)$, $\beta \in (0,\alpha)$, and $L > 0$
  such that $\bar{F}_p(x) = Lx^{-\alpha} + O(x^{-\beta})$. Then, for every
  $\varepsilon > 0$ there is a constant $D > 0$ such that
  \begin{equation*}
    \bar{\PP}_p\Big( |\hat{\alpha}_n - \alpha|^2 > D\max(n^{-\alpha},n^{-2\beta}) \Big) \leq \varepsilon.
  \end{equation*}
\end{prp}

According to Proposition \ref{pro:mle-weaker-bias}, when convergence of $x^{\alpha}\bar{F}_p(x) \to L$ happens at the rate $x^{-\beta}$, $\beta \in (0,\alpha/2)$, we can still estimate the tail-index $\alpha$, but at a deteriorated rate of $n^{-2\beta}$ in place of $n^{-\alpha}\log(n)$. Rates for the estimation of the unseen $U_{n,\lambda n}$can be adapted in the same way. It is however unclear how to obtain matching minimax lower bounds under this weaker modeling assumption. We cannot exclude that using Poisson-Kingman partitions as fuzzy hypotheses produces too ``nice'' hypotheses to get matching rates. Another obvious, but more challenging question, would be
to investigate the necessity of regular variation itself. We know that regular variations are not necessary, as shown in the works by \citet{Che(16),Gig(18),Bal(20)}, but it would be of interest to know if weaker and non-parametric power-law modeling assumptions are feasible, assuming for instance that $x^{-\alpha_1}\lesssim \bar{F}_p(x) \lesssim x^{-\alpha_2}$ for some $0 < \alpha_1 < \alpha_2 < 1 $. We conjecture that this is only possible if $\lambda$ is small and/or the loss function weak enough, such as the NMSE considered in \citet{Orl(17)}.

The unseen-species problem belongs to a broad class of discrete functional estimation problems, commonly referred to as ``species-sampling problems'' (\citet{Bal(21)}). These problems uses $n$ random samples from an unknown underlying distribution $p$, and call for estimating  features of $p$ or features of $\lambda n$ additional unobservable samples. Recent noteworthy works on ``species sampling problems'' concerned with the estimation of, e.g., the support size (\citet{val(13),Wu(19)}), the entropy (\citet{Jia(15),Wu(16)}), the missing mass (\citet{Oha(12),Mos(19),Ben(17),Fad(18)}) and the risk of disclosure (\citet{Cam(20)}). Interest in these quantities first appeared in ecology, and it has grown in the recent years driven by challenging applications in biological and physical sciences, statistical machine learning, engineering, theoretical computer science, information theory and forensic DNA analysis. Most of the works of ``species sampling problems'' do not rely on any assumption on the underlying distribution $p$. The sole exceptions are \citet{Oha(12)} and \citet{Fad(18)}, which contain a preliminary study of the problem of consistent estimation of the missing mass under the assumption of regularly varying tails for $p$, and leave open the problem of minimax optimal estimation. Our work paves the way to study ``species sampling problems'', beyond the unseen-species problem, under the assumption that $p$ has regularly varying tails. We believe that the techniques developed in our work are of direct application to the problem of disclosure risk assessment considered in \citet{Cam(20)} and to the problem of estimating the number of unseen species with multiplicity considered in \citet{Hao(20)}.


\section*{Acknowledgements}

The authors are grateful to the Editor (Professor Mark Podolsky), the Associate Editor and three anonymous Referees for all their comments, corrections, and numerous suggestions that improved remarkably the paper. Stefano Favaro and Zacharie Naulet received funding from the European Research Council (ERC) under the European Union's Horizon 2020 research and innovation programme under grant agreement No 817257. Stefano Favaro gratefully acknowledge the financial support from the Italian Ministry of Education, University and Research (MIUR), ``Dipartimenti di Eccellenza" grant 2018-2022.

\appendix

\numberwithin{equation}{section}
\numberwithin{thm}{section}
\numberwithin{lem}{section}
\numberwithin{prp}{section}

\section{Upper bounds for the estimators of the tail-index and the unseen}
\label{sec:upper-bounds-estim}

This appendix contains the proof of the Theorems~\ref{thm:1}, \ref{thm:3},
\ref{thm:6} and Proposition~\ref{pro:mle-weaker-bias}.

\subsection{Proof of Theorem~\ref{thm:1}}
\label{sec:thm:1}

It is a well-known fact that $\alpha$-stable Poisson Kingman paritition model is
equivalent to the \textit{two parameters Chinese Restaurant process} with choice
of parameters $(\alpha,0)$; see (\citet[see Section~3.2 for the definition of
the two parameters process, and Section~4.2 for the equivalence]{Pit(06)}). Let
write $M_{n,\ell}$ for the numbers of blocks of size $\ell$ in the partition
$\Pi_n$. Obviously $\M_n \coloneqq (M_{n,1},\dots,M_{n,n})$ entirely determines
the law of the partition $\Pi_n$, so we might as well assume we observe $\M_n$.
Then by \cite[Theorem~3.2]{Pit(06)}, the probability of observing
$\M_n = (m_1,\dots,m_n)$ under the parameter $\alpha$ is given by
\begin{equation*}
  p_{\alpha}(m_1,\dots,m_n)%
  = \frac{\alpha^{k-1}\prod_{i=0}^{k-2}(1+i) \prod_{j=1}^n\big( \prod_{i=0}^{j-2}(1-\alpha+i) \big)^{m_j}}{\prod_{i=0}^{n-2}(1+i)},
\end{equation*}
where $k = \sum_{j=1}^nm_j$. So any maximizer of $\alpha \mapsto p_{\alpha}(\M_n)$ is equivalently a
maximizer of
\begin{equation*}
  F(\alpha)%
  \coloneqq
  (K_n - 1)\log(\alpha) + \sum_{j=1}^n \sum_{i=0}^{j-2}M_{n,j}\log(1 - \alpha + i).
\end{equation*}
Differentiating the previous with respect to $\alpha$ gives that $\hat{\alpha}_n$ if it
exists must be solution to
\begin{align*}
  0%
  &= \frac{K_n-1}{\alpha}%
    - \sum_{j=1}^n \sum_{i=0}^{j-2} \frac{M_{n,j}}{1 - \alpha+i}
  = \frac{K_n-1}{\alpha}%
    - \sum_{i=1}^{n-1} \frac{C_{n,i+1}}{i - \alpha}.
\end{align*}
This proves that $\hat{\alpha}_n$ must be solution to $\varphi(\alpha) = 0$. Differentiating $F$
one more times gives
\begin{align*}
  F''(\alpha)%
  &= -\frac{K_n-1}{\alpha^2} - \sum_{i=1}^{n-1}\frac{C_{n,i+1}}{(i-\alpha)^2}.
\end{align*}
This establish that $F''(\alpha) < 0$ for all $\alpha \in (0,1)$, so if there is a solution
to $F'(\alpha) = 0$, equivalently to $\varphi_n(\alpha) = 0$, then it is a unique maximiser. We
observe that $C_{n,1} \geq C_{n,2} \geq \dots$. Unless $C_{n,2} = 0$ or $K_n = 1$, it is true that
$\lim_{\alpha \to 0}\varphi_n(\alpha) = -(K_n-1) < 0$ and $\lim_{\alpha\to 1}\varphi_n(\alpha) = +\infty$. Then it must be that
$\varphi_n(\alpha) = 0$ has exactly one solution provided $C_{n,2} \ne 0$ and $K_n \ne 1$. Finally, we note
that $C_{n,1} = K_n$, so $C_{n,2} = 0$ is equivalent $K_n = n$, which is
equivalent to $\Pi_n$ being the partition consisting only on singletons. Obviously
$K_n \ne 1$ is equivalent to $\Pi_n$ not being the partition with a single block.

\subsection{Proof of Theorem~\ref{thm:3}}
\label{sec:thm:3}

We remark that $\alpha \mapsto \alpha^{-1}\varphi_n(\alpha)$ is the score function
of the model $\Set{\PP_{\rho_{\alpha}} \given \alpha \in (0,1)}$. We define
$\Phi_{n,p}(\alpha) \coloneqq \bar{\EE}_p[\varphi_n(\alpha)]$ and
$\alpha^{*}_n(p)$ the solution to the equation $\Phi_{n,p}(\alpha) = 0$.
We note that the model is misspecified, such that there is no reason to have
$\alpha^{*}_n(p) = \alpha$ when $p \in \mathcal{P}(\alpha,L,L')$. We show in a
first time in Proposition~\ref{pro:control-bias-klmin} that distributions in
$\mathcal{P}(\alpha,L,L')$ are sufficiently well behaved to control the bias
$|\alpha_n^{*}(p) - \alpha|$. In a second time, we establish in
Proposition~\ref{pro:6} a uniform concentration result for
$|\hat{\alpha}_n- \alpha_n^{*}(p)|$ [using classical arguments for maximum
likelihood estimation; see \cite{vaart:2000}]. The proof of Theorem~\ref{thm:3}
follows by combining the results of Propositions~\ref{pro:control-bias-klmin}
and~\ref{pro:6}.

\begin{prp}
  \label{pro:control-bias-klmin}
  For all $0 < \underline{\alpha} < \overline{\alpha} < 1$ and all
  $\overline{L}$,
  \begin{equation*}
    \limsup_{n\to\infty}%
    \sup_{\substack{(\alpha,L,L')\in \mathbb{A}\\p\in \mathcal{P}(\alpha,L,L')}}%
    \frac{|\alpha_n^{*}(p) - \alpha|}{n^{-\alpha/2}\sqrt{\log(n)} }%
    \leq \frac{5L'}{L}.
  \end{equation*}
\end{prp}
\begin{proof}
  We first show that $|\Phi_{n,p}(\alpha)|$ must remain small if
  $p\in \mathcal{P}(\alpha,L,L')$. By Lemma~\ref{lem:asympt-cumul-freq}, for all
  $p$ and all $\alpha \in (0,1)$
  \begin{equation*}
    \Phi_{n,p}(\alpha)%
    = 1+ \sum_{k=2}^n \frac{k \alpha}{k-(1+\alpha)}\binom{n}{k}\int_0^1 \bar{F}_p(x)x^{k-1}(1-x)^{n-k}\intd x%
    - n \int_0^1 \bar{F}_p(x)(1-x)^{n-1}\intd x.
  \end{equation*}
  We decompose $\bar{F}_p$ in two parts such that
  $\bar{F}_p(x) = Lx^{-\alpha} + \delta(x)$, where by assumption we have that
  $|\delta(x)| \leq L' \sqrt{x^{-\alpha}\log(e/x)}$. We recall that
  $\int_0^1 x^{-\alpha + k -1}(1-x)^{n-k}\intd x = \frac{\Gamma(k-\alpha)\Gamma(n-k+1)}{\Gamma(n+1-\alpha)}$
  for all $k=1,\dots,n$ and all $n\geq 1$. Therefore,
  \begin{align*}
    \Phi_{n,p}(\alpha)
    &= 1 + \frac{L\alpha n!}{\Gamma(n+1-\alpha)}\sum_{k=2}^n \frac{\Gamma(k-\alpha-1)}{\Gamma(k)}%
      - \frac{L n! \Gamma(1-\alpha)}{\Gamma(n+1-\alpha)}\\
    &\quad%
      + \sum_{k=2}^n \frac{k\alpha}{k-(1+\alpha)}\binom{n}{k}%
      \int_0^1 \delta(x)x^{k-1}(1-x)^{n-k}\intd x%
      -n \int_0^1\delta(x)(1-x)^{n-1}\intd x.
  \end{align*}
  We note that
  \begin{align*}
    \sum_{k=2}^n \frac{\Gamma(k-\alpha-1)}{\Gamma(k)}%
    &= \frac{\Gamma(1-\alpha)}{\alpha} - \frac{\Gamma(n-\alpha)}{\alpha \Gamma(n)}
  \end{align*}
  which gives the estimate
  \begin{align*}
    |\Phi_{n,p}(\alpha)|%
    &\leq 1 + \frac{L\alpha n! \Gamma(n-\alpha)}{\Gamma(n+1-\alpha)\Gamma(n)}%
      + n \int_0^1|\delta(x)|(1-x)^{n-1}\intd x\\
    &\quad%
      +\sum_{k=2}^n \frac{k\alpha}{k-(1+\alpha)}\binom{n}{k} \int_0^1|\delta(x)|x^{k-1}(1-x)^{n-k}\intd x.%
  \end{align*}
  Then, we remark that the function $x \mapsto x/[x - (1+\alpha)]$ is monotone
  decreasing on $(1+\alpha,\infty)$, which implies that
  $k/(k-(1+\alpha)) \leq 2/(1 - \alpha)$ for all $k\geq 2$. It follows the
  estimate
  \begin{align*}
    &\sum_{k=2}^n \frac{k\alpha}{k-(1+\alpha)}\binom{n}{k} \int_0^1|\delta(x)|x^{k-1}(1-x)^{n-k}\intd x\\%
    &\qquad \leq \frac{2\alpha}{1-\alpha}\int_0^1 |\delta(x)| \sum_{k=1}^n \binom{n}{k}x^{k-1}(1-x)^{n-k} \intd x\\
    &\qquad= \frac{2\alpha}{1-\alpha}\int_0^1 |\delta(x)| \frac{1 - (1-x)^n}{x}\intd x
  \end{align*}
  and then
  \begin{align*}
    |\Phi_{n,p}(\alpha)|%
    &\leq 1 + \frac{L\alpha n! \Gamma(n-\alpha)}{\Gamma(n+1-\alpha)\Gamma(n)}%
      + n \int_0^1|\delta(x)|(1-x)^{n-1}\intd x\\
    &\quad%
      +\frac{2\alpha}{1-\alpha}\int_0^1|\delta(x)| \frac{1 - (1-x)^n}{x}\intd x.
  \end{align*}
  By Stirling's formula and by the Lemma~\ref{lem:5}, we obtain that,
  \begin{align*}
    \limsup_{n\to \infty}\sup_{\substack{(\alpha,L,L')\in \mathbb{A}\\p\in \mathcal{P}(\alpha,L,L')}}%
    \frac{(1-\alpha)|\Phi_{n,p}(\alpha)|}{5L'\Gamma(1-\alpha/2)\cdot n^{\alpha/2} \sqrt{\log(n)} }%
    \leq 1.
  \end{align*}

  Finally, we show that $\Phi_{n,p}(\alpha') = 0$ implies that $\alpha'$ must be
  close to $\alpha$ when $p \in \mathcal{P}(\alpha,L,L')$. In fact, from the
  definition of $\varphi_n$ we see that
  $\Phi_{n,p}'(\alpha') = \sum_{k=1}^{n-1} \frac{k}{(k-\alpha')^2}\bar{\EE}_p[C_{n,k+1}] \geq \frac{\bar{\EE}_p[C_{n,2}]}{(1 - \alpha')^2}$.
  By Lemma~\ref{lem:asympt-cumul-freq}, for $n$ large enough it will be true
  that $\bar{\EE}_p[C_{n,2}] \gg 0$, so that $\Phi_{n,p}$ is monotone increasing
  on $(0,1)$. It follows from a Taylor expansion that for any
  $\alpha' \leq \alpha - \varepsilon_{n,\alpha}$ [$\varepsilon_{n,\alpha}$ a decreasing sequence
  of real numbers going to zero, to be chosen accordingly]
  \begin{align}
    \label{eq:3}
    \Phi_{n,p}(\alpha')%
    &\leq \Phi_{n,p}(\alpha - \varepsilon_{n,\alpha})%
      \leq \Phi_{n,p}(\alpha) - \inf_{z \in (0,\varepsilon_{n,\alpha})}\Phi_{n,p}'(\alpha -z)\varepsilon_{n,\alpha}%
      \leq \Phi_{n,p}(\alpha) - \frac{\bar{\EE}_p[C_{n,2}] \varepsilon_{n,\alpha}}{(1 - \alpha)^2}.
  \end{align}
  Similarly if $\alpha' \geq \alpha + \varepsilon_{n,\alpha}$,
  \begin{align*}
    \Phi_{n,p}(\alpha')%
    \geq \Phi_{n,p}(\alpha) + \frac{\bar{\EE}_p[C_{n,2}]\varepsilon_{n,\alpha}}{(1 - \alpha + \varepsilon_{n,\alpha})^2}.
  \end{align*}
  It follows,
  \begin{align*}
    \inf_{|\alpha' - \alpha| \geq \varepsilon_{n,\alpha}}|\Phi_{n,p}(\alpha')|%
    &\geq \frac{\bar{\EE}_p[C_{n,2}]\varepsilon_{n,\alpha}}{(1 - \alpha + \varepsilon_{n,\alpha})^2}%
      - |\Phi_{n,p}(\alpha)|.
  \end{align*}
  By the equation~\eqref{eq:3} and the Lemma~\ref{lem:unif-expect-cumulfreq},
  \begin{align*}
    &\liminf_{n\to \infty}%
    \inf_{\substack{(\alpha,L,L')\in \mathbb{A}\\p\in \mathcal{P}(\alpha,L,L')}}%
    \inf_{|\alpha' - \alpha| \geq \varepsilon_{n,\alpha}}|\Phi_{n,p}(\alpha')|\\%
    &\qquad\geq%
      \liminf_{n\to \infty}%
      \inf_{\substack{(\alpha,L,L')\in \mathbb{A}\\p\in \mathcal{P}(\alpha,L,L')}}%
    \frac{5L'\Gamma(1-\alpha/2)\cdot n^{\alpha/2} \sqrt{\log(n)}}{1-\alpha} \Big[%
     \frac{L n^{\alpha/2}\Gamma(1-\alpha)\cdot \varepsilon_{n,\alpha}}{ 5L'\Gamma(1-\alpha/2)\sqrt{\log(n)} }
    -1\Big]
  \end{align*}
  which will be strictly greater than zero whenever
  $\varepsilon_{n,\alpha} \geq \frac{5L'}{L}n^{-\alpha/2}\sqrt{\log(n)}$. Since
  by construction $\Phi_{n,p}(\alpha_n^{*}(p)) = 0$, it is necessary that
  \begin{align*}
    \limsup_{n\to\infty}%
    \sup_{\substack{(\alpha,L,L')\in \mathbb{A}\\p\in \mathcal{P}(\alpha,L,L')}}%
    \frac{|\alpha_n^{*}(p) - \alpha|}{n^{-\alpha/2}\sqrt{\log(n)} }%
    \leq \frac{5L'}{L}.
  \end{align*}
\end{proof}

\begin{prp}
  \label{pro:6}
  For every $\varepsilon > 0$ there exists a constant $B$, depending only on
  $\varepsilon$ and $\mathbb{A}$, such that
  \begin{equation*}
    \limsup_{n\to \infty} \sup_{\substack{(\alpha,L,L')\in \mathbb{A} \\ p\in \mathcal{P}(\alpha,L,L^{\prime})}}
    \bar{\PP}_p\Big(|\hat{\alpha}_n - \alpha_n^{*}(p)| > Bn^{-\alpha/2} \Big) \leq \varepsilon.
  \end{equation*}
\end{prp}
\begin{proof}
  For all $0 < \delta \leq 1/2$ and all $p$, let define the events
  $\Omega_{n}(p;\delta) \coloneqq \Set{\bm{X} \given \big|C_{n,2} - \bar{\EE}_{p}[C_{n,2}]\big| \leq \delta \bar{\EE}_p[C_{n,2}] }$.
  Then, by Chebychev's inequality and because $M_{n,2} \leq C_{n,2}$
  almost-surely, the Lemma~\ref{lem:var-cumulfreq} guarantees that
  \begin{equation*}
    \bar{\PP}_p(\Omega_n(p;\delta)^c)%
    \leq \frac{\overline{\var}_p(C_{n,2})}{\delta^2 \bar{\EE}_p[C_{n,2}]^2}%
    \leq \frac{2 \bar{\EE}_p[M_{n,2}]}{\delta^2 \bar{\EE}_p[C_{n,2}]^2}%
    \leq \frac{2}{\delta^2\bar{\EE}_p[C_{n,2}] }.
  \end{equation*}
  We deduce from the last display and the Lemma~\ref{lem:unif-expect-cumulfreq}
  that whenever $\delta^2n^{\alpha} \to \infty$ as $n \to \infty$,\begin{equation*}
    \lim_{n\to \infty} \sup_{\substack{(\alpha,L,L')\in \mathbb{A} \\ p\in \mathcal{P}(\alpha,L,L^{\prime})}}
    \bar{\PP}_p(\Omega_n(p;\delta)^c) = 0.
  \end{equation*}
  We define other events
  $\Omega_n'(p;\varepsilon) \coloneqq \Set{\bm{X} \given |\varphi(\alpha_n^{*}(p))|^2 \leq A_p\bar{\EE}_p[C_{n,1}]\varepsilon^{-1}}$,
  where we let $A_p \coloneqq \max\{1 , \alpha_n^{*}(p)/(1-\alpha_n^{*}(p))^2\}$.
  We remark that
  $\bar{\EE}_p[\varphi(\alpha_n^{*}(p))] = \Phi_p(\alpha_n^{*}(p)) = 0$. Hence
  by Chebychev's inequality and the Proposition~\ref{pro:7}, for all $p$,
  \begin{equation*}
    \bar{\PP}_p(\Omega_n'(p;\varepsilon)^c)%
    \leq \frac{\overline{\var}_p(\varphi(\alpha_{*}(p)))}{A_p \bar{\EE}_p[C_{n,1}]\varepsilon^{-1}}%
    \leq \varepsilon.
  \end{equation*}
  We finish the proof by showing that on the event
  $\Omega_n(p;\delta) \cap \Omega_n'(p;\varepsilon)$, it must be that
  $|\hat{\alpha}_n - \alpha_n^{*}(p)| \leq x_n$ for some $x_n \geq 0$ to be
  determined at the end of the proof. In fact, observe that
  $\varphi'(\alpha) = \sum_{k=1}^{n-1} \frac{k}{(k-\alpha)^{2}}C_{n,k+1} \geq C_{n,2}/(1-\alpha)^2$.
  Therefore $\varphi$ is strictly monotone increasing on $(0,1)$ when we are on
  the event $\Omega_n(p;\delta)$. Consequently, for all
  $\alpha \geq \alpha^{*}_n + x_n$, by a Taylor expansion of $\varphi$
  [remarking that $\varphi'' \geq 0$ on $\Omega_n(p;\delta)$ with the same
  argument as before]
  \begin{equation*}
    \varphi(\alpha)%
    \geq \varphi(\alpha_n^{*}(p) + x_n)%
    \geq \varphi(\alpha_n^{*}(p)) + \varphi'(\alpha_n^{*}(p))x_n%
    \geq \varphi(\alpha_n^{*}(p)) + \frac{C_{n,2}}{(1-\alpha_n^{*}(p))^2}x_n;
  \end{equation*}
  Similarly for all $\alpha \leq \alpha_n^{*} - x_n$,
  \begin{equation*}
    \varphi(\alpha)%
    \leq \varphi(\alpha_n^{*}(p)) - \frac{C_{n,2}}{(1 -\alpha_n^{*}(p) + x_n)^2}x_n.
  \end{equation*}
  It follows that on the event
  $\Omega_n(p;\delta) \cap \Omega_n'(p;\varepsilon)$, for all $\alpha$ such that
  $|\alpha - \alpha_n^{*}(p)| \geq x_n$ we must have
  \begin{equation*}
    |\varphi(\alpha)|%
    \geq \frac{(1-\delta)\bar{\EE}_p[C_{n,2}]}{(1 - \alpha_n^{*}(p) + x_n)^2}x_n%
    - \sqrt{\frac{A_p\bar{\EE}_p[C_{n,1}]}{\varepsilon}}.
  \end{equation*}
  By choosing $x_n = Bn^{-\alpha/2}$ for a sufficiently large
  constant $B$, the rhs of the last display will be strictly greater than $0$
  uniformly for all $p \in \mathcal{P}(\alpha,L,L')$ for all
  $(\alpha,L,L') \in \mathbb{A}$, because of
  Proposition~\ref{pro:control-bias-klmin} [which guarantees that $A_p$ remains
  bounded for $n$ large enough] and because of
  Lemma~\ref{lem:unif-expect-cumulfreq} [which controls $\bar{\EE}_p[C_{n,k}]$
  uniformly in $p$]. But, by construction $\varphi(\hat{\alpha}_n) = 0$, so on
  $\Omega_n(p;\delta)\cap \Omega_n'(p;\varepsilon)$ it is necessary that
  $|\hat{\alpha}_n - \alpha_n^{*}(p)| \leq B n^{-\alpha/2}$.
\end{proof}

\begin{lem}
  \label{lem:5}
  For all $0 < \underline{\alpha} < \overline{\alpha} < 1$, the following
  uniform limit results are true:
  \begin{align*}
    \lim_{n\to \infty} \sup_{(\alpha,L,L')\in \mathbb{A}} \Bigg|\frac{\int_0^1L'x^{-\alpha/2}\sqrt{\log(e/x)}\frac{1-(1-x)^n}{x}\intd x }{\frac{L'\Gamma(1-\alpha/2)}{\alpha/2}\cdot n^{\alpha/2}\sqrt{\log(n)}}%
    - 1\Bigg| = 0.
  \end{align*}
  and
  \begin{align*}
    \lim_{n\to\infty}\sup_{(\alpha,L,L')\in \mathbb{A}}\Bigg|%
    \frac{n\int_0^1L'x^{-\alpha/2} \sqrt{\log(e/x)}(1-x)^{n-1}\intd x}{L'\Gamma(1-\alpha/2)\cdot n^{\alpha/2}\sqrt{\log(n)} } - 1\Bigg| = 0.
  \end{align*}
\end{lem}
\begin{proof}
  We prove the result for the first integral which is the most complicated. The
  proof for the second integral is similar. By a simple change of variable
  \begin{align*}
    \int_0^1 L' x^{-\alpha/2}\sqrt{\log(e/x)} \frac{1 - (1-x)^n}{x}\intd x%
    &= L' n^{\alpha/2}\sqrt{\log(en)}\int_0^n \sqrt{\frac{\log(en/y)}{\log(en)}} \frac{1 - (1-y/n)^n}{y^{1+\alpha/2}}\intd y.
  \end{align*}
  We decompose the integral in the rhs of the last display as
  $I_1 + I_2 + I_3 + I_4$, where
  \begin{align*}
    I_1
    &\coloneqq%
      \int_0^{\infty} \frac{1 - e^{-y}}{y^{1+\alpha/2}}\intd y\\
    I_2
    &\coloneqq%
      - \int_n^{\infty}\frac{1 - e^{-y}}{y^{1+\alpha/2}}\intd y\\
    I_3
    &\coloneqq%
      \int_0^n\Big[ \frac{1 - (1-y/n)^n}{y^{1+\alpha/2}} - \frac{1 - e^{-y}}{y^{1+\alpha/2}} \Big] \intd y\\
    I_4%
    &\coloneqq%
      \int_0^n \Big[\sqrt{\frac{\log(en/y)}{\log(en)}} - 1\Big] \frac{1 - (1-y/n)^n}{y^{1+\alpha/2}}\intd y.
  \end{align*}
  By direct computations it is found that $I_1 = 2\Gamma(1-\alpha/2)/\alpha$.
  Then,
  \begin{align*}
    |I_2|%
    &\leq \int_n^{\infty} \frac{1}{y^{1+\alpha/2}}\intd y%
      = \frac{2n^{-\alpha/2}}{\alpha}.
  \end{align*}
  Regarding $I_3$, we can rewrite,
  \begin{align*}
    I_3%
    &= \int_0^n \frac{e^{-y}[1 - e^{n\log(1-y/n) + y}]}{y^{1+\alpha/2}}\intd y
  \end{align*}
  Since $n\log(1-y/n)\leq y$ for all $0 < y < n$, we see that $I_3 \geq 0$.
  Furthermore,
  $n\log(1-y/n) + y \geq \frac{-y}{1-y/n} + y = \frac{-y^2/n}{1-y/n}$. This
  establishes that
  \begin{align*}
    |I_3|%
    &\leq \frac{1}{n}\int_0^{n/2} \frac{y^{1-\alpha/2}e^{-y}}{1 - y/n}\intd y%
      + \int_{n/2}^n \frac{e^{-y}}{y^{1+\alpha/2}}\intd y\\
    &\leq \frac{2}{n}\int_0^{n/2} y^{1-\alpha/2}e^{-y} \intd y%
      + e^{-n/2}\int_{n/2}^{\infty} \frac{1}{y^{1+\alpha/2}}\intd y\\
    &\leq \frac{2\Gamma(2-\alpha/2)}{n} + \frac{2^{\alpha/2+1} e^{-n/2}}{\alpha n^{\alpha/2}}.
  \end{align*}
  Now regarding $I_4$,
  \begin{align*}
    I_4%
    &= \frac{1}{\sqrt{\log(en)}}\int_0^n \Big[\sqrt{\log(en) - \log(y)} - \sqrt{\log(en)} \Big]%
      \frac{1 - (1-y/n)^n}{y^{1+\alpha/2}}\intd y\\
    &= \frac{1}{\sqrt{\log(en)}}\int_0^n \frac{-\log(y)}{\sqrt{\log(en/y)} + \sqrt{\log(en)}}%
      \frac{1 - (1-y/n)^n}{y^{1+\alpha/2}}\intd y.
  \end{align*}
  We deduce that,
  \begin{align*}
    |I_4|%
    &\leq \frac{1}{\log(en)}\int_0^1 \log\Big(\frac{1}{y}\Big) \frac{1 - (1-y/n)^n}{y^{1+\alpha/2}}\intd y%
      + \frac{1}{\log(en)}\int_1^n \frac{\log(y)}{y^{1+\alpha/2}}\intd y\\
    &\leq \frac{1}{(1 - 1/n)\log(en)}\int_0^1 \log\Big(\frac{1}{y}\Big)y^{-\alpha/2}\intd y%
      + \frac{1}{\log(en)}\int_1^n \frac{\log(y)}{y^{1+\alpha/2}}\intd y,
  \end{align*}
  where the second line follows because
  $1 - (1-y/n)^n = 1 - \exp[n\log(1-y/n)] \leq 1 - \exp[-y/(1-y/n)] \leq y/[1 - y/n]$
  for all $y \in [0,1]$. Finally,
  \begin{align*}
    |I_4|%
    &\leq \frac{1}{\log(en)} \frac{1}{1-1/n}\frac{1}{(1-\alpha/2)^2}%
      + \frac{1}{\log(en)} \frac{4}{\alpha^2}.
  \end{align*}
  Finally, the conclusion follows by gathering all the previous estimates.
\end{proof}

\begin{lem}
  \label{lem:var-cumulfreq}
  For all $p$, all $n\geq 1$ and all $1 \leq k \leq n$,
  \begin{equation*}
    \overline{\var}_p(C_{n,k})%
    \leq k \cdot \bar{\EE}_p[M_{n,k}].
  \end{equation*}
\end{lem}
\begin{proof}
  We bound the variance of $C_{n,k}$ using an Efron-Stein argument. By the
  Theorem~3.1 in \cite{boucheron:lugosi:massart:2013}, it is the case that
  \begin{equation*}
    \overline{\var}_p(C_{n,k})%
    \leq \sum_{i=1}^n\inf_{Z_i}\bar{\EE}_p[(C_{n,k} - Z_i)^2]
  \end{equation*}
  where the infimum is over all $(X_1,\dots,X_{i-1},X_{i+1},X_n)$-measurable
  random variables. In particular, defining
  $Y_{n,j}^{(i)} \coloneqq Y_{n,j} - \1_{X_i = j}$, and
  $C_{n,k}^{(i)} \coloneqq \sum_{j\geq 1}^n\1_{Y_{n,j}^{(i)}\geq k}$, we see
  that $C_{n,k}^{(i)}$ does not depend on $X_i$. Hence the last display
  entails that
  \begin{equation*}
    \overline{\var}_p(C_{n,k})%
    \leq \sum_{i=1}^n\bar{\EE}_p[(C_{n,k} - C_{n,k}^{(i)})^2].
  \end{equation*}
  But,
  \begin{align*}
    C_{n,k} - C_{n,k}^{(i)}%
    &= \sum_{j\geq 1}(\1_{Y_{n,j}\geq k} - \1_{Y_{n,j}^{(i)}\geq k})\\
    &= \sum_{j\geq 1}\1_{Y_{n,j} = k}\1_{X_i = j},
  \end{align*}
  where the second line follows because $Y_{n,j} \geq Y_{n,j}^{(i)} \geq Y_{n,j} - 1$ almost-surely, hence
  $\1_{Y_{n,j}\geq k} - \1_{Y_{n,j}^{(i)}\geq k}$ is either zero or one, and it is one
  iff $Y_{n,j}= k$ and $Y_{n,j}^{(i)} = k-1$; that is iff $Y_{n,j}=k$ and
  $X_i = j$. It follows,
  \begin{align*}
    \overline{\var}_{p}(C_{n,k})
    &\leq \sum_{i=1}^n\bar{\EE}_{p}\Big[ \sum_{j\geq 1}\1_{Y_{n,j}=k}\1_{X_i=j}\sum_{j'\geq 1}\1_{Y_{n,j'}=k}\1_{X_i= j'}\Big]\\
    &= \sum_{j\geq 1}\bar{\EE}_{p}\Big[\1_{Y_{n,j}=k}\sum_{i=1}^n\1_{X_i=j} \Big]\\
    &= \sum_{j\geq 1}\bar{\EE}_{p}[\1_{Y_{n,j}=k}Y_{n,j}]\\
    &= k \bar{\EE}_{p}\Big[\sum_{j\geq 1}\1_{Y_{n,j}=k} \Big]\\
    &= k\cdot \EE_{p}[M_{n,k}].
  \end{align*}
\end{proof}

\begin{prp}
  \label{pro:7}
  The following is true for all $\alpha \in (0,1)$ and all $p$
  \begin{equation*}
    \overline{\var}_{p}(\varphi_n(\alpha))%
    \leq \max\Big\{1, \Big( \frac{\alpha}{1 - \alpha} \Big)^2\Big\}\bar{\EE}_p[C_{n,1}].
  \end{equation*}
\end{prp}
\begin{proof}
  We bound the variance of $\varphi_n(\alpha')$ using the same Efron-Stein
  argument as in the proof of Lemma~\ref{lem:var-cumulfreq}. For simplicity we
  write $\varphi \equiv \varphi_n(\alpha)$; and we let $C_{n,k}^{(i)}$ having the
  same definition as in Lemma~\ref{lem:var-cumulfreq}. Then defining
  $\varphi_i \coloneqq \sum_{k=1}^{n-1}\frac{\alpha}{k-\alpha}C_{n,k+1}^{(i)} - C_{n,1}^{(i)}$,
  following the arguments in Lemma~\ref{lem:var-cumulfreq}%
  \begin{align*}
    \overline{\var}_{p}(\varphi)%
    &\leq%
      \sum_{i=1}^n\EE_{p}[(\varphi - \varphi_i)^2]
  \end{align*}
  and
  \begin{align*}
    (\varphi - \varphi_i)^2%
    &= \Big(\sum_{k=1}^{n-1} \frac{\alpha'}{k-\alpha'}(C_{n,k+1} - C_{n,k+1}^{(i)}) - (C_{n,1} - C_{n,1}^{(i)})  \Big)^2\\
    &= \Big\{ \sum_{j\geq 1}\1_{X_i = j}\Big(\sum_{k=1}^{n-1} \frac{\alpha'}{k -\alpha'}\1_{Y_{n,j}=k+1} - \1_{Y_{n,j}=1} \Big) \Big\}^2.
  \end{align*}
  It follows,
  \begin{align*}
    \sum_{i=1}^n (\varphi - \varphi_i)^2%
    &= \sum_{j\geq 1}\Big(\sum_{i=1}^n\1_{X_i=j}\Big)\Big(\sum_{k=1}^{n-1} \frac{\alpha'}{k -\alpha'}\1_{Y_{n,j}=k+1} - \1_{Y_{n,j}=1} \Big)^2\\
    &= \sum_{j\geq 1}Y_{n,j}\Big(\sum_{k=1}^{n-1} \frac{\alpha'}{k -\alpha'}\1_{Y_{n,j}=k+1} - \1_{Y_{n,j}=1} \Big)^2\\
    &= \sum_{j\geq 1}\1_{Y_{n,j} = 1} + \sum_{j\geq 1}Y_{n,j}\1_{Y_{n,j}\ne 1}\Big(\sum_{k=1}^{n-1} \frac{\alpha'}{k -\alpha'}\1_{Y_{n,j}=k+1}\Big)^2\\
    &= M_{n,1}%
      + \sum_{k=1}^{n-1}\sum_{k'=1}^{n-1} \frac{\alpha'}{k-\alpha'}\frac{\alpha'}{k'-\alpha'} \sum_{j\geq 1}Y_{n,j}\1_{Y_{n,j}\ne 1}\1_{Y_{n,j}=k+1}\1_{Y_{n,j}=k'+1}\\
    &= M_{n,1} + \sum_{k=1}^{n-1} \Big( \frac{\alpha'}{k - \alpha'} \Big)^2 \sum_{j\geq 1}Y_{n,j}\1_{Y_{n,j}=k+1}.
  \end{align*}
  Hence we have shown,
  \begin{align*}
    \sum_{i=1}^n \bar{\EE}_{p}[(\varphi - \varphi_i)^2]%
    &= \bar{\EE}_{p}[M_{n,1}] + \sum_{k=1}^{n-1} \Big( \frac{\alpha'}{k - \alpha'} \Big)^2 (k+1)\bar{\EE}_{p}[M_{n,k+1}]\\
    &\leq \max\Big\{1,2\Big(\frac{\alpha'}{1 - \alpha'}\Big)^2\Big\}\bar{\EE}_p[C_{n,1}],
  \end{align*}
  where the second line follows because $(k+1)/(k-\alpha')^2 \leq 2/(1-\alpha')^2$ for all
  $k \geq 1$ and because $\sum_{\ell=1}^nM_{n,\ell} =  C_{n,1}$.
\end{proof}

\subsection{Proof of Theorem~\ref{thm:6}}
\label{sec:proof-theor-refthm:6}

We want to bound the probability of the event
\begin{align*}
  E_n \coloneqq%
  \Set*{ \ell_{\alpha}(\hat{U}_{n,\lambda n},U_{n,\lambda n}) \leq D\min\Big(1, \frac{\log(n)\log(\lambda)^2}{n^{\alpha}} \Big) }
\end{align*}
for $(\alpha,L,L') \in \mathbb{A}$ and $p \in \mathcal{P}(\alpha,L,L')$ be
arbitrary. First consider the case where
$\log(\lambda) > \frac{2C n^{\alpha/2}}{\sqrt{\log(n)}}$. Then we see that
$\log(\lambda) > \frac{C n^{\hat{\alpha}_n/2}}{\sqrt{\log(n)}}\cdot 2n^{(\alpha - \hat{\alpha}_n)/2}$.
We deduce that on the event that
$|\hat{\alpha}_n - \alpha| \leq \frac{2\log(2)}{\log(n)}$ is must be that
$\hat{U}_{n,m} = 0$, that is
$\ell_{\alpha}(\hat{U}_{n,\lambda n},U_{n,\lambda n}) = \frac{U_{n,\lambda n}^2}{(\lambda n)^{2\alpha}} \leq \frac{K_{n+\lambda n}^2}{(\lambda n)^{2\alpha}}$.
It follows for $D' = D \min(1, 4C^2)$,
\begin{align*}
  \bar{\PP}_p(E_n^c)%
  &\leq%
    \bar{\PP}_p(K_{n+\lambda n} \geq (\lambda n)^{\alpha} \sqrt{D'} )%
    + \bar{\PP}_p\Big(| \hat{\alpha}_n - \alpha| > \frac{2\log(2)}{\log(n)}\Big)\\
  &\leq \frac{\overline{\var}_p(K_{n+\lambda n})}{(\lambda n)^{2\alpha} D'}%
    + \bar{\PP}_p\Big(| \hat{\alpha}_n - \alpha| > \frac{2\log(2)}{\log(n)}\Big).
\end{align*}
The first term of the rhs of the last display can be made smaller than
$\varepsilon/2$ by choosing $D$ (and thus $D'$) large enough, in virtue of
Lemma~\ref{lem:var-cumulfreq} and~\ref{lem:unif-expect-cumulfreq} [remarking
that $K_n = C_{n,1}$]. The second term of the rhs of the last display is
controlled via the Theorem~\ref{thm:3} and will be smaller than $\varepsilon/2$
whenever $n$ gets large enough.

Now consider the case where
$\log(\lambda) \leq 2C n^{\alpha/2}/ \sqrt{\log(n)}$. We define the event
$A_n \coloneqq \Set{\log(\lambda) \leq C n^{\hat{\alpha}_n/2}/ \sqrt{\log(n)}}$,
such that $\hat{U}_{n,m} = 0$ on the complement $A_n^c$. With the exact same
argument as above
\begin{align*}
  \bar{\PP}_p\Big(E_n^c \bigcap A_n^c\Big)%
  &\leq \bar{\PP}_p(K_{n+\lambda n} \geq (\lambda n)^{\alpha} \sqrt{D'} )
  \leq \frac{\varepsilon}{3}
\end{align*}
by choosing $D$ sufficiently large. While if we are on the event $A_n$, it is
the case that $\hat{U}_{n,\lambda n} = K_n((1+\lambda)^{\hat{\alpha}_n} - 1)$.
Equivalently, it is the case that
$\hat{U}_{n,\lambda n} - U_{n,\lambda n} = K_n(1 + \lambda)^{\hat{\alpha}_n} - K_{n+\lambda n}$.
Define another event
$B_n \coloneqq \Set{|\hat{\alpha}_n - \alpha|^2 \leq B\log(n)n^{-\alpha} }$ for
a large enough constant $B > 0$. By choosing $B$ large enough, by the
Theorem~\ref{thm:3} it will be that $\bar{\PP}_p(B_n^c) \leq \varepsilon/3$. So it is enough to analyze what happens on the
event $A_n \cap B_n$. But on $A_n \cap B_n$ for $n$ large enough,
\begin{align*}
  |\hat{U}_{n,\lambda n} - U_{n,\lambda n}|%
  &= |K_n(1 + \lambda)^{\alpha}[(1 + \lambda)^{\hat{\alpha}_n - \alpha} - 1  ]%
    + [K_n(1 + \lambda)^{\alpha} - K_{n+\lambda n}]|\\
  &\leq 2e^{2C\sqrt{B}}\cdot K_n(1 + \lambda)^{\alpha}|\hat{\alpha}_n - \alpha|\log(1+\lambda) %
    + |K_n(1 + \lambda)^{\alpha} - K_{n+\lambda n}|\\
  &\lesssim \frac{K_n \lambda^{\alpha}\log(\lambda)\sqrt{\log(n)} }{n^{\alpha/2}}%
    + |\bar{\EE}_p[K_n](1 + \lambda)^{\alpha} - \bar{\EE}_p[K_{n+\lambda n}] |\\%
  &\qquad%
    + \lambda^{\alpha}|K_n - \bar{\EE}_p[K_n]|%
    + |K_{n+\lambda n} - \bar{\EE}_p[K_{n+\lambda n}]|.
\end{align*}
Again remarking that $K_n = C_{n,1}$, we deduce from
Lemma~\ref{lem:unif-expect-cumulfreq} and a couple applications of Chebychev's
inequality that for $D$ large enough,
\begin{align*}
  \bar{\PP}_p(E_n^c \cap A_n \cap B_n)%
  &\leq \frac{\varepsilon}{3}.
\end{align*}
Combining all the estimates above, we have demonstrated that in any cases
$\bar{\PP}_p(E_n^c) \leq \varepsilon$ by choosing $D$ sufficiently large.
Inspection of the proof shows that this last estimate holds uniformly over
$(\alpha,L,L') \in \mathbb{A}$ and $p \in \mathcal{P}(\alpha,L,L')$.

\subsection{Proof of Proposition~\ref{pro:mle-weaker-bias}}
\label{sec:proof-prop-weaker-bias}

The proof of the Proposition is nearly identical to the proof of the
Theorem~\ref{thm:3}. In particular, the Proposition~\ref{pro:6} is unchanged, and to
finish the proof it is enough to adapt the proof of the
Proposition~\ref{pro:control-bias-klmin} to establish that
\begin{equation*}
  \limsup_{n\to\infty} n^{\beta}|\alpha_n^{*}(p) - \alpha| \leq C
\end{equation*}
for some $C > 0$ if $\bar{F}_p(x) = Lx^{-\alpha} + O(x^{-\beta})$ as $x \to 0$.
This can be established by adapting the proofs of
Proposition~\ref{pro:control-bias-klmin} and Lemma~\ref{lem:5} line by line.

\section{Minimax lower bound for estimating the tail-index}
\label{sec:minimax-lower-bound}

This appendix contains the proof of Theorem~\ref{thm:5}.

\subsection{Guidelines and construction of the hypotheses}
\label{sec:guid-constr-hypoth}

The lower bound is obtained from the result of Theorem~\ref{thm:2}. For
$0 < \alpha_0 < \alpha_1 < 1$ to be chosen accordingly, we consider the priors
$\nu_0$ an $\nu_1$ consisting on PK distributions with respective L\'evy measures
\begin{equation*}
  \rho_0(\intd x)%
  = \frac{\alpha_0 x^{-1-\alpha_0}\intd x}{\Gamma(1-\alpha_0)},\qquad%
  \rho_1(\intd x)%
  = \frac{\alpha_0 x^{-1-\alpha_1}\,\intd x}{\Gamma(1-\alpha_0)}.
\end{equation*}
Then, we wish to apply Theorem~\ref{thm:2} to the functional
$\delta(\Pi,p) \equiv \delta(p) = \lim_{x\to 0}\frac{\log \bar{F}_{p}(x)}{-\log(x)}$.
It is a well-known fact (\citet{Kin(75),Gne(07)}) that
$\nu_0(\delta(P) = \alpha_0) = 1$ and $\nu_1(\delta(P) = \alpha_1) = 1$. Hence,
taking $2\varepsilon = \alpha_1 - \alpha_0$, we can take $\kappa = 0$ and it
suffices to show that for suitable choice of constants we have
\begin{equation}
  \label{eq:36}
  \min_{j=0,1}\nu_j\big( P\in \mathcal{P}^{*} \big) \geq 1 - \xi,
\end{equation}
where we let,
\begin{equation*}
  \mathcal{P}^{*} \coloneqq%
  \bigcup_{(\alpha,L,L')\in \mathbb{A}}\mathcal{P}(\alpha,L,L'),
\end{equation*}
and also,
\begin{equation}
  \label{eq:38}
  \KL(\PP_{\rho_0}^n;\PP_{\rho_1}^n) \leq \frac{1}{2}.
\end{equation}
Combining the results of the two next sections shows that this will be the case
if $\overline{L}$ is sufficiently large, if $\alpha_0 = \underline{\alpha}$, and if for
$c > 0$ small enough but constant we take
\begin{align*}
  \varepsilon%
  = \frac{c}{n^{\underline{\alpha}/2}\log(n)}.
\end{align*}
Summarizing, \eqref{eq:36} and \eqref{eq:38} combined with the
Theorem~\ref{thm:2} guarantee that whenever $n$ is large enough
\begin{align*}
  \inf_{\hat{\alpha}_n}\sup_{\substack{(\alpha,L,L')\in \mathbb{A}\\p \in \mathcal{P}(\alpha,L,L')}}\bar{\PP}_p\Big( |\hat{\alpha}_n - \alpha| \geq \varepsilon \Big) \geq \frac{1}{4}.
\end{align*}
The final result follows because
$\varepsilon \geq \frac{c}{n^{\alpha/2}\log(n)}$ for all $(\alpha,L,L')\in \mathbb{A}$.

\subsection{Bounding the KL divergence between hypotheses}
\label{sec:bound-kl-diverg}

In this section we establish that for the choice of hypotheses and constants
made in Section~\ref{sec:guid-constr-hypoth}, it is true that \eqref{eq:38}
holds.

\begin{prp}
  \label{pro:8}
  There is a constant $C > 0$ depending only on $\alpha_0$ such that,
  \begin{align*}
    \KL(\PP_{\rho_0}^n;\PP_{\rho_1}^n)%
    &\leq C\varepsilon^2 n^{\alpha_0}\log(n)^2.
  \end{align*}
\end{prp}
\begin{proof}
  By direct computations, we obtain that
  $\mu_0(u,\ell) = \frac{\alpha_0\Gamma(\ell-\alpha_0)}{\Gamma(1-\alpha_0)}u^{-\ell + \alpha_0}$ and
  $\mu_1(u,\ell) = \frac{\alpha_0\Gamma(\ell-\alpha_1)}{\Gamma(1-\alpha_0)}u^{-\ell + \alpha_1}$. It follows,
  \begin{align*}
    -\log \frac{\mu_1(u,\ell)}{\mu_0(u,\ell)}%
    + \frac{\mu_1(u,\ell)}{\mu_0(u,\ell)} - 1
    &=%
      - \log \frac{\Gamma(\ell - \alpha_1)}{\Gamma(\ell - \alpha_0)}%
      + \frac{\Gamma(\ell - \alpha_1)}{\Gamma(\ell - \alpha_0)}%
      - 1\\%
    &\quad+ \frac{\Gamma(\ell - \alpha_1)}{\Gamma(\ell- \alpha_0)}\Big(u^{\alpha_1 - \alpha_0} - (\alpha_1 - \alpha_0)\log u - 1\Big)\\%
    &\quad+ \Big\{\frac{\Gamma(\ell - \alpha_1)}{\Gamma(\ell- \alpha_0)} -1 \Big\} (\alpha_1 - \alpha_0)\log u.
  \end{align*}
  We let $\psi$ denote the digamma function; that is the derivative of the map
  $z\mapsto \log \Gamma(z)$. Then, by a Taylor expansion argument
  \begin{align*}
    \sup_{\ell=1,\dots,n}\Big|%
    \frac{\Gamma(\ell-\alpha_1)}{\Gamma(\ell-\alpha_0)} - 1%
    \Big|%
    &\leq \sup_{\ell=1,\dots,n}\sup_{\alpha\in(\alpha_0,\alpha_1)}\Big|e^{\psi(\ell - \alpha)(\alpha_0 - \alpha_1)} - 1 \Big|\\
    &\leq (1 + o(1))\log(n)|\alpha_0 - \alpha_1|,
  \end{align*}
  where the second line follows because $\psi(z) \sim \log(z)$ as $z \to \infty$ and because
  by construction $\log(n)(\alpha_0 - \alpha_1) = o(1)$. Similarly,
  \begin{align*}
    0%
    &\leq  u^{\alpha_1 - \alpha_0} - (\alpha_1 - \alpha_0)\log u - 1
    \leq (\alpha_1 - \alpha_0)^2\log(u)^2\cdot \max(1,\, u^{\alpha_1 - \alpha_0}).
  \end{align*}
  Therefore,
  \begin{align*}
    0%
    &\leq -\log \frac{\mu_1(u,\ell)}{\mu_0(u,\ell)}%
      + \frac{\mu_1(u,\ell)}{\mu_0(u,\ell)} - 1\\
    &\leq \big(1 + o(1)\big)\big(\alpha_0 - \alpha_1\big)^2%
      \Big(\frac{\log(n)^2}{2}%
      + \log(u)^2\cdot\max(1,u^{2\varepsilon})%
      + \log(n)\cdot|\log(u)|%
      \Big)
  \end{align*}
  Thus by Proposition~\ref{thm:4}, and by the fact that
  $\sum_{\ell=1}^nM_{n,\ell} = K_n$,
  \begin{multline*}
    \KL(\PP_{\rho_0}^n;\PP_{\rho_1}^n)%
    \leq 4\varepsilon^2\big(1 + o(1)\big)%
    \Big(\frac{\EE_{\rho_0}[K_n]\log(n)^2}{2}\\%
    + \EE_{\rho_0}[K_n\log(U_n)^2\cdot\max(1,U_n^{2\varepsilon})]%
    + \log(n)\cdot\EE_{\rho_0}[K_n|\log(U_n)|]%
    \Big).
  \end{multline*}
  It is well known that under $\rho_0$ it holds $\EE_{\rho_0}[K_n] \asymp n^{\alpha_0}$; see for
  instance \citet{Pit(06)}. Furthermore,
  \begin{align*}
    \EE_{\rho_0}[K_n\log(U_n)^2\cdot\max(1,U_n^{2\varepsilon})]%
    &\leq \EE_{\rho_0}[K_n \log(U_n)^2] + \EE_{\rho_0}[K_n \cdot U_n^{2\varepsilon}\log(U_n)^2],
  \end{align*}
  and by Cauchy-Schwarz' inequality
  \begin{align*}
    \EE_{\rho_0}[K_n \cdot U_n^{2\varepsilon}\log(U_n)^2]%
    &\leq \EE_{\rho_0}[K_n^2]^{1/2}\EE_{\rho_0}[U_n^{4\varepsilon}\log(U_n)^4]^{1/2}.
  \end{align*}
  Recall that $U_n \equalDist G_n/T$ where $G_n \sim \gammaDist(n,1)$ is
  independent of $T$, so
  \begin{align*}
    \EE_{\rho_0}[U_n^{4\varepsilon}\log(U_n)^4]%
    &= \EE_{\rho_0}[G_n^{4\varepsilon} T^{-4\varepsilon}\log(G_n/T)^4]\\
    &\leq 8\EE_{\rho_0}[G_n^{4\varepsilon} T^{-4\varepsilon}\{\log(G_n)^4 + \log(T)^4 \}]\\
    &= 8 \EE_{\rho_0}[G_n^{4\varepsilon}\log(G_n)^4]\EE_{\rho_0}\Big[\frac{1}{T^{4\varepsilon}}\Big]%
      + 8\EE_{\rho_0}[G_n^{4\varepsilon}]\EE_{\rho_0}\Big[\frac{\log(T)^4}{T^{4\varepsilon}}\Big].
  \end{align*}
  Under $\rho_0$, the random variable $T$ has a $\alpha_0$-stable distribution
  (\citet{Pit(06)}). By a classical analysis of the stable distributions
  (\citet{Sam(94)}), $T^{-4\varepsilon}$ and $T^{-4\varepsilon}\log(T)^4$ have finite expectation, and
  they can be bounded by a constant depending only on $\alpha_0$ as $\varepsilon = o(1)$.
  Since, $\varepsilon = o(n^{-\alpha_0/2})$ we obtain
  \begin{align*}
    \EE_{\rho_0}[U_n^{4\varepsilon}\log(U_n)^4]%
    &\lesssim \log(n)^4.
  \end{align*}
  Finally, $\EE_{\rho_0}[K_n^2] =O(n^{2\alpha_0})$ (\citet{Pit(06)}), so
  \begin{align*}
    \EE_{\rho_0}[K_n \log(U_n)^2\cdot \max(1,U_n^{2\varepsilon})]%
    \lesssim n^{\alpha_0}\log(n)^2.
  \end{align*}
  The remaining term $\EE_{\rho_0}[K_n|\log(U_n)|]$ can be bounded using the same
  arguments, giving the conclusion.
\end{proof}

\subsection{Priors generate regular probability vectors with large
  probabibility}
\label{sec:priors-gener-regul}

In this section we establish that for the choice of hypotheses and constants
made in Section~\ref{sec:guid-constr-hypoth}, it is true that equation~\eqref{eq:36} holds. We
proceed by proving that with large $\nu_j$-probability, there exists $L,L'$ such that
$p \in \mathcal{P}(\alpha_j,L,L')$. First, observe that by direct computations, we
have,%
\begin{align}
  \label{eq:rhobar}
  \bar{\rho}_j(x)%
  &= \frac{\alpha_0/\alpha_j}{\Gamma(1 - \alpha_0)}x^{-\alpha_j}.
\end{align}
We define two random variables
$V_j \coloneqq \lim_{x\to 0} x^{\alpha_j}\bar{\rho}_j(x T)$. Then, the
Proposition~\ref{pro:10} below guarantees that for all $\epsilon > 0$ there exists $L'$
such that $\min_{j=1,2}\nu_j\{ P \in \mathcal{P}(\alpha_j,V_j,L') \} \geq 1 - \epsilon$.

\begin{prp}
  \label{pro:10}
  For all $j=0,1$ and for all $\epsilon > 0$ there exists a constant $L' > 0$ depending
  only $\epsilon$ and $\alpha_j$ such that with $\nu_j$-probability at least
  $1-\epsilon$%
  \begin{equation*}
    \sup_{x\in (0,1)} \frac{|x^{\alpha_j} \bar{F}_P(x) - V_j |}{\sqrt{x^{\alpha_j}\log(e/x)}}%
    \leq L'.
  \end{equation*}
\end{prp}
\begin{proof}
  Let $j=0$ or $j=1$. We remark that $V_j = x^{\alpha_j}\bar{\rho}_j(xT)$ by
  equation~\eqref{eq:rhobar}. Also, recall that
  $\bar{F}_{P}(x) \equalDist \sum_{i\geq 1}\1_{J_i > x T}$ conditional on
  $(J_i)_{i\geq 1}$. Then, in distribution
  \begin{align*}
    \sup_{x\in (0,1)} \frac{|x^{\alpha_j} \bar{F}_P(x) - V_{j}|}{\sqrt{x^{\alpha_j}\log(e/x)}}%
    &\leq \sup_{x\in (0,1)} \sqrt{\frac{x^{\alpha_j}}{\log(e/x)}} \Big|\sum_{i\geq 1}\1_{J_i> x T} - \bar{\rho}_j(xT)  \Big|\\
    &= \sup_{x\in (0,T)} \sqrt{\frac{(x/T)^{\alpha_j}}{\log(eT/x)}} \Big|\sum_{i\geq 1}\1_{J_i> x} - \bar{\rho}_j(x)  \Big|.
  \end{align*}
  Let
  $\bar{\rho}_j^{-1}(u) \coloneqq \inf\Set{x > 0 \given \bar{\rho}_j(x) < u}$
  and let $(\Gamma_i)_{i\geq 1}$ be a unit rate Poisson process on $\NNReals$.
  By (\citet{Fer(72)}), it is the case that
  $(J_i)_{i\geq 1} \equalDist (\bar{\rho}_j^{-1}(\Gamma_i))_{i\geq 1}$.
  Furthermore, $\bar{\rho}_j^{-1}(x) > y \iff x < \bar{\rho}_j(y)$. Hence, in
  distribution,
  \begin{align*}
    \sup_{x\in (0,1)} \frac{|x^{\alpha_j} \bar{F}_P(x) - V_j |}{\sqrt{x^{\alpha_j}\log(e/x)}}%
    &= \sup_{x\in (0,T)} \sqrt{\frac{(x/T)^{\alpha_j}}{\log(eT/x)}}\Big|\sum_{i\geq 1}\1_{\bar{\rho}^{-1}(\Gamma_i) > x} - \bar{\rho_j}(x)  \Big|\\
    &= \sup_{x\in (0,T)} \sqrt{\frac{(x/T)^{\alpha_j}}{\log(e T/x)}}\Big|\sum_{i\geq 1}\1_{\Gamma_i < \bar{\rho}_j(x)} - \bar{\rho}_j(x)  \Big|\\
    &=\sup_{y\in (\bar{\rho}_j(T),\infty )}\sqrt{\frac{(\bar{\rho}_j^{-1}(y)/T)^{\alpha_j}}{\log(eT/\bar{\rho}_j^{-1}(y))}}\Big|\sum_{i\geq 1}\1_{\Gamma_i < y} - y  \Big|.
  \end{align*}
  We deduce from equation~\eqref{eq:rhobar} that
  $\bar{\rho}_j^{-1}(y) = C_1(\alpha_j)y^{-1/\alpha_j}$, where
  $C_1(\alpha_j) \coloneqq [(\alpha_0/\alpha_j)/\Gamma(1-\alpha_0)]^{1/\alpha_j}$. Also
  $T = \bar{\rho}_j^{-1}(\bar{\rho}(T))$; ie
  $\bar{\rho}_j^{-1}(y)/T = (y/\bar{\rho}_j(T))^{-1/\alpha_j}$. Since the function
  $z \mapsto z^{\alpha_j}/\log(e/z)$ is monotone increasing in $(0,\infty)$, we have on the
  event that $\bar{\rho}_j(T) \leq A$
  \begin{align*}
    \sup_{x\in (0,1)} \frac{|x^{\alpha_j} \bar{F}_P(x) - V_j |}{\sqrt{x^{\alpha_j}\log(e/x)}}%
    &\leq \sup_{y \in (0,\infty)} \sqrt{ \frac{A/y}{\log[e(y/A)^{1/\alpha_j} ]} }%
      \Big|\sum_{i\geq 1}\1_{\Gamma_i < y} - y  \Big|\\
    &\lesssim \sup_{y\in (0,\infty)} \frac{\big|\sum_{j\geq 1}\1_{\Gamma_j < y} - y  \big|}{\sqrt{\max\{1,y\log y \}} }
  \end{align*}
  The random variable $T$ has a $\alpha_j$ stable distribution, so it is always
  possible to choose $A$ large enough so that the event $\bar{\rho}_j(T) \leq A$ occur
  with probability $1 - \epsilon/2$. Furthermore, by the Lemma~\ref{lem:3}, for $B$
  sufficiently large the rhs of the last display will be smaller than $B$
  with probability at least $1-\epsilon/2$. This concludes the proof.
\end{proof}

\begin{lem}
  \label{lem:3}
  Let $(\Gamma_i)_{i\geq 1}$ be an homogeneous Poisson process on the half real
  line. Then, for every $\epsilon > 0$ there exists $D > 0$ such that with probability
  at least $1-\epsilon$
  \begin{align*}
    \sup_{y\in (0,\infty)} \frac{\big|\sum_{j\geq 1}\1_{\Gamma_j < y} - y  \big|}{\sqrt{\max\{1,y\log y \}} } \leq D.
  \end{align*}
\end{lem}
\begin{proof}
  Let $y_0$ be the unique non-negative solution of $y \log y = 1$. Obviously,
  \begin{align*}
    \sup_{y\in (0,\infty)} \frac{\big|\sum_{j\geq 1}\1_{\Gamma_j < y} - y  \big|}{\sqrt{\max\{1,y\log y \}} }%
    &\leq \max\Big\{ \sup_{y\in (0,y_0)}\Big|\sum_{j\geq 1}\1_{\Gamma_j \leq y} - y\Big|,\, \sup_{y\geq y_0}\frac{\big|\sum_{j\geq 1}\1_{\Gamma_j < y} - y  \big|}{\sqrt{y \log y} } \Big\}.
  \end{align*}
  By Doob's maximal inequality,
  \begin{align*}
    \EE\Big[\sup_{y\in (0,y_0)}\Big|\sum_{j\geq 1}\1_{\Gamma_j \leq y} - y\Big|^2  \Big]%
    \leq 4 \EE\Big[\Big|\sum_{j\geq 1}\1_{\Gamma_j \leq y_0} - y_0\Big|^2\Big]%
    = 4y_0.
  \end{align*}
  Hence, by Chebychev's inequality we have that
  $\sup_{y\in (0,y_0)}|\sum_{j\geq 1}\1_{\Gamma_j< y} - y| \leq C$ with probability at least
  $1 - 4y_0/C^2$. For $C > 0$ large enough this will happen with probability $1-\epsilon/2$.
  For the other term, we decompose,
  \begin{align*}
    \sup_{y\geq y_0}\frac{\big|\sum_{j\geq 1}\1_{\Gamma_j < y} - y  \big|}{\sqrt{y \log y} }%
    &= \sup_{k\geq 1}\sup_{y\in [ky_0,(k+1)y_0]} \frac{\big|\sum_{j\geq 1}\1_{\Gamma_j < y} - y  \big|}{\sqrt{y \log y} }\\
    &\leq \sup_{k\geq 1} \frac{1}{\sqrt{ky_0\log(ky_0)}} \sup_{y\in [ky_0,(k+1)y_0]}\big|\sum_{j\geq 1}\1_{\Gamma_j < y} - y  \big|.
  \end{align*}
  Note that for all $y \in [ky_0,(k+1)y_0]$ we have,
  \begin{align*}
    \sum_{j\geq 1}\1_{\Gamma_j < y} - y%
    &\leq \sum_{j\geq 1}\1_{\Gamma_j < (k+1)y_0} - ky_0%
      = \sum_{j\geq 1}\1_{\Gamma_j < (k+1)y_0} - (k + 1)y_0 + y_0,
  \end{align*}
  and,
  \begin{align*}
    \sum_{j\geq 1}\1_{\Gamma_j < y} - y%
    &\geq \sum_{j\geq 1}\1_{\Gamma_j < k} -(k+1)y_0%
      = \sum_{j\geq 1}\1_{\Gamma_j < k} - ky_0 - y_0.
  \end{align*}
  Therefore,
  \begin{align*}
    \sup_{y\in [ky_0,(k+1)y_0]} \big|\sum_{j\geq 1}\1_{\Gamma_j < y} - y  \big|%
    \leq y_0 + \max_{\ell=k,k+1}\Big\{\Big|\sum_{j\geq 1}\1_{\Gamma_j < \ell} - \ell \Big| \Big\}.
  \end{align*}
  It follows, because $(k+1)y_0\log((k+1)y_0) \geq k y_0\log(ky_0)$,
  \begin{align*}
    \sup_{y\geq y_0}\frac{\big|\sum_{j\geq 1}\1_{\Gamma_j < y} - y  \big|}{\sqrt{y \log y} }%
    &\leq \sqrt{\frac{y_0}{\log(y_0)}}%
      + \sup_{k\geq 1} \frac{\big|\sum_{j\geq 1}\1_{\Gamma_j < k} - k \big|}{\sqrt{ky_0 \log(ky_0)}}.
  \end{align*}
  We finish the proof using Chernoff's bound on the Poisson distribution.
  Indeed, for any $t > 0$
  \begin{align*}
    \PP\Big\{ \sup_{k\geq 1} \frac{\big|\sum_{j\geq 1}\1_{\Gamma_j < k} - k \big|}{\sqrt{ky_0 \log(ky_0)}} > t\Big\}%
    &\leq \sum_{k\geq 1}\PP\Big\{ \frac{\big|\sum_{j\geq 1}\1_{\Gamma_j < k} - k \big|}{\sqrt{ky_0 \log(ky_0)}} > t\Big\}\\
    &\leq 2\sum_{k\geq 1}\exp\Big\{- \frac{t^2 k y_0 \log(ky_0)/2}{k + t \sqrt{ky_0\log(ky_0)}} \Big\}\\
    &= 2\sum_{k\geq 1}\exp\Big\{- \frac{t^2 y_0 \log(ky_0)/2}{1 + t y_0 \sqrt{\log(ky_0)/(ky_0)}} \Big\}\\
    &\leq 2\sum_{k\geq 1}\exp\Big\{- \frac{t^2 y_0 \log(ky_0)/2}{1 + t y_0 e^{-1/2}} \Big\}.
  \end{align*}
  Letting  $g(t) \coloneqq  \frac{y_0t^2/2}{1 + y_0t e^{-1/2}}$, and remarking
  that $g(t) > 1$ if we take $t$ sufficiently large,
  \begin{align*}
    \PP\Big\{ \sup_{k\geq 1} \frac{\big|\sum_{j\geq 1}\1_{\Gamma_j < k} - k \big|}{\sqrt{ky_0 \log(ky_0)}} > t\Big\}%
    &\leq 2 \sum_{k\geq 1}(y_0k)^{-g(t)}\\
    &= 2y_0^{-g(t)}%
      + 2 \sum_{k\geq 2}(y_0k)^{-g(t)}\\
    &\leq 2y_0^{-g(t)} + 2y_0^{-g(t)}\int_1^{\infty} x^{-g(t)}\intd x\\
    &= \frac{2 g(t)}{g(t) - 1} y_0^{-g(t)}.
  \end{align*}
  Observe that $y_0 > 1$ (this is obvious since for $y \leq 1$ it must be that
  $y\log(y) \leq 0$, hence any positive solution of $y\log(y) = 1$ must strictly
  greater than $1$; by numerical computations, $y_0 \approx 1.76$). Then, the previous
  display can be made smaller than $\epsilon/2$ by taking $t$ sufficiently large, which
  concludes the proof of the lemma.
\end{proof}

\section{Minimax lower bound for estimating the unseen}
\label{sec:minimax-lower-bound-1}

This appendix contains the proof of Theorem~\ref{thm:7}. We prove the theorem by showing that predicting the unseen over the class
$\mathcal{P}(\alpha,L,L')$ is essentially a harder problem than estimating $\alpha$. That
is, we propose to reduce the problem to the estimation of $\alpha$. In the whole
proof, for $G_1,G_2 > 0$ to be chosen accordingly, we let
\begin{equation*}
  \varepsilon
  \coloneqq \min\Big\{ \frac{G_1 \log(1+\lambda)}{n^{\underline{\alpha}/2}\log(n)}  ,\, G_2\Big\}.
\end{equation*}

\subsection{Reduction to the problem of estimating the tail-index}
\label{sec:reduct-probl-estim}

The first step of the reduction consists on noticing that estimating $U_{n,m}$
and $K_{n+m}$ are equivalent problems. Indeed $U_{n,m} = K_{n+m} - K_n$ implies
that for all estimator $\hat{U}$ and all $p \in \mathcal{P}(\alpha,L,L')$, writing
$\hat{K} = K_n + \hat{U}$,
\begin{align*}
  \bar{\PP}_{p}\big( \ell_{\alpha}(\hat{U},U_{n,m}) > \varepsilon^2 \big)
  &= \bar{\PP}_{p}\big( |\hat{U}  - U_{n,m}|^2 > (n\lambda)^{2\alpha}\varepsilon^2 \big)\\
  &= \bar{\PP}_{p}\big( |\hat{K} - K_{n + m}| > (n\lambda)^{\alpha}\varepsilon \big)
\end{align*}
The second step consists on showing that predicting $K_{n+m}$ is always as much
as hard as predicting its conditional expectation
$\tilde{K}_{n,m}(p) \coloneqq \bar{\EE}_{p}[K_{n+m} \mid \Pi_n]$.
This trick was also used by \citet{Orl(17)} where they obtain a
similar result in expectation using the convexity of the loss function and
Jensen's inequality. Here we need a result that holds in probability, but the
idea the similar. In fact, let consider the probability of the events
$\Omega_1 \coloneqq \Set{\Pi \given \bar{\EE}_{p}[|\hat{K} - K_{n+m}| \mid \Pi_n] > (n\lambda)^{\alpha}\varepsilon }$
and $\Omega_2 \coloneqq \Set{\Pi \given |\hat{K} - K_{n+m}| > (n\lambda)^{\alpha}\varepsilon}$. Then,
\begin{align*}
  \bar{\PP}_{p}(\Omega_1)%
  &= \bar{\PP}_{p}(\Omega_1 \cap \Omega_2) + \bar{\PP}_{p}(\Omega_1\cap \Omega_2^c)%
    = \bar{\PP}_{p}(\Omega_1 \cap \Omega_2)%
    \leq \bar{\PP}_{p}(\Omega_2),
\end{align*}
where the second equality follows because it must be the case that
$\Omega_1 \cap \Omega_2^c = \varnothing$. In other words,
\begin{align*}
  \bar{\PP}_{p}\big( |\hat{K} - K_{n+m}| > (n\lambda)^{\alpha}\varepsilon \big)%
  &\geq \bar{\PP}_{p}\big( \bar{\EE}_{p}[|\hat{K} - K_{n+m}| \mid \Pi_n] > (n\lambda)^{\alpha}\varepsilon \big).
\end{align*}
Furthermore, by Jensen's inequality and by the fact that $\hat{K}$ depends only
on $\Pi_n$, it is true that
$\bar{\EE}_{p}[|\hat{K} - K_{n+m}| \mid \Pi_n] \geq |\hat{K} - \bar{\EE}_{p}[K_{n+m} \mid \Pi_n]|$,
which entails
\begin{align*}
  \bar{\PP}_{p}\big( |\hat{K} - K_{n+m}| > (n\lambda)^{\alpha}\varepsilon \big)%
  &\geq \bar{\PP}_{p}\big( |\hat{K} - \tilde{K}_{n,m}(p)| > (n\lambda)^{\alpha}\varepsilon \big).
\end{align*}
The third step of the reduction trades $\tilde{K}_{n,m}(p)$ for its
expectation
$\bar{K}_{n,m}(p) \coloneqq \bar{\EE}_{p}[\tilde{K}_{n,m}]$; observe
that we also have $\bar{K}_{n,m}(p) = \bar{\EE}_{p}[K_{n+m}]$. In the
sequel we will also write $\bar{K}_n(p) \coloneqq \bar{K}_{n,0}(p)$
for the expectation of $K_n$. In fact, rather than reducing to the problem of estimating
$\bar{K}_{n,m}(p)$, it is more convenient to reduce to the problem of estimating
$R_{n,m}(p) \coloneqq \bar{K}_{n,m}(p)/\bar{K}_n(p)$. From
$\hat{K}$ we can build the estimator $\hat{R} \coloneqq \hat{K}/K_n$ for
$R_{n,m}(p)$, then by triangular inequality
\begin{align*}
  |\hat{K} - \tilde{K}_{n,m}(p)|%
  &= |K_n \hat{R} - K_n R_{n,m}(p) + K_n R_{n,m}(p) - \tilde{K}_{n,m}(p)|\\
  &\geq K_n|\hat{R} - R_{n,m}(p)| - |K_nR_{n,m}(p) - \tilde{K}_{n,m}(p)|.
\end{align*}
So we have that $K_n|\hat{R} - R_{n,m}(p)| > 2(n\lambda)^{\alpha} \varepsilon$ and
$|\tilde{K}_{n,m}(p) - K_nR_{n,m}(p)| \leq (n\lambda)^{\alpha} \varepsilon$ implies
$|\hat{K} - \tilde{K}_{n,m}(p)| > (n\lambda)^{\alpha}\varepsilon$, which means that
$\bar{\PP}_{p}( \ell_{\alpha}(\hat{U},U_{n,m}) > \varepsilon^2)$ is always at least
\begin{align}
  \label{eq:135}
  \bar{\PP}_{p}\Big( |\hat{R} - R_{n,m}(p)| > \frac{2(n\lambda)^{\alpha}\varepsilon}{K_n},\, |\tilde{K}_{n,m}(p) - K_nR_{n,m}(p)| \leq (n\lambda)^{\alpha}\varepsilon \Big).
\end{align}
But by triangular inequality again,
\begin{align*}
  |\tilde{K}_{n,m}(p) - K_nR_{n,m}(p)|%
  &\leq |\tilde{K}_{n,m}(p) - \bar{K}_n(p)R_{n,m}(p)|%
    + |\bar{K}_n(p)R_{n,m}(p) - K_nR_{n,m}(p)|\\
  &= |\tilde{K}_{n,m}(p) - \bar{K}_{n,m}(p)|%
    + R_{n,m}(p)|K_n - \bar{K}_{n}(p)|.
\end{align*}
The last display, together with Chebychev's inequality and Lemma~\ref{lem:7}, implies
that for any $p\in \mathcal{P}(\alpha,L,L')$
\begin{align*}
  &\bar{\PP}_{p}\big(|\tilde{K}_{n,m} - K_nR_{n,m}(p)| > (n\lambda)^{\alpha}\varepsilon \big)\\
  &\qquad\leq \bar{\PP}_{p}\Big(|\tilde{K}_{n,m} - \bar{K}_{n,m}(p)| > \frac{(n\lambda)^{\alpha}\varepsilon}{2} \Big)%
    + \bar{\PP}_{p}\Big(|K_n - \bar{K}_{n}(p)| > \frac{(n\lambda)^{\alpha}\varepsilon}{2R_{n,m}(p)} \Big)\\
  &\qquad\leq%
    \frac{4\overline{\var}_{p}(\tilde{K}_{n,m})}{(n\lambda)^{2\alpha}\varepsilon^2}%
    + \frac{4R_{n,m}(p)^2\overline{\var}_{p}(K_n)}{(n\lambda)^{2\alpha}\varepsilon^2}\\
  &\qquad\lesssim%
    \frac{1}{\lambda (n\lambda)^{\alpha}\varepsilon^2}%
    + \frac{R_{n,m}(p)^2}{\lambda^{2\alpha}n^{\alpha}\varepsilon^2}.
\end{align*}
When $p\in \mathcal{P}(\alpha,L,L')$, we get by Lemma~\ref{lem:unif-expect-cumulfreq} that
$R_{n,m}(p) = \lambda^{\alpha}(1+o(1))$ as $n\to \infty$ [recall that $K_n = C_{n,1}$]. Hence, there is a constant $c_1 > 0$
such that
\begin{align}
  \label{eq:154}
  \bar{\PP}_{p}\big(|\tilde{K}_{n,m} - K_nR_{n,m}(p)| > (n\lambda)^{\alpha}\varepsilon \big)%
  &\leq \frac{c_1}{n^{\alpha}\varepsilon^2}.
\end{align}
Since by assumption $\log(\lambda) \geq b \log(n)$, if the constant $b$ is large enough
we can make $\varepsilon$ sufficiently large so that the previous display smaller than
$1/16$. We note that in order to reduce the problem to the estimation of $\alpha$ it
is necessary to take $G_1 \propto \sqrt{D}$ where $D$ is the constant in Theorem~\ref{thm:5},
so that our only degree of freedom to make the rhs of \eqref{eq:154} smaller than
$1/16$ is to have $b$ large enough. Since for any events $A,B$ we have
$\bar{\PP}_{p}(A \cap B) = \bar{\PP}_{p}(A) + \bar{\PP}_{p}(B) - \bar{\PP}_{p}(A \cup B) \geq \bar{\PP}_{p}(A) + \bar{\PP}_{p}(B)-1 = \bar{\PP}_{p}(A) - \bar{\PP}_{p}(B^c)$,
we deduce from equations~\eqref{eq:135} and \eqref{eq:154} that for all $\hat{U}$ and all
$p\in \mathcal{P}(\alpha,L,L')$
\begin{align*}
  \bar{\PP}_{p}\big( \ell_{\alpha}(\hat{U},U_{n,m}) > \varepsilon^2 \big)
  &\geq  \bar{\PP}_{p}\Big( |\hat{R} - R_{n,m}(p)  | > \frac{2(n\lambda)^{\alpha}\varepsilon}{K_n} \Big)
    - \frac{1}{16}.
\end{align*}
Also, we already know that $K_n \sim L\Gamma(1 - \alpha)n^{\alpha}$ almost-surely if
$p\in \mathcal{P}(\alpha,L,L')$. So in fact if $n$ is large enough there is a
constant $c_2 > 0$ such that for all $\hat{U}$ and all
$p\in \mathcal{P}(\alpha,L,L')$
\begin{equation}
  \label{eq:155}
  \bar{\PP}_{p}\big( \ell_{\alpha}(\hat{U},U_{n,m}) > \varepsilon^2 \big)
  \geq  \bar{\PP}_{p}\big( |\hat{R} - R_{n,m}(p)  | > 2 c_2 \lambda^{\alpha}\varepsilon \big)
    - \frac{1}{8}.
\end{equation}
The fourth step consists on reducing to the problem of estimating the functional
$p \mapsto \lambda^{\alpha(p)}$. Indeed, we have already shown that
$R_{n,m}(p) = \lambda^{\alpha}(1+o(1))$ which is a consequence of
Lemma~\ref{lem:unif-expect-cumulfreq}. Indeed, by being more careful we see that
in fact $R_{n,m}(p) = (1+\lambda)^{\alpha}(1 + O(n^{-\alpha/2}\sqrt{\log(n)}))$
whenever $p\in \mathcal{P}(\alpha,L,L')$. The $\sqrt{\log(n)}$ is a bit
problematic as it would make the reduction scheme work only for
$\varepsilon \geq n^{-\alpha/2} \sqrt{\log(n)}$, that is only for
$\log(\lambda) \gg \log(n)^{3/2}$. In order to weaken to
$\log(\lambda) \gg \log(n)$, we leverage that the maximum risk over
$\mathcal{P}(\alpha,L,L')$ has to be greater than the maximum risk over
$\mathcal{P}(\alpha,L,L') \cap \mathcal{R}$, where we define
\begin{equation*}
  \mathcal{R} \coloneqq%
  \Set*{ p \given |R_{n,m}(p) - (1+\lambda)^{\alpha(p)}| \leq c_2\lambda^{\alpha(p)}\varepsilon  }.
\end{equation*}
Later we will have to take additional care to handle the fact that we restrict
the class of parameters to a better behaved subclass, but for now it is always
true that from \eqref{eq:155} and the definition of $\mathcal{R}$ that
\begin{equation}
  \label{eq:173}
  \inf_{\hat{U}}\sup_{\substack{(\alpha,L,L')\in \mathbb{A}\\p \in \mathcal{P}(\alpha,L,L')}}%
  \bar{\PP}_{p}\big(\ell_{\alpha}(\hat{U},U_{n,m}) > \varepsilon^2\big)%
  \geq \inf_{\hat{R}}\sup_{\substack{(\alpha,L,L')\in \mathbb{A}\\p \in \mathcal{P}(\alpha,L,L')\cap \mathcal{R}}}\bar{\PP}_{p}\big(|\hat{R} - (1+\lambda)^{\alpha}| > c_2 \lambda^{\alpha} \varepsilon \big) - \frac{1}{8}.
\end{equation}
Observe that, for every $c_3 > 0$,
\begin{align*}
  &|\log \hat{R} - \alpha \log (1+\lambda) | > c_3\varepsilon\\%
  &\qquad\implies \log \hat{R} - \alpha \log(1+\lambda) > c_3\varepsilon\quad \mathrm{or}\quad \log \hat{R} - \alpha \log (1+\lambda) < -c_3\varepsilon\\
  &\qquad\implies |\hat{R} - (1+\lambda)^{\alpha}| > \min\{(1+\lambda)^{\alpha}(e^{c_3\varepsilon}-1),(1+\lambda)^{\alpha}(1 - e^{c_3\varepsilon}) \}\\
  &\qquad\implies |\hat{R} - (1+\lambda)^{\alpha}| > (1+\lambda)^{\alpha}c_3\varepsilon e^{-c_3\varepsilon}\\
  &\qquad\implies  |\hat{R} - (1+\lambda)^{\alpha}| > (1+\lambda)^{\alpha}c_3e^{-c_3G_2}\varepsilon,
\end{align*}
where the last line follows because by assumption $\varepsilon \leq G_2$. Hence, by taking
$G_2$ sufficiently small and $c_3$ sufficiently large, we have
\begin{align*}
  |\log \hat{R} - \alpha \log (1+\lambda) | > c_3\varepsilon%
  \implies |\hat{R} - (1+\lambda)^{\alpha}| > c_2 \lambda^{\alpha} \varepsilon.
\end{align*}
Therefore,
\begin{align*}
  \bar{\PP}_{p}\big(|\hat{R} - (1+\lambda)^{\alpha}| > c_2 \lambda^{\alpha} \varepsilon \big)%
  &\geq \bar{\PP}_{p}\big(|\log (\hat{R}) - \alpha \log(1+\lambda) | > c_3  \varepsilon \big)\\
  &\geq \bar{\PP}_{p}\Big\{\Big|\frac{\log (\hat{R})}{\log(1+ \lambda)} - \alpha \Big| > \frac{c_3  \varepsilon}{\log(1+\lambda)} \Big\}\\
  &\geq \bar{\PP}_{p}\Big\{\Big|\frac{\log (\hat{R})}{\log(1+\lambda)} - \alpha \Big| > \frac{c_3  G_1}{n^{\underline{\alpha}/2}\log(n) } \Big\},
\end{align*}
where the last inequality follows because
$\varepsilon \leq G_1\log(1+\lambda)n^{-\underline{\alpha}/2}/\log(n)$ by construction.
Then by \eqref{eq:173}, we have obtained that
\begin{equation*}
  \inf_{\hat{U}}\sup_{\substack{(\alpha,L,L')\in \mathbb{A}\\p \in \mathcal{P}(\alpha,L,L')}}\bar{\PP}_{p}\big(\ell_{\alpha}(\hat{U},U_{n,m}) > \varepsilon^2\big)%
  \geq \inf_{\hat{\alpha}}\sup_{\substack{(\alpha,L,L')\in \mathbb{A}\\p \in \mathcal{P}(\alpha,L,L')\cap \mathcal{R}}}\bar{\PP}_{p}\Big\{|\hat{\alpha} - \alpha| > \frac{c_3  G_1}{n^{\underline{\alpha}/2}\log(n)} \Big\} - \frac{1}{8}.
\end{equation*}
We have finally established that estimating $\alpha$ is an easier task than
predicting $U_{n,m}$. However, we cannot immediately conclude using the
Theorem~\ref{thm:5}, as we have restricted the maximum over
$\mathcal{P}(\alpha,L,L')\cap \mathcal{R}$ instead of $\mathcal{P}(\alpha,L,L')$. But, in
view of the Theorem~\ref{thm:2}, it is obvious that the conclusion of Theorem~\ref{thm:5}
remains true over $\mathcal{P}(\alpha,L,L')\cap \mathcal{R}$ if we establish that
$\nu_0(\mathcal{R})$ and $\nu_1(\mathcal{R})$ can be made arbitrarily close to $1$,
which is the case is $\log(\lambda) \geq b \log(n)$ for suitably large $b$ by Lemma~\ref{lem:10}.

\subsection{Auxiliary results}
\label{sec:auxiliary-results-1}

\begin{lem}
  \label{lem:7}
 For any $p \in \mathcal{P}(\alpha,L,L')$ it holds
  \begin{equation*}
    \bar{\EE}_{p}[\tilde{K}_{n,m}]%
    =L\Gamma(1-\alpha)n^{\alpha}(1+\lambda)^{\alpha} + O\big(  (n\lambda)^{\alpha/2} \sqrt{\log(n\lambda)} ),
  \end{equation*}
  and,
  \begin{equation*}
    \overline{\var}_{p}(\tilde{K}_{n,m})%
    \leq \frac{\bar{\EE}_{p}[K_{n+2m}]}{1 + 2\lambda}%
    \lesssim \frac{(n\lambda)^{\alpha}}{\lambda}.
  \end{equation*}
\end{lem}
\begin{proof}
  We observe that $\bar{\EE}_{p}[\tilde{K}_{n,m}] = \bar{\EE}_{p}[K_{n+m}]$, so
  that the first relation simply follows from
  Lemma~\ref{lem:unif-expect-cumulfreq}. For the second relation, we proceed
  using an Efron-Stein argument, as in the Lemma~\ref{lem:var-cumulfreq}. We write
  $\Delta_j \coloneqq \sum_{i=n+1}^{n+m}\1_{X_i=j}$, and we observe that
  $Y_{n+m,j} \geq Y_{n,j}$ almost-surely. Then,
  \begin{align*}
    K_{n+m}
    &= \sum_{j\geq 1}\1_{Y_{n+m,j}\geq 1}\\
    &= \sum_{j\geq 1}\1_{Y_{n,j}\geq 1}\1_{Y_{n+m,j}\geq 1}%
      + \sum_{j\geq 1}\1_{Y_{n,j}=0}\1_{Y_{n+m,j}\geq 1}\\
    &= K_n + \sum_{j\geq 1}\1_{Y_{n,j}=0}\1_{\Delta_j \geq 1}\\
    &= K_n + \sum_{j\geq 1}\1_{Y_{n,j}=0}(1 - \1_{\Delta_j = 0}).
  \end{align*}
  It follows,
  \begin{align*}
    \tilde{K}_{n,m} = \EE[K_{n+m} \mid \Pi_n]%
    &= K_n + \sum_{j\geq 1}\1_{Y_{n,j}=0}(1 - (1-p_j)^m)\\
    &= K_n + \sum_{j\geq 1}(1 - \1_{Y_{n,j}\geq 1})(1 - (1-p_j)^m)\\
    &= K_n + \sum_{j\geq 1}(1 - (1-p_j)^m)%
      - \sum_{j \geq 1}\1_{Y_{n,j}\geq 1}(1 - (1-p_j)^m)\\
    &=  \sum_{j\geq 1}(1 - (1-p_j)^m)%
      + \sum_{j\geq 1}\1_{Y_{n,j}\geq 1}(1-p_j)^m
  \end{align*}
  Then, with the same notations as in Lemma~\ref{lem:var-cumulfreq}, we can bound
  \begin{align*}
    \overline{\var}_{p}(\tilde{K}_{n,m})%
    &= \overline{\var}_{p}\Big(\sum_{j\geq 1}\1_{Y_{n,j}\geq 1}(1-p_j)^m \Big)\\
    &\leq \sum_{i=1}^n\bar{\EE}_{p}\Big[\Big( \sum_{j\geq 1}\big(\1_{Y_{n,j}\geq 1} - \1_{Y_{n,j}^{(i)}\geq 1} \big)(1-p_j)^m \Big)^2 \Big]\\
    &= \sum_{i=1}^n\bar{\EE}_{p}\Big[\Big( \sum_{j\geq 1}\1_{Y_{n,j}=1}\1_{X_i=j}(1-p_j)^m \Big)^2 \Big]\\
    &= \sum_{i=1}^n\bar{\EE}_{p}\Big[\sum_{j\geq 1}\1_{Y_{n,j}=1}\1_{X_i=j}(1-p_j)^m\sum_{k\geq 1}\1_{Y_{n,k}=1}\1_{X_i=k}(1-p_k)^m \Big]\\
    &= \bar{\EE}_{p}\Big[\sum_{j\geq 1}\1_{Y_{n,j}=1}\sum_{i=1}^n\1_{X_i=j}(1-p_j)^{2m}\Big]\\
    &= \sum_{j\geq 1} \bar{\EE}_{p}[\1_{Y_{n,j}=1}](1-p_j)^{2m},
  \end{align*}
  where the last line follows because $Y_{n,j} = \sum_{i=1}^n\1_{X_i=j}$ by
  construction. Now $Y_{n,j}=1$ if and only if all the $X_1,\dots,X_n$ are not
  equal to $j$ but one, which can happens only in $n$ ways (either $X_1 = j$ and
  all the others are not equal to $j$, or $X_{2}=j$ and the other are not equal
  to $j$, etc..). We deduce that
  \begin{align*}
    \overline{\var}_{p}(\tilde{K}_{n,m})%
    &\leq \sum_{j\geq 1} np_j(1-p_j)^{n+2m - 1}\\
    &= \frac{n}{n+2m} \bar{\EE}_{p}[M_{n+2m,1}]\\
    &\leq \frac{\bar{\EE}_{p}[K_{n+2m}]}{1+2\lambda},
  \end{align*}
  which concludes the proof.
\end{proof}

\begin{lem}
  \label{lem:10}
  For every $\eta > 0$ there exists $C > 0$ and $n_0 > 0$ such that if
  $\varepsilon \geq Cn^{\alpha_j/2}$ and $n \geq n_0$, then
  \begin{align*}
    \min\{\nu_0(\mathcal{R}),\, \nu_1(\mathcal{R})\} \geq 1 - \eta.
  \end{align*}
\end{lem}
\begin{proof}
  Let $j=0$ or $j=1$ arbitrary. Recall that under $\nu_j$ we have
  $\alpha(P) = \alpha_j$ almost-surely. Then, we observe that
  \begin{align}
    \label{eq:184}
    |R_{n,m}(p) - (1+\lambda)^{\alpha_j}|%
    \leq c_2 \lambda^{\alpha_j}\varepsilon%
    &\iff |\bar{K}_{n,m}(p) - (1+\lambda)^{\alpha_j}\bar{K}_n(p)| \leq c_2\lambda^{\alpha_j}\bar{K}_n(p)\varepsilon.
  \end{align}
  Recall that $\bar{K}_n(p) = n\int_0^1 \bar{F}_{p}(x)(1-x)^{n-1}\intd x$; see
  for instance the proof of Proposition~\ref{pro:6}. Recall that conditional on
  $(J_i)_{i\geq 1}$ we have $\bar{F}_P(x) = \sum_{i\geq 1}\1_{J_i > x T}$ with
  $T = \sum_{i\geq 1}J_i$. Then, we can decompose $\bar{K}_{n}(P)$ as
  \begin{align}
    \label{eq:185}
    \bar{K}_n(P)%
    &= n \int_0^1 \bar{\rho}_j(xT)(1-x)^{n-1}\intd x%
      + \underbrace{n \int_0^1\{\bar{F}_{P}(x) - \bar{\rho}_j(xT)\}(1-x)^{n-1}\intd x}_{\Delta_n(P)}.
  \end{align}
  We now show that $\Delta_n(P) = O_p(n^{\alpha_j/2})$ where the $O_p$ is understood
  with respect to $\nu_j$. Indeed,
  \begin{align*}
    |\Delta_n(P)|%
    &\leq n \int_0^1\Big| \sum_{i\geq 1}\1_{J_i> xT} - \bar{\rho}_j(xT)\Big| (1-x)^{n-1}\\
    &= \frac{n}{T} \int_0^T\Big| \sum_{i\geq 1}\1_{J_i> x} - \bar{\rho}_j(x)\Big| e^{-(n-1)x/T}\intd x
  \end{align*}
  Recall that $T$ has an $\alpha_j$-stable distribution. Then, for any $\eta > 0$ we can
  find $t_0,t_1 > 0$ such that with $\nu_j$-probability at least $1-\eta$ we have
  $t_0 \leq T \leq t_1$. Moreover, on the event that $t_0 \leq T \leq t_1$, we have
  \begin{align*}
    |\Delta_n(P)|%
    \leq \frac{n}{t_0}\int_0^{\infty}\Big| \sum_{i\geq 1}\1_{J_i> x} - \bar{\rho}_j(x)\Big| e^{-(n-1)x/t_1}\intd x.
  \end{align*}
  So it is enough to show that the above integral is a $O_p(n^{\alpha_j}/2)$. But,
  under $\nu_j$, the variable $\sum_{i\geq 1}\1_{J_i > x}$ has a Poisson distribution
  with mean $\bar{\rho}_j(x)$. Therefore,
  \begin{align*}
    \EE_{\nu_j}\Big[\int_0^{\infty}\Big| \sum_{i\geq 1}\1_{J_i> x} - \bar{\rho}_j(x)\Big| e^{-(n-1)x/t_1}\intd x\Big]
    &= \int_0^{\infty}\EE_{\nu_j}\Big[\Big| \sum_{i\geq 1}\1_{J_i> x} - \bar{\rho}_j(x)\Big|\Big] e^{-(n-1)x/t_1}\intd x\\
    &\leq \int_0^{\infty} \bar{\rho}_j(x)^{1/2} e^{-(n-1)x/t_1}\intd x.
  \end{align*}
  But $\bar{\rho}_j(x) \lesssim x^{-\alpha_j}\1_{0< x \leq 1}$, which follows from
  \eqref{eq:rhobar}. This establish the claim that $\Delta_n(P) = O_p(n^{\alpha_j/2})$,
  because
  \begin{align*}
    \int_0^{1} x^{-\alpha_j/2}e^{-(n-1)x/t_1}\intd x%
    \leq \Big\{ \frac{n-1}{t_1} \Big\}^{\alpha_j/2 - 1}\Gamma(1 + \alpha_j/2).
  \end{align*}
  On the other hand, writing $C_j \coloneqq \frac{\alpha_0/\alpha_j}{\Gamma(1-\alpha_j)}$, we have
  again by \eqref{eq:rhobar} that
  \begin{align*}
    n\int_0^1\bar{\rho}_j(xT)(1 - x)^{n-1}\intd x%
    &= C_j n\int_0^1 (xT)^{-\alpha_j}(1 - (xT)^{\alpha_j})\1_{0< xT \leq 1}(1-x)^{n-1}\intd x\\
    &= C_j n\int_0^1 (xT)^{-\alpha_j}\1_{xT \leq 1}(1-x)^{n-1}\intd x%
      + C_jn \int_0^1 \1_{xT\leq 1}(1 - x)^{n-1}\\
    &= C_j n\int_0^1 (xT)^{-\alpha_j}\1_{xT \leq 1}(1-x)^{n-1}\intd x%
      + O(1).
  \end{align*}
  We consider two cases, according to whether $T \leq 1$ or not. If $T \leq 1$, then
  \begin{align*}
    n\int_0^1 (xT)^{-\alpha_j}\1_{xT \leq 1}(1-x)^{n-1}\intd x%
    &= \frac{n}{T^{\alpha_j}}\int_0^1 x^{-\alpha_j}(1-x)^{n-1}\intd x\\
    &= \frac{\Gamma(1-\alpha_j)}{T^{\alpha_j}} \frac{n!}{\Gamma(n+1 -\alpha)}\\
    &= \frac{\Gamma(1-\alpha_j)n^{\alpha_j}}{T^{\alpha_j}} + O(T^{-\alpha_j}).
  \end{align*}
  If $T > 1$, then
  \begin{align*}
    n\int_0^1 (xT)^{-\alpha_j}\1_{xT \leq 1}(1-x)^{n-1}\intd x%
    &= \frac{n}{T^{\alpha_j}}\int_0^{1/T}x^{-\alpha_j}(1-x)^{n-1}\intd x\\
    &= \frac{n}{T^{\alpha_j}}\int_0^{1}x^{-\alpha_j}(1-x)^{n-1}\intd x - \frac{n}{T^{\alpha_j}}\int_{1/T}^1x^{-\alpha_j}(1-x)^{n-1}\intd x\\
    &= \frac{\Gamma(1-\alpha_j)n^{\alpha_j}}{T^{\alpha_j}} + O(T^{-\alpha_j}) - \frac{n}{T^{\alpha_j}}\int_{1/T}^1x^{-\alpha_j}(1-x)^{n-1}\intd x\\
    &= \frac{\Gamma(1-\alpha_j)n^{\alpha_j}}{T^{\alpha_j}} + O(1) - \frac{n}{T^{\alpha_j}}\int_{1/T}^1x^{-\alpha_j}(1-x)^{n-1}\intd x
  \end{align*}
  On the event that $t_0 \leq T \leq t_1$, it is clear that the integral in the rhs of
  the last display is a $o(1)$ as $n \to \infty$, so that we have shown that in both
  the cases we have
  \begin{align*}
    n \int_0^1\bar{\rho}_j(xT)(1 - x)^{n-1}\intd x%
    &= \frac{C_j\Gamma(1-\alpha_j)n^{\alpha_j}}{T^{\alpha_j}} + O_p(1)\\
    &= \frac{\alpha_0/\alpha_j}{T^{\alpha_j}}n^{\alpha_j} + O_p(1).
  \end{align*}
  Combining this result with \eqref{eq:185}, it follows that under $\nu_j$,
  \begin{align*}
    \bar{K}_n(P) = \frac{\alpha_0/\alpha_j}{T^{\alpha_j}}n^{\alpha_j} + O_p(n^{\alpha_j/2}).
  \end{align*}
  Since $\bar{K}_{n,m}(P) = \bar{K}_{n+m}(P)$ and since $m = n\lambda$, we
  deduce that under $\nu_j$
  \begin{align*}
    \bar{K}_{n,m}(P) - (1+\lambda)^{\alpha_j}\bar{K}_n(P)%
    &= \frac{\alpha_0/\alpha_j}{T^{\alpha_j}}\Big\{ (n+m)^{\alpha_j} - (1+\lambda)^{\alpha_j}n^{\alpha_j}  \Big\}%
      + O_p\big( n^{\alpha_j/2} \lambda ^{\alpha_j} \big)\\
    &= O_p\Big( \frac{\lambda^{\alpha_j} \bar{K}_n(P)}{n^{\alpha_j/2}} \Big)
  \end{align*}
  It follows that if $\varepsilon \geq C n^{\alpha_j/2}$ for a large enough constant $C > 0$, then
  with $\nu_j$ probability at least $1-\eta$ it will be true that
  $|\bar{K}_{n,m}(P) - (1+\lambda)^{\alpha_j}\bar{K}_n(P)| \leq c_2 \lambda^{\alpha_j}\bar{K}_n(P)\varepsilon$,
  and then the conclusion follows from \eqref{eq:184}.
\end{proof}

\section{Method of the two fuzzy hypotheses for partition models}
\label{sec:method-two-fuzzy}

This section contains the proof of Theorem~\ref{thm:2} and of Proposition~\ref{thm:4}.

\subsection{Proof of Theorem~\ref{thm:2}}
\label{sec:thm:2}

The proof is a minor adaptation of \cite[Theorems 2.14 and 2.15]{tsybakov2009}.
For any estimator $\hat{\delta}$,
\begin{align*}
  \sup_{p\in \mathcal{P}'} \bar{\PP}_{p}\Big(|\hat{\delta}(\Pi_n) - \delta(\Pi,p)|\geq  \varepsilon \Big)%
  &\geq \frac{1}{2}\sum_{j=0}^1\int_{\mathcal{P}'}\bar{\PP}_{p}\Big(|\hat{\delta}(\Pi_n) - \delta(\Pi,p)|\geq  \varepsilon \Big) \frac{\nu_j(\intd p)}{\nu_j(\mathcal{P}')}\\%
  &\geq \frac{1}{2} \sum_{j=0}^1 \int \bar{\PP}_{p}\Big(|\hat{\delta}(\Pi_n) - \delta(\Pi,p)|\geq  \varepsilon \Big)\nu_j(\intd p) - \xi\\
  &= \frac{1}{2}\sum_{j=0}^1F_j\Big(|\hat{\delta}(\Pi_n) - \delta(\Pi,P)|\geq  \varepsilon \Big) - \xi\\
\end{align*}
But, because for any events $A,B$ we have
$F_j(A\cap B) = F_j(A) + F_j(B) - F_j(A\cup B) \geq F_j(A) + F_j(B)-1 = F_j(A) - F_j(B^c)$,
we get%
\begin{align*}
  F_0\Big(|\hat{\delta}(\Pi_n) - \delta(\Pi,P)|\geq  \varepsilon \Big)%
  &\geq F_0\big(\hat{\delta}(\Pi_n) \geq c+\varepsilon,\ \delta(\Pi,P) \leq c \big)\\
  &\geq F_0\big(\hat{\delta}(\Pi_n) \geq c+\varepsilon \big) - F_0(\delta(\Pi,P) > c  )\\
  &\geq F_0\big(\hat{\delta}(\Pi_n) \geq c+\varepsilon \big) - \kappa\\
  &= \nu_0(\PP_{P}^n)\big(\hat{\delta}(\Pi_n) \geq c+\varepsilon \big) - \kappa,
\end{align*}
and,
\begin{align*}
  F_1\Big(|\hat{\delta}(\Pi_n) - \delta(\Pi,P)|\geq  \varepsilon \Big)
  &\geq F_1\big(\hat{\delta}(\Pi_n) < c+\varepsilon,\ \delta(\Pi,P) \geq c + 2\varepsilon \big)\\
  &\geq F_1\big(\hat{\delta}(\Pi_n) < c+\varepsilon\big) - F_1\big(\delta(\Pi,P) < c + 2\varepsilon \big)\\
  &\geq \nu_1(\PP_{P}^n)\big(\hat{\delta}(\Pi_n) < c+\varepsilon\big) - \kappa.
\end{align*}
Then,
\begin{align*}
  \sup_{p\in \mathcal{P}'} \bar{\PP}_{p}\big(|\hat{\delta}(\Pi_n) - \delta(\Pi,p)|\geq  \varepsilon \big)%
  &\geq \frac{1}{2}\Big( \nu_0(\PP_{P}^n)(\hat{\delta} \geq c +\varepsilon) + \nu_1(\PP_{P}^n)(\hat{\delta}< c+ \varepsilon)  \Big) - \kappa - \xi\\
  &\geq \frac{1}{2}\Big(1 - \|\nu_0(\PP_{P}^n) - \nu_1(\PP_{P}^n)\|_{\mathrm{TV}}\Big) - \kappa - \xi.
\end{align*}

\subsection{Proof of Proposition~\ref{thm:4}}
\label{sec:thm:4}

We note that $\sum_{j=1}^nM_{n,j} \leq n$ and $\sum_{j=1}^j jM_{n,j} = n$ and we will
write $\mathcal{M}_n$ the set of vectors satisfying these constraints. By
\cite[Chapter~2]{Pit(06)} the statistic $\M_n$ is sufficient for $\rho$, so
that the KL divergence between the laws of $\Pi_n$ is the same as the KL
divergence between the laws of $\M_n$. We will abusively write $\PP_{\rho}^n$ to
denote either the law of $\Pi_n$ or the law of $\M_n$ under the model $\rho$. We
write $\QQ_{\rho}^n$ the joint-distribution of $(\M_n,U_n)$ under the model $\rho$. We
obtain our first inequality using the chain-rule for relative entropies
\begin{align*}
  \KL(\PP_{\rho_0}^n;\PP_{\rho_1}^n)%
  &\leq \KL(\QQ_{\rho_0}^n;\QQ_{\rho_1}^n).
\end{align*}
Let define the so-called Laplace exponents $\psi_j(u) \coloneqq \int_0^{\infty}(1 - e^{-u x})\rho_j(\intd x)$, $j=0,1$. In view of
Propositions~\ref{pro:1} and~\ref{pro:2}, for all $\bm{m} \in \mathcal{M}_n$,
\begin{align*}
  \QQ_{\rho_0}^n\big(\M_n = (m_1,\dots,m_n),\, U_n \in \intd u\big)%
  &= \frac{n!}{\prod_{i=1}^n(i!)^{m_i}m_i!} \frac{u^{n-1}e^{-\psi_0(u)}}{\Gamma(n)}\prod_{\ell=1}^n\mu_0(u,\ell)^{m_{\ell}}\,\intd u,
\end{align*}
and similarly under $\rho_1$. Therefore,
\begin{align*}
  \KL(\PP_{\rho_0}^n;\PP_{\rho_1}^n)%
  &\leq \EE_{\rho_0}^n\Big[- \psi_0(U_n) + \psi_1(U_n) + \sum_{\ell=1}^n M_{n,\ell}\log
    \frac{\mu_0(U_n,\ell)}{\mu_1(U_n,\ell)} \Big].
\end{align*}
Hence, to finish the proof it is enough to show that
\begin{align}
  \label{eq:28}
  \EE_{\rho_0}^n[\psi_j(U_n)]%
  &= \sum_{\ell=1}^n \EE_{\rho_0}^n\Big[M_{n,\ell} \frac{\mu_j(U_n,\ell)}{\mu_0(U_n,\ell)}  \Big],\qquad j=0,1.
\end{align}

We now establish \eqref{eq:28}. The starting point is the following equation,
which can be obtained from \cite[Lemma~4.3]{Pit(06)}.
\begin{equation}
  \label{eq:163}
  \PP_{\rho}^n\big( \Pi_n = \Set{A_1,\dots,A_k} \mid T \big)%
  =  \int_{\NNReals^k}\frac{f_T(\rho;T - \sum_{j=1}^ky_j)\1\Set{\sum_{j=1}^ky_j <
  T}}{T^nf_T(\rho;T)}%
  \prod_{j=1}^ky_j^{|A_i|} \prod_{j=1}^k\rho(\intd y_j),
\end{equation}
where $f_T(\rho;\cdot)$ is the density (wrt Lebesge measure) of the variable $T$ under the
model $\rho$. From here, the strategy runs as follows. We first obtain from
\eqref{eq:163} the law of $\M_n \mid U_n$ (Proposition~\ref{pro:1}), then the laws of
$T \mid U_n$ and $U_n$ (Proposition~\ref{pro:2}), from which we determine
$\EE_{\rho}^n[M_{n,\ell} \mid U_n]$ (Proposition~\ref{pro:3}). Finally, we show in Proposition~\ref{pro:5} that
$\EE_{\rho_0}^n[\psi_j(U_n)]$ can be conveniently reexpressed in term of
$\EE_{\rho_0}^n[M_{n,\ell} \mid U_n]$, leading to the equation~\eqref{eq:28}. An important quantity
that will show up everywhere in this section is the function $h_n$ such that
\begin{equation*}
  h_n(\rho;u)%
  \coloneqq \int_0^{\infty}x^n f_T(\rho;x)e^{-u x}\intd x,\qquad u \geq 0.
\end{equation*}

\begin{prp}
  \label{pro:1}
  For all $\bm{m} \in \mathcal{M}_n$, $j=0,1$,
  \begin{align*}
    \PP_{\rho_j}^n\big( \M_n = (m_1,\dots,m_k) \mid U_n \big)%
    =  \frac{n!}{\prod_{i=1}^n(i!)^{m_i}m_i!} \frac{e^{-\psi_j(U_n)}}{h_n(\rho_j;U_n)}\prod_{\ell=1}^n \mu(U_n,\ell)^{m_{\ell}}.
  \end{align*}
\end{prp}
\begin{proof}
  By equation~\eqref{eq:163} and by Proposition~\ref{pro:2},
  \begin{align*}
    &\PP_{\rho_j}^n\big( \Pi_n = \Set{A_1,\dots,A_k} \mid U_n \big)\\%
    &\qquad= %
      \int_{\NNReals^k}\int_{\NNReals} \frac{f_T(\rho_j;t -\sum_{j=1}^ky_k)\1\Set{t -
      \sum_{j=1}^k y_j > 0}}{h_{n}(\rho_j;U_n)}e^{-U_nt}\intd t \prod_{j=1}^k
      y_j^{|A_j|}\prod_{j=1}^k \rho_j(\intd y_j)\\
    &\qquad=  \int_{\NNReals^k} \int_{\NNReals} \frac{f_T(\rho_j;t)e^{-U_n
      t} \,\intd t}{h_{n}(\rho_j;U_n)}%
      \prod_{j=1}^ky_j^{|A_j|}e^{-U_n y_j} \prod_{j=1}^k \rho_j(\intd y_j)\\
    &\qquad=  \frac{e^{-\psi_j(U_n)}}{h_n(\rho_j;U_n)}\prod_{j=1}^k \mu_j(U_n,|A_j|).
  \end{align*}
  The conclusion follows from the previous result by taking a partition
  $\Set{A_1,\dots,A_k}$ with $m_1$ blocks of size $1$, $m_2$ blocks of size $2$,
  etc.. and multiplying by the number of such partitions. It can be seen that
  there are exactely $\frac{n!}{\prod_{i=1}^n(i!)^{m_i}m_i!}$ such partitions
  \cite[see for instance][Section~2]{Pit(06)}.
\end{proof}

\begin{prp}
  \label{pro:2}
  The following are true.
  \begin{align*}
    \PP_{\rho_j}^n(T \in \intd t \mid U_n)%
    &= \frac{t^nf_T(\rho_j;t)e^{-U_nt}\intd t}{h_n(\rho_j;U_n)},
  \end{align*}
  and,
  \begin{align*}
    \PP_{\rho_j}^n (U_n \in \intd u)%
    &= \frac{u^{n-1}h_n(\rho_j;u)\intd u}{\Gamma(n)}.
  \end{align*}
\end{prp}
\begin{proof}
  We have the first trivial relation,
  \begin{align*}
    \PP_{\rho_j}^n(0 < T < a,\, 0 < U_n < b)%
    &= \PP_{\rho_j}^n(0 < T < a,\, 0 < G_n < b T)\\
    &= \int_0^a f_T(\rho_j;t) \int_0^{bt} \frac{u^{n-1}e^{-u}}{\Gamma(n)}\intd u \intd
   t\\
    &= \int_0^a\int_0^b t^n f_T(\rho_j;t) \frac{u^{n-1}e^{-ut}}{\Gamma(n)}\intd u \intd t.
  \end{align*}
  Whence the result after conditioning.
\end{proof}

\begin{prp}
  \label{pro:3}
  It holds
  \begin{align*}
    \EE_{\rho_j}^n[M_{n,\ell} \mid U_n]%
    &= \frac{h_{n-\ell}(\rho_j;U_n)}{h_n(\rho_j;U_n)} \binom{n}{\ell} \mu_j(U_n,\ell).
  \end{align*}
\end{prp}
\begin{proof}
  We obtain the result in three steps. We compute
  $\EE_{\rho_j}^n[M_{n,\ell} \mid (J_j)_{j\geq 1}]$, then
  $\EE_{\rho_j}^n[M_{n,\ell} \mid T]$, and finally
  $\EE_{\rho_j}^n[M_{n,\ell} \mid U_n]$. Obviously,
  \begin{align*}
    \EE_{\rho_j}^n[M_{n,\ell} \mid (J_j)_{j\geq 1}]%
    &= \sum_{j\geq 1}\EE_{\rho_j}^n[\1_{Y_{n,j}=\ell} \mid (J_j)_{j\geq 1}]\\
    &= \binom{n}{\ell} \sum_{i\geq 1} \Big(\frac{J_i}{T}\Big)^{\ell}\Big(1 -
      \frac{J_i}{T}\Big)^{n-\ell}\\
    &= T^{-n}\binom{n}{\ell} \sum_{i\geq 1} J_i^{\ell}(T - J_i)^{n-\ell}.
  \end{align*}
  It follows using Slivnyak-Mecke's theorem \citet{Dal(07)}, that for any $t > 0$
  \begin{align*}
    \EE_{\rho_j}^n[M_{n,\ell}\1\Set{0 < T < t}]%
    &= \binom{n}{\ell}\EE_{\rho_j}^n\Big[\sum_{i\geq 1} \frac{J_i^{\ell}(\sum_{k\geq
      1}J_k - J_i)^{n-\ell}\1\Set{0 < \sum_{k\geq 1}J_k < t} }{(\sum_{k\geq 1}J_k
      )^n} \Big]\\
    &= \binom{n}{\ell} \int_{\NNReals} x^{\ell}
      \EE_{\rho_j}^n\Big[\frac{T^{n-\ell}\1\Set{0 <T + x < t}}{(T+x)^n}
      \Big]\,\rho_j(\intd x)\\
    &= \binom{n}{\ell} \int_{\NNReals} f_T(\rho_j;u)u^{n-\ell}\int_{\NNReals}
      \frac{x^{\ell}\1\Set{0 < u + x < t}}{(u+x)^n}\rho_j(\intd x)\intd u\\
    &= \binom{n}{\ell} \int_0^t \int_0^v\frac{f_T(\rho_j;v -x)}{v^n} x^{\ell}(v
      -x)^{n-\ell} \rho_j(\intd x)\intd v.
  \end{align*}
  So, by Lebesgue's differentiation theorem, we deduce that,
  \begin{align*}
    \EE_{\rho_j}^n[M_{n,\ell} \mid T =t]%
    &= \lim_{\varepsilon\to 0} \frac{\EE_{\rho_j}^n[M_{n,\ell}\1\Set{t-\varepsilon <
      T < t+\varepsilon}]}{\PP_{\rho_j}^n(t-\varepsilon<T<t+\varepsilon)}\\
    &= \lim_{\varepsilon\to 0} \frac{\binom{n}{\ell}
      \int_{t-\varepsilon}^{t+\varepsilon} \int_0^v \frac{f_T(\rho_j;v -x)}{v^n}
      x^{\ell}(v -x)^{n-\ell} \rho_j(\intd x)\intd
      v}{\int_{t-\varepsilon}^{t+\varepsilon}f_T(\rho_j;v)\intd v}\\
    &= \binom{n}{\ell} \int_0^t \frac{f_T(\rho_j;t -x)}{t^nf_T(\rho_j;t)}x^{\ell}(t-x)^{n-\ell}
      \rho_j(\intd x).
  \end{align*}
  So it remains to compute
  $\EE_{\rho_j}^n[M_{n,\ell} \mid U_n] = \EE_{\rho_j}^n[ \EE_{\rho_j}^n[M_{n,\ell} \mid T]
  \mid U_n ]$. The law of $T \mid U_n$ has been determined in Proposition~\ref{pro:2}. Using this result, we have
  \begin{align*}
    \EE_{\rho_j}^n[M_{n,\ell} \mid U_n]%
    &= \frac{1}{h_n(\rho_j;U_n)} \binom{n}{\ell}%
      \int_{\NNReals} \int_0^t f_T(\rho_j;t - x)e^{-U_n t}x^{\ell}(t - x)^{n-\ell}
      \rho_j(\intd x) \intd t\\
    &= \frac{1}{h_n(\rho_j;U_n)}\binom{n}{\ell}\int_{\NNReals} v^{n-\ell}f_T(\rho_j;v)e^{-U_n v}\intd v
      \int_{\NNReals}x^{\ell} e^{-U_n x} \rho_j(\intd x)\\
    &= \frac{h_{n-\ell}(\rho_j;U_n)}{h_n(\rho_jU_n)} \binom{n}{\ell} \mu_j(U_n,\ell).
  \end{align*}
\end{proof}

\begin{prp}
  \label{pro:5}
  For $j=0,1$,
  \begin{align*}
    \EE_{\rho_0}^n[\psi_j(U_n)]%
    &= \sum_{\ell=1}^n\binom{n}{\ell}\EE_{\rho_0}^n\Big[
      \frac{h_{n-\ell}(\rho_0;U_n)}{h_n(\rho_0;U_n)} \mu_j(U_n;\ell)\Big]%
      = \sum_{\ell=1}^n \EE_{\rho_0}^n\Big[M_{n,\ell} \frac{\mu_j(U_n,\ell)}{\mu_0(U_n,\ell)}  \Big].
  \end{align*}
\end{prp}
\begin{proof}
  We start from the rhs, and we will show that it is equal to the lhs. Using
  the expression of the density of $U_n$ in Proposition~\ref{pro:2}, we have that
  \begin{align*}
    \EE_{\rho_0}^n\Big[ \frac{h_{n-\ell}(\rho_0;U_n)}{h_n(\rho_0;U_n)}\mu_j(U_n,\ell)\Big]%
    &= \frac{1}{\Gamma(n)} \int_{\NNReals}u^{n-1} h_{n-\ell}(\rho_0;u) \mu_j(u,\ell)\,\intd u\\
    &=\frac{1}{\Gamma(n)}
      \int_{\NNReals^3}u^{n-1}x^{n-\ell} f_T(\rho_0;x)e^{-ux}
      y^{\ell}e^{-u y} \rho_j(\intd y)\intd x \intd u.
  \end{align*}
  Consequently,
  \begin{align*}
    &\sum_{\ell=1}^n \binom{n}{\ell} \EE_{\rho_0}^n\Big[
    \frac{h_{n-\ell}(\rho_0;U_n)}{h_n(\rho_0;U_n)}\mu_j(U_n,\ell)\Big]\\%
    &\qquad\qquad= \frac{1}{\Gamma(n)}
      \int_{\NNReals^3} u^{n-1}e^{-u(x+y)}f_T(\rho_0;x)\Big\{
      (x+y)^n - x^n\Big\}\rho_j(\intd y)\intd x \intd u.
  \end{align*}
  We remark that,
  \begin{align*}
    \frac{1}{\Gamma(n)}\int_{\NNReals}f_T(\rho_0;x) (x+y)^n\int_{\NNReals}
    u^{n-1}e^{-u(x+y)}\intd u \intd x%
    &= \int_{\NNReals}f_T(\rho_0;x) \intd x
      = \EE_{\rho_0}^n[1].
  \end{align*}
  and,
  \begin{align*}
    \frac{1}{\Gamma(n)}\int_{\NNReals} u^{n-1}\int_{\NNReals}x^n
    e^{-u(x+y)}f_T(\rho_0;x)\intd x \intd u%
    &= \frac{1}{\Gamma(n)} \int_{\NNReals}u^{n-1}e^{-uy} h_n(\rho_0;u)\intd u\\
    &= \EE_{\rho_0}^n[e^{-U_ny}].
  \end{align*}
  So by Fubini,
  \begin{align*}
    \sum_{\ell=1}^n\binom{n}{\ell}\EE_{\rho_0}^n\Big[
    \frac{h_{n-\ell}(\rho_0;U_n)}{h_n(\rho_0;U_n)}\mu_j(U_n,\ell)\Big]
    &= \int_{\NNReals}\EE_{\rho_0}^n[1 - e^{-U_ny}]\rho_j(\intd y)\\
    &= \EE_{\rho_0}^n\Big[ \int_{\NNReals}(1 -e^{-U_ny})\rho_j(\intd y) \Big]\\
    &= \EE_{\rho_0}^n[\psi_j(U_n)].
  \end{align*}
\end{proof}

\section{Results on regular variations}

\subsection{Folklore results}
\label{sec:folkl-results-regul}

In this section, we state some classical results about the behaviour of $K_n$,
$\Set{M_{n,\ell}\given \ell=1,\dots,n }$ and $\Set{C_{n,k}\given k=1,\dots,n }$
when $\bar{F}_p$ has regularly varying tails. Those results are used multiple
times in the proofs of our theorems. The following lemma is a restatement of
Corollary~21 in \cite{Gne(07)}.

\begin{lem}
  \label{lem:1}
  If there is $\alpha\in (0,1)$ and a non-negative slowly varying function at
  infinity $\ell$ such that $\bar{F}_p(x) \sim x^{-\alpha}\ell(1/x)$ when
  $x \to 0$, then almost-surely as $n \to \infty$
  \begin{equation*}
    K_n \sim \Gamma(1-\alpha)\cdot n^{\alpha}\ell(n),%
  \end{equation*}
  and, for all $\ell \in \Nats$ almost-surely as $n\to \infty$,
  \begin{equation*}
    M_{n,\ell} \sim \frac{\alpha \Gamma(\ell - \alpha)}{\ell!}\cdot n^{\alpha}\ell(n).
  \end{equation*}
\end{lem}

With the same arguments \cite{Gne(07)}, we can easily establish the following
for the random variables $\Set{C_{n,k} \given k=1,\dots,n }$. It is to be noted
that equation~(66) in \cite{Kar(67)} establishes that
$\lim_{n\to\infty}\frac{C_{n,k}}{\bar{\EE}_{p}[C_{n,k}]} = 1$ almost-surely
under no condition on $p$. Surprisingly \cite{Kar(67)} does not
explicit the asymptotic behaviour of $\bar{\EE}_{p}[C_{n,k}]$ under regular
variations; we fill this gap in the following lemma.

\begin{lem}
  \label{lem:asympt-cumul-freq}
  For all $p$ and all $k\in \Nats$,
  $\lim_{n\to \infty} \frac{C_{n,k}}{\bar{\EE}_p[C_{n,k}]} \to 1$ almost-surely,
  where
  \begin{equation*}
    \bar{\EE}_{p}[C_{n,k}]%
    =%
    k \binom{n}{k} \int_0^1 \bar{F}_{p}(x)x^{k-1}(1-x)^{n-k}\intd x.
  \end{equation*}
  Furthermore, if there is $\alpha\in (0,1)$ and a non-negative slowly varying
  function at infinity $\ell$ such that $\bar{F}_p(x) \sim x^{-\alpha}\ell(1/x)$
  when $x \to 0$, then for all $k \in \Nats$
  \begin{equation*}
    \bar{\EE}_p[C_{n,k}] \sim n^{\alpha}\ell(n)\cdot \frac{\Gamma(k-\alpha)}{\Gamma(k)}
  \end{equation*}
  as $n\to \infty$.
\end{lem}
\begin{proof}
  The almost-sure convergence result is proved in \cite{Kar(67)}. Next,
  the starting point is
  \begin{align*}
    \bar{\EE}_{p}[M_{n,\ell}]%
    &= \sum_{j\geq
      1}\bar{\EE}_{p}[\1_{Y_{n,j}=\ell}]\\
    &= \binom{n}{\ell}\sum_{j\geq
      1}p_j^{\ell}(1-p_j)^{n-\ell}\\
    &= \binom{n}{\ell}\sum_{j\geq 1}\int_0^{p_j}\Big\{\ell
      x^{\ell-1}(1-x)^{n-\ell} - (n-\ell)x^{\ell}(1-x)^{n-\ell-1}\Big\}\intd x.
  \end{align*}
  Since $\bar{F}_{p}(x) = \sum_{j\geq 1}\1_{p_j > x}$
  \begin{align*}
    \bar{\EE}_{p}[M_{n,\ell}]%
    &=\binom{n}{\ell}\int_0^{1}\bar{F}_{p}(x)x^{\ell-1}(1-x)^{n-\ell - 1}(\ell -nx)\intd x.
  \end{align*}
  The first expression for the expectation follows after remarking that
  \begin{equation*}
    \sum_{\ell=k}^n \binom{n}{\ell} x^{\ell-1}(1-x)^{n-\ell-1}(\ell - nx)%
    = k \binom{n}{k}x^{k-1}(1-x)^{n-k}.
  \end{equation*}
  Furthermore,
  \begin{equation*}
    \int_0^1 \bar{F}_p(x)x^{k-1}(1 - x)^{n-k}\intd x%
    = \frac{1}{n^k}\int_0^{\infty}\bar{F}_p(un^{-1})u^{k-1}(1 - un^{-1})^{n-k}\1_{u\leq n}\intd u
  \end{equation*}
  But if $\bar{F}_p(x) \sim x^{-\alpha}\ell(x)$, by Potter's bound
  \cite[Theorem~1.5.6]{Bin(89)} for all $\delta > 0$ we can find
  $A,\varepsilon > 0$ such that
  $\bar{F}_p(un^{-1})/[n^{\alpha}\ell(n)] \leq A u^{-\alpha - \delta}$ for all
  $un^{-1} \leq \varepsilon$. If $un^{-1} > \varepsilon$, then
  $\bar{F}_p(un^{-1}) \leq C$ for some universal $C$, and
  $(1- un^{-1})^{n-k} \leq e^{-u (n-k)/n}$. Since
  $n^{\alpha}\ell(n) \to \infty$, this shows that
  $u \mapsto \bar{F}_p(un^{-1})u^{k-1}(1 - un^{-1})^{n-k}\1_{u\leq n}/[n^{\alpha}\ell(n)]$
  is dominated by an integrable function on $\NNReals$. Therefore by the
  dominated convergence theorem,
  \begin{align*}
    &\lim_{n\to \infty}\frac{1}{n^{\alpha}\ell(n)}\int_0^{\infty}\bar{F}_p(un^{-1})u^{k-1}(1 - un^{-1})^{n-k}\1_{u\leq n}\intd u\\%
    &\qquad= \int_0^{\infty} \lim_{n\to \infty} \frac{\bar{F}_p(un^{-1})}{n^{\alpha}\ell(n)}u^{k-1}(1 - un^{-1})^{n-k}\1_{u\leq n} \intd u\\
    &\qquad=\int_0^{\infty}u^{-\alpha + k - 1}e^{-u}\intd u\\
    &\qquad= \Gamma(k-\alpha).
  \end{align*}
  The conclusion follows by Stirling's formula.
\end{proof}

\subsection{Uniform results}
\label{sec:uniform-results}

\begin{lem}
  \label{lem:unif-expect-cumulfreq}
  For all $0 < \underline{\alpha} < \overline{\alpha} < 1$ and all
  $\overline{L} > 0$, and for all $k \geq 1$,
  \begin{equation*}
    \limsup_{n\to \infty}\sup_{\substack{(\alpha,L,L')\in \mathbb{A}\\p\in \mathcal{P}(\alpha,L,L')}}%
    \frac{\big|\EE_p[C_{n,k}] - L n^{\alpha}\frac{\Gamma(k-\alpha)}{\Gamma(k)} \big|}{L' n^{\alpha/2} \sqrt{\log(n)}\Gamma(k - \alpha/2)}%
    \leq \Gamma(k).
  \end{equation*}
\end{lem}
\begin{proof}
  If $p \in \mathcal{P}(\alpha,L,L')$, then we can rewrite
  $\bar{F}_p(x) = Lx^{-\alpha} + \delta(x)$ where $\delta$ is a function
  satisfying $|\delta(x)| \leq L' \sqrt{x^{-\alpha}\log(e/x)}$. Then, by the
  Lemma~\ref{lem:asympt-cumul-freq}
  \begin{align*}
    \bar{\EE}_p[C_{n,k}]%
    &= L k\binom{n}{k}\int_0^1x^{k-\alpha-1}(1-x)^{n-k}\intd x%
      + k\binom{n}{k}\int_0^1\delta(x)x^{k-1}(1-x)^{n-k}\intd x\\
    &= L \cdot \frac{n!}{\Gamma(n+1-\alpha)}\cdot \frac{\Gamma(k-\alpha)}{\Gamma(k)}%
      + k\binom{n}{k}\int_0^1\delta(x)x^{k-1}(1-x)^{n-k}\intd x.
  \end{align*}
  On one hand, by Stirling's formula, there exists $n_0$ such that for all
  $n \geq n_0$
  \begin{align*}
    \Big|\frac{L n!}{\Gamma(n+1-\alpha)} - Ln^{\alpha}  \Big|%
    &\leq 2Ln^{-1+\alpha}.
  \end{align*}
  On the other hand, since $|\delta(x)| \leq L' \sqrt{x^{-\alpha}\log(e/x)}$,
  \begin{equation*}
    \Big|\int_0^1\delta(x)x^{k-1}(1-x)^{n-k}\intd x \Big|%
    \leq L' \int_0^1 \sqrt{\log(e/x)}x^{-\alpha/2 + k - 1}(1-x)^{n-k}\intd x.
  \end{equation*}
  We remark that for all $B > 0$, the function
  $B\log(x) + \frac{1}{2}\log\log(e/x)$ has a unique maximum at
  $x = e^{1-1/(2B)}$. From this observation, we deduce that
  $\sqrt{\log(e/x)} \leq \frac{e^B}{e^{1/2}\sqrt{2B}}x^{-B}$ for all $B > 0$.
  Therefore taking $B > 0$ small enough, we find that
  \begin{align*}
    \Big|\int_0^1\delta(x)x^{k-1}(1-x)^{n-k}\intd x \Big|%
    &\leq \frac{L' e^B}{e^{1/2}\sqrt{2B}}\int_0^1 x^{-B - \alpha/2 + k -1}(1-x)^{n-k}\intd x\\
    &= \frac{L' e^B}{e^{1/2}\sqrt{2B}} \frac{\Gamma(n-k+1)\Gamma(k - \alpha/2 - B)}{\Gamma(n + 1 -\alpha/2 - B)}.
  \end{align*}
  By an asymptotic expansion of $\log \Gamma(z - B)$ as $B \to 0$ and using that
  $\frac{\Gamma'(z)}{\Gamma(z)} \sim \log(z)$ as $z \to \infty$, we find that
  \begin{align*}
    \limsup_{n\to\infty} \sup_{(\alpha,L,L')\in \mathbb{A}}%
    \frac{\big|k\binom{n}{k}\int_0^1\delta(x)x^{k-1}(1-x)^{n-k}\intd x \big|  }{L' n^{\alpha/2}\sqrt{\log(n)}\Gamma(k-\alpha/2)/\Gamma(k) }%
    &\leq \limsup_{n\to \infty} \frac{e^B\cdot e^{B \log(n)}}{e^{1/2}\sqrt{2B \cdot \log(n)}}\\
    &\leq 1,
  \end{align*}
  by choosing $2B\cdot \log(n) = 1$.
\end{proof}

\section{Additional results regarding simulations}
\label{sec:addit-inform-about}

\subsection{Proof of Lemma~\ref{lem:double-zipf}}
\label{sec:proof-lemma}

Recall that
$p = C(\alpha,\beta,J)\big[\sum_{j=1}^Jj^{-1/\alpha}\delta_j + \sum_{j>J}j^{-1/\beta}\delta_j \big]$,
with $\beta < \alpha$. Now suppose
$x < C(\alpha,\beta,J)J^{-1/\beta}$. We have
that
$C(\alpha,\beta,J)j^{-1/\alpha} \geq C(\alpha,\beta,J)j^{-1/\beta} \geq C(\alpha,\beta,J)J^{-1/\beta} > x$
for all $j\in \Set{1,\dots,J}$. Therefore when
$x < C(\alpha,\beta,J)J^{-1/\beta}$
\begin{equation*}
  \bar{F}_p(x)%
  = \Big\lfloor \frac{C(\alpha,\beta,J)^{\beta}}{x^{\beta}}\Big\rfloor
\end{equation*}
from which we deduce that
$\lim_{x\to 0} x^{\beta}\bar{F}_p(x) = C(\alpha,\beta,J)^{\beta}$ and using that
$ C(\alpha,\beta,J)^{\beta} - x^{\beta} \leq  x^{\beta}\lfloor C(\alpha,\beta,J)^{\beta}/x^{\beta} \rfloor \leq C(\alpha,\beta,J)^{\beta}$
\begin{equation*}
  \frac{|x^{\beta}\bar{F}_p(x) - C(\alpha,\beta,J)^{\beta}|}{\sqrt{x^{\beta}\log(e/x)}}%
  \leq \frac{x^{\beta}}{\sqrt{x^{\beta/2}\log(e/x)}}%
  \leq \frac{C(\alpha,\beta,J)^{\beta}}{\sqrt{J - J \log C(\alpha,\beta,J) + (J/\beta)\log J}}.
\end{equation*}
We now consider the case where $x \geq C(\alpha,\beta,J)J^{-1/\beta}$. In this
case we have that
\begin{equation*}
  \bar{F}_p(x)%
  = \Big\lfloor \frac{C(\alpha,\beta,J)^{\alpha}}{x^{\alpha}} \Big\rfloor.
\end{equation*}
With the same reasoning as above,
\begin{align*}
  |x^{\beta}\bar{F}_p(x) - C(\alpha,\beta,J)^{\beta}|%
  &\leq%
    x^{\beta} +  C(\alpha,\beta,J)^{\beta}\Big|\frac{[C(\alpha,\beta,J)/x]^{\alpha}}{[C(\alpha,\beta,J)/x]^{\beta}} - 1 \Big|.
\end{align*}
Observe that
$C(\alpha,\beta,J) \geq \sum_{j\geq 1}j^{-1/\beta} = \zeta(1/\beta)$ where
$\zeta$ is Riemann's Zeta function, because $\alpha > \beta$. But a classical
analysis of the Zeta function shows that $\zeta(1/\beta) > 1$ for all
$0 < \beta < 1$. It follows that for all $x \in [0,1]$ we must have
$x \leq C(\alpha,\beta,J)$ and we only need to analyze the last display for
$C(\alpha,\beta,J)J^{-1/\beta} \leq x \leq C(\alpha,\beta,J)$. In this case,
\begin{equation*}
  \frac{|x^{\beta}\bar{F}_p(x) - C(\alpha,\beta,J)^{\beta}|}{\sqrt{x^{\beta}\log(e/x)}}%
  \leq%
  \sqrt{\frac{x^{\beta}}{\log(e/x)}}%
  + C(\alpha,\beta,J)^{\alpha} \frac{x^{\beta-\alpha}(\alpha-\beta)\log \frac{C(\alpha,\beta,J)}{x}}{\sqrt{x^{\beta}\log(e/x)}}
\end{equation*}
which completes the proof.

\subsection{Elementary diagnoses for checking the power-law assumption}
\label{sec:elem-diagn-check}

As mentioned in the introduction of the paper, it is impossible to guarantee
that $p \in \mathcal{P}(\alpha,L,L')$ on the basis of the data alone. This is
because one cannot extrapolate the behaviour of $\bar{F}_p(x)$ near zero from a
finite amount of data. If that was possible, estimating the unseen would be
feasible without extra modeling assumption. Consequently, without exterior
information on $\bar{F}_p$, coming for instance from biological or physical
insights, the result of the estimation cannot be, in general, trusted. That is
said, power-law data have readily identifiable signatures that allows to assess
if the power-law assumption is plausible or not. In other word, we cannot
guarantee that the assumption holds, but it is often easy to diagnose when it
does not hold.

\begin{figure}[!htb]
  \centering

  \includegraphics[width=11cm]{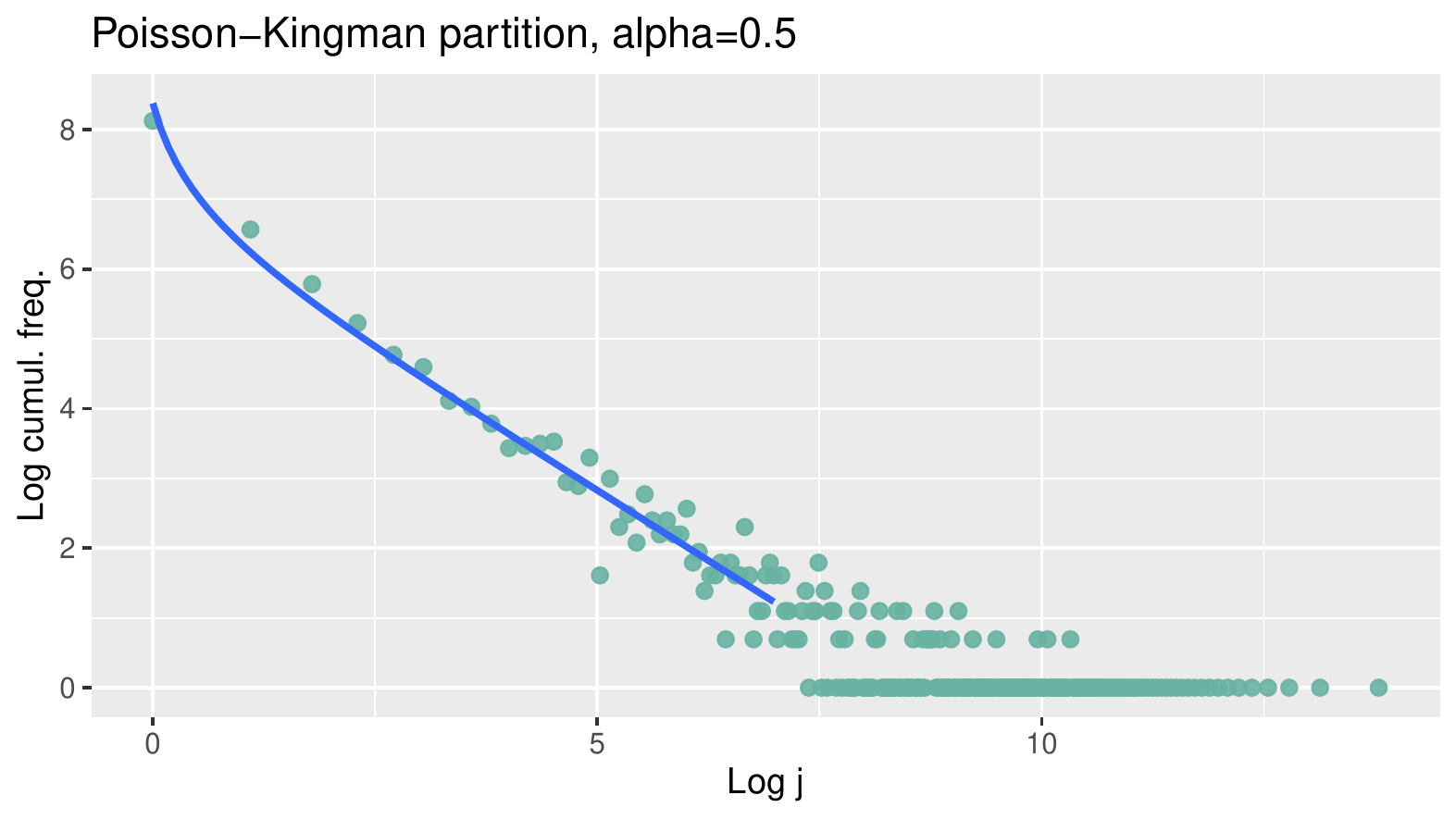}

  \caption{In green, $C_{n,j}$ as a function of $j$ on a log-log scale. The
    cumulative frequencies $C_{n,j}$ come from a Poisson-Kingman partition with
    Lévy measure $\rho_{\alpha}$, with $\alpha = 0.7$ and $n=1000000$. In blue
    is the best fit for
    $\log(C_{n,j}) = \log \frac{\Gamma(j-A)}{\Gamma(j)} + B$, which is the
    expected theoretical behaviour of the cumulative frequencies for power-law
    data.}
  \label{fig:diag1}
\end{figure}

\begin{figure}[!htb]
  \centering

  \includegraphics[width=15cm]{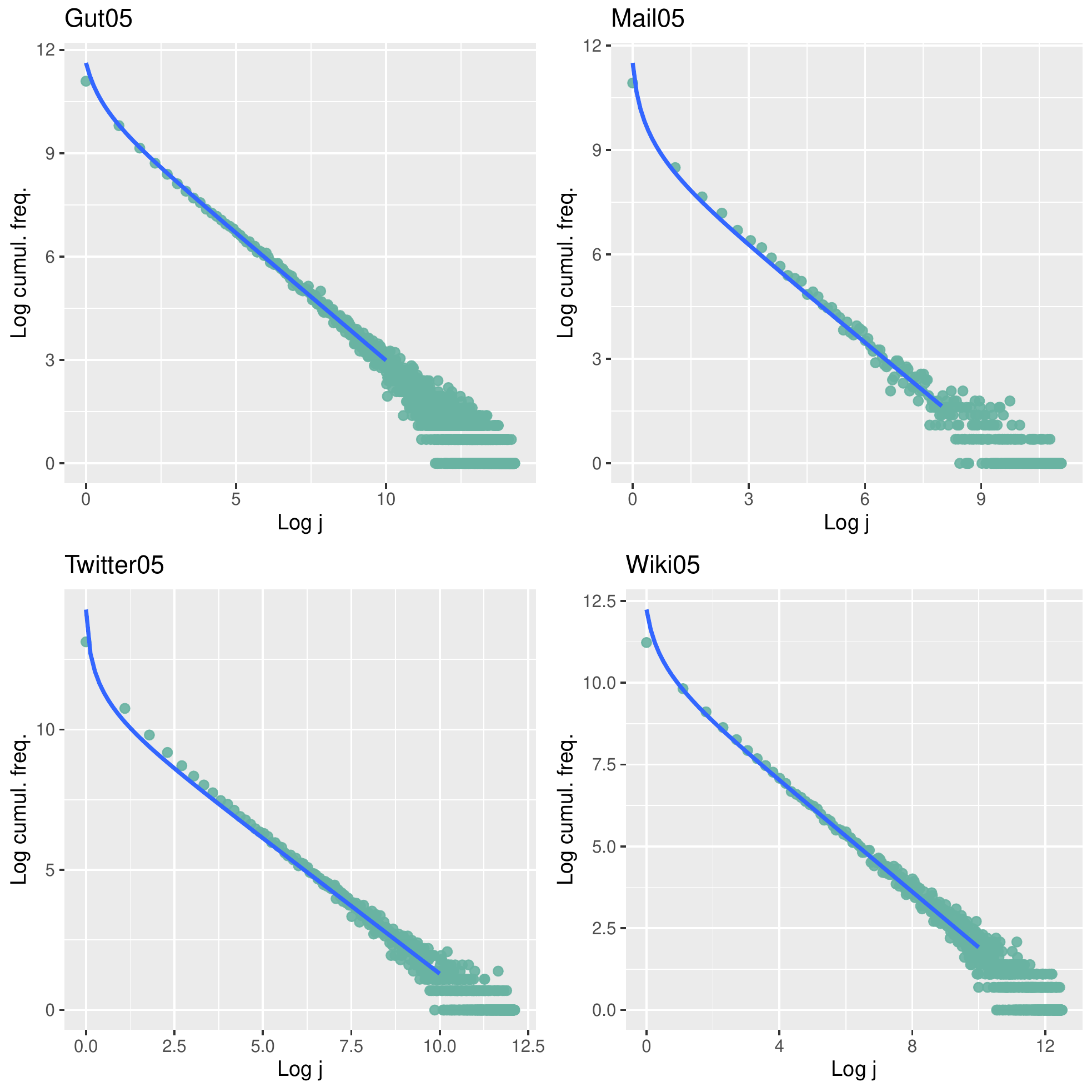}

  \caption{In green, $C_{n,j}$ as a function of $j$ on a log-log scale for the
    four datasets used in the numerical illustration section. In blue are the
    corresponding best fit for
    $\log(C_{n,j}) = \log \frac{\Gamma(j-A)}{\Gamma(j)} + B$, which is the
    expected theoretical behaviour of the cumulative frequencies for power-law
    data.}
  \label{fig:diag2}
\end{figure}

One possible diagnostic consists in inspecting the behaviour of the statistic
$j \mapsto C_{n,j}$. If the data is power-law, for small $j$ the variables
$C_{n,j}$ should not deviate too much from their expectation (in virtue of the
Lemma~\ref{lem:unif-expect-cumulfreq}). Then, we expect that
$C_{n,j} \propto \frac{\Gamma(j - \alpha)}{\Gamma(j)}$. To illustrate this, we
have plotted in Figure~\ref{fig:diag1} $\log(C_{n,j})$ as a function of
$\log(j)$ for a partition generated from a Poisson-Kingman process with Lévy
measure $\rho_{\alpha}$, with $\alpha=0.7$ and $n= 1000000$. We have added in
blue the best fit (with respect to mean square error) obtained for
$\log(C_{n,j}) = \log \frac{\Gamma(j-A)}{\Gamma(j)} + B$. This confirms
that in the ideal scenario, the empirical cumulative frequencies are in good
agreement with their expected values. In Figure~\ref{fig:diag2}, we did the same
thing for the four real datasets used in Section~\ref{sec:real-data}. In all
cases, the empirical cumulative frequencies seem compatible with the power-law
assumption.

\end{document}

%% file: unseen_large.pgf
\begingroup%
\makeatletter%
\begin{pgfpicture}%
\pgfpathrectangle{\pgfpointorigin}{\pgfqpoint{5.492558in}{4.246309in}}%
\pgfusepath{use as bounding box, clip}%
\begin{pgfscope}%
\pgfsetbuttcap%
\pgfsetmiterjoin%
\definecolor{currentfill}{rgb}{1.000000,1.000000,1.000000}%
\pgfsetfillcolor{currentfill}%
\pgfsetlinewidth{0.000000pt}%
\definecolor{currentstroke}{rgb}{1.000000,1.000000,1.000000}%
\pgfsetstrokecolor{currentstroke}%
\pgfsetdash{}{0pt}%
\pgfpathmoveto{\pgfqpoint{0.000000in}{-0.000000in}}%
\pgfpathlineto{\pgfqpoint{5.492558in}{-0.000000in}}%
\pgfpathlineto{\pgfqpoint{5.492558in}{4.246309in}}%
\pgfpathlineto{\pgfqpoint{0.000000in}{4.246309in}}%
\pgfpathclose%
\pgfusepath{fill}%
\end{pgfscope}%
\begin{pgfscope}%
\pgfsetbuttcap%
\pgfsetmiterjoin%
\definecolor{currentfill}{rgb}{1.000000,1.000000,1.000000}%
\pgfsetfillcolor{currentfill}%
\pgfsetlinewidth{0.000000pt}%
\definecolor{currentstroke}{rgb}{0.000000,0.000000,0.000000}%
\pgfsetstrokecolor{currentstroke}%
\pgfsetstrokeopacity{0.000000}%
\pgfsetdash{}{0pt}%
\pgfpathmoveto{\pgfqpoint{0.432558in}{0.450309in}}%
\pgfpathlineto{\pgfqpoint{5.392558in}{0.450309in}}%
\pgfpathlineto{\pgfqpoint{5.392558in}{4.146309in}}%
\pgfpathlineto{\pgfqpoint{0.432558in}{4.146309in}}%
\pgfpathclose%
\pgfusepath{fill}%
\end{pgfscope}%
\begin{pgfscope}%
\pgfsetbuttcap%
\pgfsetroundjoin%
\definecolor{currentfill}{rgb}{0.000000,0.000000,0.000000}%
\pgfsetfillcolor{currentfill}%
\pgfsetlinewidth{0.803000pt}%
\definecolor{currentstroke}{rgb}{0.000000,0.000000,0.000000}%
\pgfsetstrokecolor{currentstroke}%
\pgfsetdash{}{0pt}%
\pgfsys@defobject{currentmarker}{\pgfqpoint{0.000000in}{-0.048611in}}{\pgfqpoint{0.000000in}{0.000000in}}{%
\pgfpathmoveto{\pgfqpoint{0.000000in}{0.000000in}}%
\pgfpathlineto{\pgfqpoint{0.000000in}{-0.048611in}}%
\pgfusepath{stroke,fill}%
}%
\begin{pgfscope}%
\pgfsys@transformshift{0.612467in}{0.450309in}%
\pgfsys@useobject{currentmarker}{}%
\end{pgfscope}%
\end{pgfscope}%
\begin{pgfscope}%
\definecolor{textcolor}{rgb}{0.000000,0.000000,0.000000}%
\pgfsetstrokecolor{textcolor}%
\pgfsetfillcolor{textcolor}%
\pgftext[x=0.612467in,y=0.353086in,,top]{\color{textcolor}\rmfamily\fontsize{6.000000}{7.200000}\selectfont \(\displaystyle {0}\)}%
\end{pgfscope}%
\begin{pgfscope}%
\pgfsetbuttcap%
\pgfsetroundjoin%
\definecolor{currentfill}{rgb}{0.000000,0.000000,0.000000}%
\pgfsetfillcolor{currentfill}%
\pgfsetlinewidth{0.803000pt}%
\definecolor{currentstroke}{rgb}{0.000000,0.000000,0.000000}%
\pgfsetstrokecolor{currentstroke}%
\pgfsetdash{}{0pt}%
\pgfsys@defobject{currentmarker}{\pgfqpoint{0.000000in}{-0.048611in}}{\pgfqpoint{0.000000in}{0.000000in}}{%
\pgfpathmoveto{\pgfqpoint{0.000000in}{0.000000in}}%
\pgfpathlineto{\pgfqpoint{0.000000in}{-0.048611in}}%
\pgfusepath{stroke,fill}%
}%
\begin{pgfscope}%
\pgfsys@transformshift{1.523394in}{0.450309in}%
\pgfsys@useobject{currentmarker}{}%
\end{pgfscope}%
\end{pgfscope}%
\begin{pgfscope}%
\definecolor{textcolor}{rgb}{0.000000,0.000000,0.000000}%
\pgfsetstrokecolor{textcolor}%
\pgfsetfillcolor{textcolor}%
\pgftext[x=1.523394in,y=0.353086in,,top]{\color{textcolor}\rmfamily\fontsize{6.000000}{7.200000}\selectfont \(\displaystyle {2000}\)}%
\end{pgfscope}%
\begin{pgfscope}%
\pgfsetbuttcap%
\pgfsetroundjoin%
\definecolor{currentfill}{rgb}{0.000000,0.000000,0.000000}%
\pgfsetfillcolor{currentfill}%
\pgfsetlinewidth{0.803000pt}%
\definecolor{currentstroke}{rgb}{0.000000,0.000000,0.000000}%
\pgfsetstrokecolor{currentstroke}%
\pgfsetdash{}{0pt}%
\pgfsys@defobject{currentmarker}{\pgfqpoint{0.000000in}{-0.048611in}}{\pgfqpoint{0.000000in}{0.000000in}}{%
\pgfpathmoveto{\pgfqpoint{0.000000in}{0.000000in}}%
\pgfpathlineto{\pgfqpoint{0.000000in}{-0.048611in}}%
\pgfusepath{stroke,fill}%
}%
\begin{pgfscope}%
\pgfsys@transformshift{2.434321in}{0.450309in}%
\pgfsys@useobject{currentmarker}{}%
\end{pgfscope}%
\end{pgfscope}%
\begin{pgfscope}%
\definecolor{textcolor}{rgb}{0.000000,0.000000,0.000000}%
\pgfsetstrokecolor{textcolor}%
\pgfsetfillcolor{textcolor}%
\pgftext[x=2.434321in,y=0.353086in,,top]{\color{textcolor}\rmfamily\fontsize{6.000000}{7.200000}\selectfont \(\displaystyle {4000}\)}%
\end{pgfscope}%
\begin{pgfscope}%
\pgfsetbuttcap%
\pgfsetroundjoin%
\definecolor{currentfill}{rgb}{0.000000,0.000000,0.000000}%
\pgfsetfillcolor{currentfill}%
\pgfsetlinewidth{0.803000pt}%
\definecolor{currentstroke}{rgb}{0.000000,0.000000,0.000000}%
\pgfsetstrokecolor{currentstroke}%
\pgfsetdash{}{0pt}%
\pgfsys@defobject{currentmarker}{\pgfqpoint{0.000000in}{-0.048611in}}{\pgfqpoint{0.000000in}{0.000000in}}{%
\pgfpathmoveto{\pgfqpoint{0.000000in}{0.000000in}}%
\pgfpathlineto{\pgfqpoint{0.000000in}{-0.048611in}}%
\pgfusepath{stroke,fill}%
}%
\begin{pgfscope}%
\pgfsys@transformshift{3.345249in}{0.450309in}%
\pgfsys@useobject{currentmarker}{}%
\end{pgfscope}%
\end{pgfscope}%
\begin{pgfscope}%
\definecolor{textcolor}{rgb}{0.000000,0.000000,0.000000}%
\pgfsetstrokecolor{textcolor}%
\pgfsetfillcolor{textcolor}%
\pgftext[x=3.345249in,y=0.353086in,,top]{\color{textcolor}\rmfamily\fontsize{6.000000}{7.200000}\selectfont \(\displaystyle {6000}\)}%
\end{pgfscope}%
\begin{pgfscope}%
\pgfsetbuttcap%
\pgfsetroundjoin%
\definecolor{currentfill}{rgb}{0.000000,0.000000,0.000000}%
\pgfsetfillcolor{currentfill}%
\pgfsetlinewidth{0.803000pt}%
\definecolor{currentstroke}{rgb}{0.000000,0.000000,0.000000}%
\pgfsetstrokecolor{currentstroke}%
\pgfsetdash{}{0pt}%
\pgfsys@defobject{currentmarker}{\pgfqpoint{0.000000in}{-0.048611in}}{\pgfqpoint{0.000000in}{0.000000in}}{%
\pgfpathmoveto{\pgfqpoint{0.000000in}{0.000000in}}%
\pgfpathlineto{\pgfqpoint{0.000000in}{-0.048611in}}%
\pgfusepath{stroke,fill}%
}%
\begin{pgfscope}%
\pgfsys@transformshift{4.256176in}{0.450309in}%
\pgfsys@useobject{currentmarker}{}%
\end{pgfscope}%
\end{pgfscope}%
\begin{pgfscope}%
\definecolor{textcolor}{rgb}{0.000000,0.000000,0.000000}%
\pgfsetstrokecolor{textcolor}%
\pgfsetfillcolor{textcolor}%
\pgftext[x=4.256176in,y=0.353086in,,top]{\color{textcolor}\rmfamily\fontsize{6.000000}{7.200000}\selectfont \(\displaystyle {8000}\)}%
\end{pgfscope}%
\begin{pgfscope}%
\pgfsetbuttcap%
\pgfsetroundjoin%
\definecolor{currentfill}{rgb}{0.000000,0.000000,0.000000}%
\pgfsetfillcolor{currentfill}%
\pgfsetlinewidth{0.803000pt}%
\definecolor{currentstroke}{rgb}{0.000000,0.000000,0.000000}%
\pgfsetstrokecolor{currentstroke}%
\pgfsetdash{}{0pt}%
\pgfsys@defobject{currentmarker}{\pgfqpoint{0.000000in}{-0.048611in}}{\pgfqpoint{0.000000in}{0.000000in}}{%
\pgfpathmoveto{\pgfqpoint{0.000000in}{0.000000in}}%
\pgfpathlineto{\pgfqpoint{0.000000in}{-0.048611in}}%
\pgfusepath{stroke,fill}%
}%
\begin{pgfscope}%
\pgfsys@transformshift{5.167104in}{0.450309in}%
\pgfsys@useobject{currentmarker}{}%
\end{pgfscope}%
\end{pgfscope}%
\begin{pgfscope}%
\definecolor{textcolor}{rgb}{0.000000,0.000000,0.000000}%
\pgfsetstrokecolor{textcolor}%
\pgfsetfillcolor{textcolor}%
\pgftext[x=5.167104in,y=0.353086in,,top]{\color{textcolor}\rmfamily\fontsize{6.000000}{7.200000}\selectfont \(\displaystyle {10000}\)}%
\end{pgfscope}%
\begin{pgfscope}%
\definecolor{textcolor}{rgb}{0.000000,0.000000,0.000000}%
\pgfsetstrokecolor{textcolor}%
\pgfsetfillcolor{textcolor}%
\pgftext[x=2.912558in,y=0.223457in,,top]{\color{textcolor}\rmfamily\fontsize{10.000000}{12.000000}\selectfont \(\displaystyle \lambda\)}%
\end{pgfscope}%
\begin{pgfscope}%
\pgfsetbuttcap%
\pgfsetroundjoin%
\definecolor{currentfill}{rgb}{0.000000,0.000000,0.000000}%
\pgfsetfillcolor{currentfill}%
\pgfsetlinewidth{0.803000pt}%
\definecolor{currentstroke}{rgb}{0.000000,0.000000,0.000000}%
\pgfsetstrokecolor{currentstroke}%
\pgfsetdash{}{0pt}%
\pgfsys@defobject{currentmarker}{\pgfqpoint{-0.048611in}{0.000000in}}{\pgfqpoint{-0.000000in}{0.000000in}}{%
\pgfpathmoveto{\pgfqpoint{-0.000000in}{0.000000in}}%
\pgfpathlineto{\pgfqpoint{-0.048611in}{0.000000in}}%
\pgfusepath{stroke,fill}%
}%
\begin{pgfscope}%
\pgfsys@transformshift{0.432558in}{1.207884in}%
\pgfsys@useobject{currentmarker}{}%
\end{pgfscope}%
\end{pgfscope}%
\begin{pgfscope}%
\definecolor{textcolor}{rgb}{0.000000,0.000000,0.000000}%
\pgfsetstrokecolor{textcolor}%
\pgfsetfillcolor{textcolor}%
\pgftext[x=0.100000in, y=1.178949in, left, base]{\color{textcolor}\rmfamily\fontsize{6.000000}{7.200000}\selectfont \(\displaystyle {0.015}\)}%
\end{pgfscope}%
\begin{pgfscope}%
\pgfsetbuttcap%
\pgfsetroundjoin%
\definecolor{currentfill}{rgb}{0.000000,0.000000,0.000000}%
\pgfsetfillcolor{currentfill}%
\pgfsetlinewidth{0.803000pt}%
\definecolor{currentstroke}{rgb}{0.000000,0.000000,0.000000}%
\pgfsetstrokecolor{currentstroke}%
\pgfsetdash{}{0pt}%
\pgfsys@defobject{currentmarker}{\pgfqpoint{-0.048611in}{0.000000in}}{\pgfqpoint{-0.000000in}{0.000000in}}{%
\pgfpathmoveto{\pgfqpoint{-0.000000in}{0.000000in}}%
\pgfpathlineto{\pgfqpoint{-0.048611in}{0.000000in}}%
\pgfusepath{stroke,fill}%
}%
\begin{pgfscope}%
\pgfsys@transformshift{0.432558in}{2.021900in}%
\pgfsys@useobject{currentmarker}{}%
\end{pgfscope}%
\end{pgfscope}%
\begin{pgfscope}%
\definecolor{textcolor}{rgb}{0.000000,0.000000,0.000000}%
\pgfsetstrokecolor{textcolor}%
\pgfsetfillcolor{textcolor}%
\pgftext[x=0.100000in, y=1.992965in, left, base]{\color{textcolor}\rmfamily\fontsize{6.000000}{7.200000}\selectfont \(\displaystyle {0.020}\)}%
\end{pgfscope}%
\begin{pgfscope}%
\pgfsetbuttcap%
\pgfsetroundjoin%
\definecolor{currentfill}{rgb}{0.000000,0.000000,0.000000}%
\pgfsetfillcolor{currentfill}%
\pgfsetlinewidth{0.803000pt}%
\definecolor{currentstroke}{rgb}{0.000000,0.000000,0.000000}%
\pgfsetstrokecolor{currentstroke}%
\pgfsetdash{}{0pt}%
\pgfsys@defobject{currentmarker}{\pgfqpoint{-0.048611in}{0.000000in}}{\pgfqpoint{-0.000000in}{0.000000in}}{%
\pgfpathmoveto{\pgfqpoint{-0.000000in}{0.000000in}}%
\pgfpathlineto{\pgfqpoint{-0.048611in}{0.000000in}}%
\pgfusepath{stroke,fill}%
}%
\begin{pgfscope}%
\pgfsys@transformshift{0.432558in}{2.835917in}%
\pgfsys@useobject{currentmarker}{}%
\end{pgfscope}%
\end{pgfscope}%
\begin{pgfscope}%
\definecolor{textcolor}{rgb}{0.000000,0.000000,0.000000}%
\pgfsetstrokecolor{textcolor}%
\pgfsetfillcolor{textcolor}%
\pgftext[x=0.100000in, y=2.806981in, left, base]{\color{textcolor}\rmfamily\fontsize{6.000000}{7.200000}\selectfont \(\displaystyle {0.025}\)}%
\end{pgfscope}%
\begin{pgfscope}%
\pgfsetbuttcap%
\pgfsetroundjoin%
\definecolor{currentfill}{rgb}{0.000000,0.000000,0.000000}%
\pgfsetfillcolor{currentfill}%
\pgfsetlinewidth{0.803000pt}%
\definecolor{currentstroke}{rgb}{0.000000,0.000000,0.000000}%
\pgfsetstrokecolor{currentstroke}%
\pgfsetdash{}{0pt}%
\pgfsys@defobject{currentmarker}{\pgfqpoint{-0.048611in}{0.000000in}}{\pgfqpoint{-0.000000in}{0.000000in}}{%
\pgfpathmoveto{\pgfqpoint{-0.000000in}{0.000000in}}%
\pgfpathlineto{\pgfqpoint{-0.048611in}{0.000000in}}%
\pgfusepath{stroke,fill}%
}%
\begin{pgfscope}%
\pgfsys@transformshift{0.432558in}{3.649933in}%
\pgfsys@useobject{currentmarker}{}%
\end{pgfscope}%
\end{pgfscope}%
\begin{pgfscope}%
\definecolor{textcolor}{rgb}{0.000000,0.000000,0.000000}%
\pgfsetstrokecolor{textcolor}%
\pgfsetfillcolor{textcolor}%
\pgftext[x=0.100000in, y=3.620998in, left, base]{\color{textcolor}\rmfamily\fontsize{6.000000}{7.200000}\selectfont \(\displaystyle {0.030}\)}%
\end{pgfscope}%
\begin{pgfscope}%
\pgfpathrectangle{\pgfqpoint{0.432558in}{0.450309in}}{\pgfqpoint{4.960000in}{3.696000in}}%
\pgfusepath{clip}%
\pgfsetrectcap%
\pgfsetroundjoin%
\pgfsetlinewidth{1.505625pt}%
\definecolor{currentstroke}{rgb}{1.000000,0.000000,0.000000}%
\pgfsetstrokecolor{currentstroke}%
\pgfsetdash{}{0pt}%
\pgfpathmoveto{\pgfqpoint{0.658013in}{0.618309in}}%
\pgfpathlineto{\pgfqpoint{0.703559in}{0.987154in}}%
\pgfpathlineto{\pgfqpoint{0.840198in}{1.521893in}}%
\pgfpathlineto{\pgfqpoint{1.067930in}{1.992371in}}%
\pgfpathlineto{\pgfqpoint{1.523394in}{2.526898in}}%
\pgfpathlineto{\pgfqpoint{1.978858in}{2.870888in}}%
\pgfpathlineto{\pgfqpoint{2.434321in}{3.121867in}}%
\pgfpathlineto{\pgfqpoint{2.889785in}{3.321196in}}%
\pgfpathlineto{\pgfqpoint{3.345249in}{3.493334in}}%
\pgfpathlineto{\pgfqpoint{3.800713in}{3.633392in}}%
\pgfpathlineto{\pgfqpoint{4.256176in}{3.760226in}}%
\pgfpathlineto{\pgfqpoint{4.711640in}{3.874950in}}%
\pgfpathlineto{\pgfqpoint{5.167104in}{3.978309in}}%
\pgfusepath{stroke}%
\end{pgfscope}%
\begin{pgfscope}%
\pgfsetrectcap%
\pgfsetmiterjoin%
\pgfsetlinewidth{0.803000pt}%
\definecolor{currentstroke}{rgb}{0.000000,0.000000,0.000000}%
\pgfsetstrokecolor{currentstroke}%
\pgfsetdash{}{0pt}%
\pgfpathmoveto{\pgfqpoint{0.432558in}{0.450309in}}%
\pgfpathlineto{\pgfqpoint{0.432558in}{4.146309in}}%
\pgfusepath{stroke}%
\end{pgfscope}%
\begin{pgfscope}%
\pgfsetrectcap%
\pgfsetmiterjoin%
\pgfsetlinewidth{0.803000pt}%
\definecolor{currentstroke}{rgb}{0.000000,0.000000,0.000000}%
\pgfsetstrokecolor{currentstroke}%
\pgfsetdash{}{0pt}%
\pgfpathmoveto{\pgfqpoint{5.392558in}{0.450309in}}%
\pgfpathlineto{\pgfqpoint{5.392558in}{4.146309in}}%
\pgfusepath{stroke}%
\end{pgfscope}%
\begin{pgfscope}%
\pgfsetrectcap%
\pgfsetmiterjoin%
\pgfsetlinewidth{0.803000pt}%
\definecolor{currentstroke}{rgb}{0.000000,0.000000,0.000000}%
\pgfsetstrokecolor{currentstroke}%
\pgfsetdash{}{0pt}%
\pgfpathmoveto{\pgfqpoint{0.432558in}{0.450309in}}%
\pgfpathlineto{\pgfqpoint{5.392558in}{0.450309in}}%
\pgfusepath{stroke}%
\end{pgfscope}%
\begin{pgfscope}%
\pgfsetrectcap%
\pgfsetmiterjoin%
\pgfsetlinewidth{0.803000pt}%
\definecolor{currentstroke}{rgb}{0.000000,0.000000,0.000000}%
\pgfsetstrokecolor{currentstroke}%
\pgfsetdash{}{0pt}%
\pgfpathmoveto{\pgfqpoint{0.432558in}{4.146309in}}%
\pgfpathlineto{\pgfqpoint{5.392558in}{4.146309in}}%
\pgfusepath{stroke}%
\end{pgfscope}%
\end{pgfpicture}%
\makeatother%
\endgroup%

%% file: mc_estimate_risk_alpha_double_zipf.pgf
\begingroup%
\makeatletter%
\begin{pgfpicture}%
\pgfpathrectangle{\pgfpointorigin}{\pgfqpoint{5.613333in}{4.116679in}}%
\pgfusepath{use as bounding box, clip}%
\begin{pgfscope}%
\pgfsetbuttcap%
\pgfsetmiterjoin%
\definecolor{currentfill}{rgb}{1.000000,1.000000,1.000000}%
\pgfsetfillcolor{currentfill}%
\pgfsetlinewidth{0.000000pt}%
\definecolor{currentstroke}{rgb}{1.000000,1.000000,1.000000}%
\pgfsetstrokecolor{currentstroke}%
\pgfsetdash{}{0pt}%
\pgfpathmoveto{\pgfqpoint{0.000000in}{0.000000in}}%
\pgfpathlineto{\pgfqpoint{5.613333in}{0.000000in}}%
\pgfpathlineto{\pgfqpoint{5.613333in}{4.116679in}}%
\pgfpathlineto{\pgfqpoint{0.000000in}{4.116679in}}%
\pgfpathclose%
\pgfusepath{fill}%
\end{pgfscope}%
\begin{pgfscope}%
\pgfsetbuttcap%
\pgfsetmiterjoin%
\definecolor{currentfill}{rgb}{1.000000,1.000000,1.000000}%
\pgfsetfillcolor{currentfill}%
\pgfsetlinewidth{0.000000pt}%
\definecolor{currentstroke}{rgb}{0.000000,0.000000,0.000000}%
\pgfsetstrokecolor{currentstroke}%
\pgfsetstrokeopacity{0.000000}%
\pgfsetdash{}{0pt}%
\pgfpathmoveto{\pgfqpoint{0.513581in}{0.320679in}}%
\pgfpathlineto{\pgfqpoint{5.473581in}{0.320679in}}%
\pgfpathlineto{\pgfqpoint{5.473581in}{4.016679in}}%
\pgfpathlineto{\pgfqpoint{0.513581in}{4.016679in}}%
\pgfpathclose%
\pgfusepath{fill}%
\end{pgfscope}%
\begin{pgfscope}%
\pgfsetbuttcap%
\pgfsetroundjoin%
\definecolor{currentfill}{rgb}{0.000000,0.000000,0.000000}%
\pgfsetfillcolor{currentfill}%
\pgfsetlinewidth{0.803000pt}%
\definecolor{currentstroke}{rgb}{0.000000,0.000000,0.000000}%
\pgfsetstrokecolor{currentstroke}%
\pgfsetdash{}{0pt}%
\pgfsys@defobject{currentmarker}{\pgfqpoint{0.000000in}{-0.048611in}}{\pgfqpoint{0.000000in}{0.000000in}}{%
\pgfpathmoveto{\pgfqpoint{0.000000in}{0.000000in}}%
\pgfpathlineto{\pgfqpoint{0.000000in}{-0.048611in}}%
\pgfusepath{stroke,fill}%
}%
\begin{pgfscope}%
\pgfsys@transformshift{0.577997in}{0.320679in}%
\pgfsys@useobject{currentmarker}{}%
\end{pgfscope}%
\end{pgfscope}%
\begin{pgfscope}%
\definecolor{textcolor}{rgb}{0.000000,0.000000,0.000000}%
\pgfsetstrokecolor{textcolor}%
\pgfsetfillcolor{textcolor}%
\pgftext[x=0.577997in,y=0.223457in,,top]{\color{textcolor}\rmfamily\fontsize{10.000000}{12.000000}\selectfont \(\displaystyle {0}\)}%
\end{pgfscope}%
\begin{pgfscope}%
\pgfsetbuttcap%
\pgfsetroundjoin%
\definecolor{currentfill}{rgb}{0.000000,0.000000,0.000000}%
\pgfsetfillcolor{currentfill}%
\pgfsetlinewidth{0.803000pt}%
\definecolor{currentstroke}{rgb}{0.000000,0.000000,0.000000}%
\pgfsetstrokecolor{currentstroke}%
\pgfsetdash{}{0pt}%
\pgfsys@defobject{currentmarker}{\pgfqpoint{0.000000in}{-0.048611in}}{\pgfqpoint{0.000000in}{0.000000in}}{%
\pgfpathmoveto{\pgfqpoint{0.000000in}{0.000000in}}%
\pgfpathlineto{\pgfqpoint{0.000000in}{-0.048611in}}%
\pgfusepath{stroke,fill}%
}%
\begin{pgfscope}%
\pgfsys@transformshift{1.383192in}{0.320679in}%
\pgfsys@useobject{currentmarker}{}%
\end{pgfscope}%
\end{pgfscope}%
\begin{pgfscope}%
\definecolor{textcolor}{rgb}{0.000000,0.000000,0.000000}%
\pgfsetstrokecolor{textcolor}%
\pgfsetfillcolor{textcolor}%
\pgftext[x=1.383192in,y=0.223457in,,top]{\color{textcolor}\rmfamily\fontsize{10.000000}{12.000000}\selectfont \(\displaystyle {50}\)}%
\end{pgfscope}%
\begin{pgfscope}%
\pgfsetbuttcap%
\pgfsetroundjoin%
\definecolor{currentfill}{rgb}{0.000000,0.000000,0.000000}%
\pgfsetfillcolor{currentfill}%
\pgfsetlinewidth{0.803000pt}%
\definecolor{currentstroke}{rgb}{0.000000,0.000000,0.000000}%
\pgfsetstrokecolor{currentstroke}%
\pgfsetdash{}{0pt}%
\pgfsys@defobject{currentmarker}{\pgfqpoint{0.000000in}{-0.048611in}}{\pgfqpoint{0.000000in}{0.000000in}}{%
\pgfpathmoveto{\pgfqpoint{0.000000in}{0.000000in}}%
\pgfpathlineto{\pgfqpoint{0.000000in}{-0.048611in}}%
\pgfusepath{stroke,fill}%
}%
\begin{pgfscope}%
\pgfsys@transformshift{2.188386in}{0.320679in}%
\pgfsys@useobject{currentmarker}{}%
\end{pgfscope}%
\end{pgfscope}%
\begin{pgfscope}%
\definecolor{textcolor}{rgb}{0.000000,0.000000,0.000000}%
\pgfsetstrokecolor{textcolor}%
\pgfsetfillcolor{textcolor}%
\pgftext[x=2.188386in,y=0.223457in,,top]{\color{textcolor}\rmfamily\fontsize{10.000000}{12.000000}\selectfont \(\displaystyle {100}\)}%
\end{pgfscope}%
\begin{pgfscope}%
\pgfsetbuttcap%
\pgfsetroundjoin%
\definecolor{currentfill}{rgb}{0.000000,0.000000,0.000000}%
\pgfsetfillcolor{currentfill}%
\pgfsetlinewidth{0.803000pt}%
\definecolor{currentstroke}{rgb}{0.000000,0.000000,0.000000}%
\pgfsetstrokecolor{currentstroke}%
\pgfsetdash{}{0pt}%
\pgfsys@defobject{currentmarker}{\pgfqpoint{0.000000in}{-0.048611in}}{\pgfqpoint{0.000000in}{0.000000in}}{%
\pgfpathmoveto{\pgfqpoint{0.000000in}{0.000000in}}%
\pgfpathlineto{\pgfqpoint{0.000000in}{-0.048611in}}%
\pgfusepath{stroke,fill}%
}%
\begin{pgfscope}%
\pgfsys@transformshift{2.993581in}{0.320679in}%
\pgfsys@useobject{currentmarker}{}%
\end{pgfscope}%
\end{pgfscope}%
\begin{pgfscope}%
\definecolor{textcolor}{rgb}{0.000000,0.000000,0.000000}%
\pgfsetstrokecolor{textcolor}%
\pgfsetfillcolor{textcolor}%
\pgftext[x=2.993581in,y=0.223457in,,top]{\color{textcolor}\rmfamily\fontsize{10.000000}{12.000000}\selectfont \(\displaystyle {150}\)}%
\end{pgfscope}%
\begin{pgfscope}%
\pgfsetbuttcap%
\pgfsetroundjoin%
\definecolor{currentfill}{rgb}{0.000000,0.000000,0.000000}%
\pgfsetfillcolor{currentfill}%
\pgfsetlinewidth{0.803000pt}%
\definecolor{currentstroke}{rgb}{0.000000,0.000000,0.000000}%
\pgfsetstrokecolor{currentstroke}%
\pgfsetdash{}{0pt}%
\pgfsys@defobject{currentmarker}{\pgfqpoint{0.000000in}{-0.048611in}}{\pgfqpoint{0.000000in}{0.000000in}}{%
\pgfpathmoveto{\pgfqpoint{0.000000in}{0.000000in}}%
\pgfpathlineto{\pgfqpoint{0.000000in}{-0.048611in}}%
\pgfusepath{stroke,fill}%
}%
\begin{pgfscope}%
\pgfsys@transformshift{3.798776in}{0.320679in}%
\pgfsys@useobject{currentmarker}{}%
\end{pgfscope}%
\end{pgfscope}%
\begin{pgfscope}%
\definecolor{textcolor}{rgb}{0.000000,0.000000,0.000000}%
\pgfsetstrokecolor{textcolor}%
\pgfsetfillcolor{textcolor}%
\pgftext[x=3.798776in,y=0.223457in,,top]{\color{textcolor}\rmfamily\fontsize{10.000000}{12.000000}\selectfont \(\displaystyle {200}\)}%
\end{pgfscope}%
\begin{pgfscope}%
\pgfsetbuttcap%
\pgfsetroundjoin%
\definecolor{currentfill}{rgb}{0.000000,0.000000,0.000000}%
\pgfsetfillcolor{currentfill}%
\pgfsetlinewidth{0.803000pt}%
\definecolor{currentstroke}{rgb}{0.000000,0.000000,0.000000}%
\pgfsetstrokecolor{currentstroke}%
\pgfsetdash{}{0pt}%
\pgfsys@defobject{currentmarker}{\pgfqpoint{0.000000in}{-0.048611in}}{\pgfqpoint{0.000000in}{0.000000in}}{%
\pgfpathmoveto{\pgfqpoint{0.000000in}{0.000000in}}%
\pgfpathlineto{\pgfqpoint{0.000000in}{-0.048611in}}%
\pgfusepath{stroke,fill}%
}%
\begin{pgfscope}%
\pgfsys@transformshift{4.603971in}{0.320679in}%
\pgfsys@useobject{currentmarker}{}%
\end{pgfscope}%
\end{pgfscope}%
\begin{pgfscope}%
\definecolor{textcolor}{rgb}{0.000000,0.000000,0.000000}%
\pgfsetstrokecolor{textcolor}%
\pgfsetfillcolor{textcolor}%
\pgftext[x=4.603971in,y=0.223457in,,top]{\color{textcolor}\rmfamily\fontsize{10.000000}{12.000000}\selectfont \(\displaystyle {250}\)}%
\end{pgfscope}%
\begin{pgfscope}%
\pgfsetbuttcap%
\pgfsetroundjoin%
\definecolor{currentfill}{rgb}{0.000000,0.000000,0.000000}%
\pgfsetfillcolor{currentfill}%
\pgfsetlinewidth{0.803000pt}%
\definecolor{currentstroke}{rgb}{0.000000,0.000000,0.000000}%
\pgfsetstrokecolor{currentstroke}%
\pgfsetdash{}{0pt}%
\pgfsys@defobject{currentmarker}{\pgfqpoint{0.000000in}{-0.048611in}}{\pgfqpoint{0.000000in}{0.000000in}}{%
\pgfpathmoveto{\pgfqpoint{0.000000in}{0.000000in}}%
\pgfpathlineto{\pgfqpoint{0.000000in}{-0.048611in}}%
\pgfusepath{stroke,fill}%
}%
\begin{pgfscope}%
\pgfsys@transformshift{5.409166in}{0.320679in}%
\pgfsys@useobject{currentmarker}{}%
\end{pgfscope}%
\end{pgfscope}%
\begin{pgfscope}%
\definecolor{textcolor}{rgb}{0.000000,0.000000,0.000000}%
\pgfsetstrokecolor{textcolor}%
\pgfsetfillcolor{textcolor}%
\pgftext[x=5.409166in,y=0.223457in,,top]{\color{textcolor}\rmfamily\fontsize{10.000000}{12.000000}\selectfont \(\displaystyle {300}\)}%
\end{pgfscope}%
\begin{pgfscope}%
\pgfsetbuttcap%
\pgfsetroundjoin%
\definecolor{currentfill}{rgb}{0.000000,0.000000,0.000000}%
\pgfsetfillcolor{currentfill}%
\pgfsetlinewidth{0.803000pt}%
\definecolor{currentstroke}{rgb}{0.000000,0.000000,0.000000}%
\pgfsetstrokecolor{currentstroke}%
\pgfsetdash{}{0pt}%
\pgfsys@defobject{currentmarker}{\pgfqpoint{-0.048611in}{0.000000in}}{\pgfqpoint{-0.000000in}{0.000000in}}{%
\pgfpathmoveto{\pgfqpoint{-0.000000in}{0.000000in}}%
\pgfpathlineto{\pgfqpoint{-0.048611in}{0.000000in}}%
\pgfusepath{stroke,fill}%
}%
\begin{pgfscope}%
\pgfsys@transformshift{0.513581in}{0.885301in}%
\pgfsys@useobject{currentmarker}{}%
\end{pgfscope}%
\end{pgfscope}%
\begin{pgfscope}%
\definecolor{textcolor}{rgb}{0.000000,0.000000,0.000000}%
\pgfsetstrokecolor{textcolor}%
\pgfsetfillcolor{textcolor}%
\pgftext[x=0.100000in, y=0.837076in, left, base]{\color{textcolor}\rmfamily\fontsize{10.000000}{12.000000}\selectfont \(\displaystyle {0.002}\)}%
\end{pgfscope}%
\begin{pgfscope}%
\pgfsetbuttcap%
\pgfsetroundjoin%
\definecolor{currentfill}{rgb}{0.000000,0.000000,0.000000}%
\pgfsetfillcolor{currentfill}%
\pgfsetlinewidth{0.803000pt}%
\definecolor{currentstroke}{rgb}{0.000000,0.000000,0.000000}%
\pgfsetstrokecolor{currentstroke}%
\pgfsetdash{}{0pt}%
\pgfsys@defobject{currentmarker}{\pgfqpoint{-0.048611in}{0.000000in}}{\pgfqpoint{-0.000000in}{0.000000in}}{%
\pgfpathmoveto{\pgfqpoint{-0.000000in}{0.000000in}}%
\pgfpathlineto{\pgfqpoint{-0.048611in}{0.000000in}}%
\pgfusepath{stroke,fill}%
}%
\begin{pgfscope}%
\pgfsys@transformshift{0.513581in}{1.606532in}%
\pgfsys@useobject{currentmarker}{}%
\end{pgfscope}%
\end{pgfscope}%
\begin{pgfscope}%
\definecolor{textcolor}{rgb}{0.000000,0.000000,0.000000}%
\pgfsetstrokecolor{textcolor}%
\pgfsetfillcolor{textcolor}%
\pgftext[x=0.100000in, y=1.558307in, left, base]{\color{textcolor}\rmfamily\fontsize{10.000000}{12.000000}\selectfont \(\displaystyle {0.004}\)}%
\end{pgfscope}%
\begin{pgfscope}%
\pgfsetbuttcap%
\pgfsetroundjoin%
\definecolor{currentfill}{rgb}{0.000000,0.000000,0.000000}%
\pgfsetfillcolor{currentfill}%
\pgfsetlinewidth{0.803000pt}%
\definecolor{currentstroke}{rgb}{0.000000,0.000000,0.000000}%
\pgfsetstrokecolor{currentstroke}%
\pgfsetdash{}{0pt}%
\pgfsys@defobject{currentmarker}{\pgfqpoint{-0.048611in}{0.000000in}}{\pgfqpoint{-0.000000in}{0.000000in}}{%
\pgfpathmoveto{\pgfqpoint{-0.000000in}{0.000000in}}%
\pgfpathlineto{\pgfqpoint{-0.048611in}{0.000000in}}%
\pgfusepath{stroke,fill}%
}%
\begin{pgfscope}%
\pgfsys@transformshift{0.513581in}{2.327763in}%
\pgfsys@useobject{currentmarker}{}%
\end{pgfscope}%
\end{pgfscope}%
\begin{pgfscope}%
\definecolor{textcolor}{rgb}{0.000000,0.000000,0.000000}%
\pgfsetstrokecolor{textcolor}%
\pgfsetfillcolor{textcolor}%
\pgftext[x=0.100000in, y=2.279538in, left, base]{\color{textcolor}\rmfamily\fontsize{10.000000}{12.000000}\selectfont \(\displaystyle {0.006}\)}%
\end{pgfscope}%
\begin{pgfscope}%
\pgfsetbuttcap%
\pgfsetroundjoin%
\definecolor{currentfill}{rgb}{0.000000,0.000000,0.000000}%
\pgfsetfillcolor{currentfill}%
\pgfsetlinewidth{0.803000pt}%
\definecolor{currentstroke}{rgb}{0.000000,0.000000,0.000000}%
\pgfsetstrokecolor{currentstroke}%
\pgfsetdash{}{0pt}%
\pgfsys@defobject{currentmarker}{\pgfqpoint{-0.048611in}{0.000000in}}{\pgfqpoint{-0.000000in}{0.000000in}}{%
\pgfpathmoveto{\pgfqpoint{-0.000000in}{0.000000in}}%
\pgfpathlineto{\pgfqpoint{-0.048611in}{0.000000in}}%
\pgfusepath{stroke,fill}%
}%
\begin{pgfscope}%
\pgfsys@transformshift{0.513581in}{3.048994in}%
\pgfsys@useobject{currentmarker}{}%
\end{pgfscope}%
\end{pgfscope}%
\begin{pgfscope}%
\definecolor{textcolor}{rgb}{0.000000,0.000000,0.000000}%
\pgfsetstrokecolor{textcolor}%
\pgfsetfillcolor{textcolor}%
\pgftext[x=0.100000in, y=3.000769in, left, base]{\color{textcolor}\rmfamily\fontsize{10.000000}{12.000000}\selectfont \(\displaystyle {0.008}\)}%
\end{pgfscope}%
\begin{pgfscope}%
\pgfsetbuttcap%
\pgfsetroundjoin%
\definecolor{currentfill}{rgb}{0.000000,0.000000,0.000000}%
\pgfsetfillcolor{currentfill}%
\pgfsetlinewidth{0.803000pt}%
\definecolor{currentstroke}{rgb}{0.000000,0.000000,0.000000}%
\pgfsetstrokecolor{currentstroke}%
\pgfsetdash{}{0pt}%
\pgfsys@defobject{currentmarker}{\pgfqpoint{-0.048611in}{0.000000in}}{\pgfqpoint{-0.000000in}{0.000000in}}{%
\pgfpathmoveto{\pgfqpoint{-0.000000in}{0.000000in}}%
\pgfpathlineto{\pgfqpoint{-0.048611in}{0.000000in}}%
\pgfusepath{stroke,fill}%
}%
\begin{pgfscope}%
\pgfsys@transformshift{0.513581in}{3.770225in}%
\pgfsys@useobject{currentmarker}{}%
\end{pgfscope}%
\end{pgfscope}%
\begin{pgfscope}%
\definecolor{textcolor}{rgb}{0.000000,0.000000,0.000000}%
\pgfsetstrokecolor{textcolor}%
\pgfsetfillcolor{textcolor}%
\pgftext[x=0.100000in, y=3.722000in, left, base]{\color{textcolor}\rmfamily\fontsize{10.000000}{12.000000}\selectfont \(\displaystyle {0.010}\)}%
\end{pgfscope}%
\begin{pgfscope}%
\pgfpathrectangle{\pgfqpoint{0.513581in}{0.320679in}}{\pgfqpoint{4.960000in}{3.696000in}}%
\pgfusepath{clip}%
\pgfsetrectcap%
\pgfsetroundjoin%
\pgfsetlinewidth{1.505625pt}%
\definecolor{currentstroke}{rgb}{0.121569,0.466667,0.705882}%
\pgfsetstrokecolor{currentstroke}%
\pgfsetdash{}{0pt}%
\pgfpathmoveto{\pgfqpoint{0.739036in}{1.271632in}}%
\pgfpathlineto{\pgfqpoint{1.061114in}{0.488679in}}%
\pgfpathlineto{\pgfqpoint{1.383192in}{0.691062in}}%
\pgfpathlineto{\pgfqpoint{1.705270in}{1.317830in}}%
\pgfpathlineto{\pgfqpoint{2.027347in}{1.934069in}}%
\pgfpathlineto{\pgfqpoint{2.349425in}{2.377526in}}%
\pgfpathlineto{\pgfqpoint{2.671503in}{2.655353in}}%
\pgfpathlineto{\pgfqpoint{2.993581in}{2.930676in}}%
\pgfpathlineto{\pgfqpoint{3.315659in}{3.146854in}}%
\pgfpathlineto{\pgfqpoint{3.637737in}{3.383875in}}%
\pgfpathlineto{\pgfqpoint{3.959815in}{3.363237in}}%
\pgfpathlineto{\pgfqpoint{4.281893in}{3.531774in}}%
\pgfpathlineto{\pgfqpoint{4.603971in}{3.848679in}}%
\pgfpathlineto{\pgfqpoint{4.926049in}{3.748656in}}%
\pgfpathlineto{\pgfqpoint{5.248127in}{3.829493in}}%
\pgfusepath{stroke}%
\end{pgfscope}%
\begin{pgfscope}%
\pgfsetrectcap%
\pgfsetmiterjoin%
\pgfsetlinewidth{0.803000pt}%
\definecolor{currentstroke}{rgb}{0.000000,0.000000,0.000000}%
\pgfsetstrokecolor{currentstroke}%
\pgfsetdash{}{0pt}%
\pgfpathmoveto{\pgfqpoint{0.513581in}{0.320679in}}%
\pgfpathlineto{\pgfqpoint{0.513581in}{4.016679in}}%
\pgfusepath{stroke}%
\end{pgfscope}%
\begin{pgfscope}%
\pgfsetrectcap%
\pgfsetmiterjoin%
\pgfsetlinewidth{0.803000pt}%
\definecolor{currentstroke}{rgb}{0.000000,0.000000,0.000000}%
\pgfsetstrokecolor{currentstroke}%
\pgfsetdash{}{0pt}%
\pgfpathmoveto{\pgfqpoint{5.473581in}{0.320679in}}%
\pgfpathlineto{\pgfqpoint{5.473581in}{4.016679in}}%
\pgfusepath{stroke}%
\end{pgfscope}%
\begin{pgfscope}%
\pgfsetrectcap%
\pgfsetmiterjoin%
\pgfsetlinewidth{0.803000pt}%
\definecolor{currentstroke}{rgb}{0.000000,0.000000,0.000000}%
\pgfsetstrokecolor{currentstroke}%
\pgfsetdash{}{0pt}%
\pgfpathmoveto{\pgfqpoint{0.513581in}{0.320679in}}%
\pgfpathlineto{\pgfqpoint{5.473581in}{0.320679in}}%
\pgfusepath{stroke}%
\end{pgfscope}%
\begin{pgfscope}%
\pgfsetrectcap%
\pgfsetmiterjoin%
\pgfsetlinewidth{0.803000pt}%
\definecolor{currentstroke}{rgb}{0.000000,0.000000,0.000000}%
\pgfsetstrokecolor{currentstroke}%
\pgfsetdash{}{0pt}%
\pgfpathmoveto{\pgfqpoint{0.513581in}{4.016679in}}%
\pgfpathlineto{\pgfqpoint{5.473581in}{4.016679in}}%
\pgfusepath{stroke}%
\end{pgfscope}%
\end{pgfpicture}%
\makeatother%
\endgroup%

%% file: mc_estimate_risk_unseen_double_zipf.pgf
\begingroup%
\makeatletter%
\begin{pgfpicture}%
\pgfpathrectangle{\pgfpointorigin}{\pgfqpoint{5.474443in}{4.116679in}}%
\pgfusepath{use as bounding box, clip}%
\begin{pgfscope}%
\pgfsetbuttcap%
\pgfsetmiterjoin%
\definecolor{currentfill}{rgb}{1.000000,1.000000,1.000000}%
\pgfsetfillcolor{currentfill}%
\pgfsetlinewidth{0.000000pt}%
\definecolor{currentstroke}{rgb}{1.000000,1.000000,1.000000}%
\pgfsetstrokecolor{currentstroke}%
\pgfsetdash{}{0pt}%
\pgfpathmoveto{\pgfqpoint{0.000000in}{0.000000in}}%
\pgfpathlineto{\pgfqpoint{5.474443in}{0.000000in}}%
\pgfpathlineto{\pgfqpoint{5.474443in}{4.116679in}}%
\pgfpathlineto{\pgfqpoint{0.000000in}{4.116679in}}%
\pgfpathclose%
\pgfusepath{fill}%
\end{pgfscope}%
\begin{pgfscope}%
\pgfsetbuttcap%
\pgfsetmiterjoin%
\definecolor{currentfill}{rgb}{1.000000,1.000000,1.000000}%
\pgfsetfillcolor{currentfill}%
\pgfsetlinewidth{0.000000pt}%
\definecolor{currentstroke}{rgb}{0.000000,0.000000,0.000000}%
\pgfsetstrokecolor{currentstroke}%
\pgfsetstrokeopacity{0.000000}%
\pgfsetdash{}{0pt}%
\pgfpathmoveto{\pgfqpoint{0.374692in}{0.320679in}}%
\pgfpathlineto{\pgfqpoint{5.334692in}{0.320679in}}%
\pgfpathlineto{\pgfqpoint{5.334692in}{4.016679in}}%
\pgfpathlineto{\pgfqpoint{0.374692in}{4.016679in}}%
\pgfpathclose%
\pgfusepath{fill}%
\end{pgfscope}%
\begin{pgfscope}%
\pgfsetbuttcap%
\pgfsetroundjoin%
\definecolor{currentfill}{rgb}{0.000000,0.000000,0.000000}%
\pgfsetfillcolor{currentfill}%
\pgfsetlinewidth{0.803000pt}%
\definecolor{currentstroke}{rgb}{0.000000,0.000000,0.000000}%
\pgfsetstrokecolor{currentstroke}%
\pgfsetdash{}{0pt}%
\pgfsys@defobject{currentmarker}{\pgfqpoint{0.000000in}{-0.048611in}}{\pgfqpoint{0.000000in}{0.000000in}}{%
\pgfpathmoveto{\pgfqpoint{0.000000in}{0.000000in}}%
\pgfpathlineto{\pgfqpoint{0.000000in}{-0.048611in}}%
\pgfusepath{stroke,fill}%
}%
\begin{pgfscope}%
\pgfsys@transformshift{0.439108in}{0.320679in}%
\pgfsys@useobject{currentmarker}{}%
\end{pgfscope}%
\end{pgfscope}%
\begin{pgfscope}%
\definecolor{textcolor}{rgb}{0.000000,0.000000,0.000000}%
\pgfsetstrokecolor{textcolor}%
\pgfsetfillcolor{textcolor}%
\pgftext[x=0.439108in,y=0.223457in,,top]{\color{textcolor}\rmfamily\fontsize{10.000000}{12.000000}\selectfont \(\displaystyle {0}\)}%
\end{pgfscope}%
\begin{pgfscope}%
\pgfsetbuttcap%
\pgfsetroundjoin%
\definecolor{currentfill}{rgb}{0.000000,0.000000,0.000000}%
\pgfsetfillcolor{currentfill}%
\pgfsetlinewidth{0.803000pt}%
\definecolor{currentstroke}{rgb}{0.000000,0.000000,0.000000}%
\pgfsetstrokecolor{currentstroke}%
\pgfsetdash{}{0pt}%
\pgfsys@defobject{currentmarker}{\pgfqpoint{0.000000in}{-0.048611in}}{\pgfqpoint{0.000000in}{0.000000in}}{%
\pgfpathmoveto{\pgfqpoint{0.000000in}{0.000000in}}%
\pgfpathlineto{\pgfqpoint{0.000000in}{-0.048611in}}%
\pgfusepath{stroke,fill}%
}%
\begin{pgfscope}%
\pgfsys@transformshift{1.244302in}{0.320679in}%
\pgfsys@useobject{currentmarker}{}%
\end{pgfscope}%
\end{pgfscope}%
\begin{pgfscope}%
\definecolor{textcolor}{rgb}{0.000000,0.000000,0.000000}%
\pgfsetstrokecolor{textcolor}%
\pgfsetfillcolor{textcolor}%
\pgftext[x=1.244302in,y=0.223457in,,top]{\color{textcolor}\rmfamily\fontsize{10.000000}{12.000000}\selectfont \(\displaystyle {50}\)}%
\end{pgfscope}%
\begin{pgfscope}%
\pgfsetbuttcap%
\pgfsetroundjoin%
\definecolor{currentfill}{rgb}{0.000000,0.000000,0.000000}%
\pgfsetfillcolor{currentfill}%
\pgfsetlinewidth{0.803000pt}%
\definecolor{currentstroke}{rgb}{0.000000,0.000000,0.000000}%
\pgfsetstrokecolor{currentstroke}%
\pgfsetdash{}{0pt}%
\pgfsys@defobject{currentmarker}{\pgfqpoint{0.000000in}{-0.048611in}}{\pgfqpoint{0.000000in}{0.000000in}}{%
\pgfpathmoveto{\pgfqpoint{0.000000in}{0.000000in}}%
\pgfpathlineto{\pgfqpoint{0.000000in}{-0.048611in}}%
\pgfusepath{stroke,fill}%
}%
\begin{pgfscope}%
\pgfsys@transformshift{2.049497in}{0.320679in}%
\pgfsys@useobject{currentmarker}{}%
\end{pgfscope}%
\end{pgfscope}%
\begin{pgfscope}%
\definecolor{textcolor}{rgb}{0.000000,0.000000,0.000000}%
\pgfsetstrokecolor{textcolor}%
\pgfsetfillcolor{textcolor}%
\pgftext[x=2.049497in,y=0.223457in,,top]{\color{textcolor}\rmfamily\fontsize{10.000000}{12.000000}\selectfont \(\displaystyle {100}\)}%
\end{pgfscope}%
\begin{pgfscope}%
\pgfsetbuttcap%
\pgfsetroundjoin%
\definecolor{currentfill}{rgb}{0.000000,0.000000,0.000000}%
\pgfsetfillcolor{currentfill}%
\pgfsetlinewidth{0.803000pt}%
\definecolor{currentstroke}{rgb}{0.000000,0.000000,0.000000}%
\pgfsetstrokecolor{currentstroke}%
\pgfsetdash{}{0pt}%
\pgfsys@defobject{currentmarker}{\pgfqpoint{0.000000in}{-0.048611in}}{\pgfqpoint{0.000000in}{0.000000in}}{%
\pgfpathmoveto{\pgfqpoint{0.000000in}{0.000000in}}%
\pgfpathlineto{\pgfqpoint{0.000000in}{-0.048611in}}%
\pgfusepath{stroke,fill}%
}%
\begin{pgfscope}%
\pgfsys@transformshift{2.854692in}{0.320679in}%
\pgfsys@useobject{currentmarker}{}%
\end{pgfscope}%
\end{pgfscope}%
\begin{pgfscope}%
\definecolor{textcolor}{rgb}{0.000000,0.000000,0.000000}%
\pgfsetstrokecolor{textcolor}%
\pgfsetfillcolor{textcolor}%
\pgftext[x=2.854692in,y=0.223457in,,top]{\color{textcolor}\rmfamily\fontsize{10.000000}{12.000000}\selectfont \(\displaystyle {150}\)}%
\end{pgfscope}%
\begin{pgfscope}%
\pgfsetbuttcap%
\pgfsetroundjoin%
\definecolor{currentfill}{rgb}{0.000000,0.000000,0.000000}%
\pgfsetfillcolor{currentfill}%
\pgfsetlinewidth{0.803000pt}%
\definecolor{currentstroke}{rgb}{0.000000,0.000000,0.000000}%
\pgfsetstrokecolor{currentstroke}%
\pgfsetdash{}{0pt}%
\pgfsys@defobject{currentmarker}{\pgfqpoint{0.000000in}{-0.048611in}}{\pgfqpoint{0.000000in}{0.000000in}}{%
\pgfpathmoveto{\pgfqpoint{0.000000in}{0.000000in}}%
\pgfpathlineto{\pgfqpoint{0.000000in}{-0.048611in}}%
\pgfusepath{stroke,fill}%
}%
\begin{pgfscope}%
\pgfsys@transformshift{3.659887in}{0.320679in}%
\pgfsys@useobject{currentmarker}{}%
\end{pgfscope}%
\end{pgfscope}%
\begin{pgfscope}%
\definecolor{textcolor}{rgb}{0.000000,0.000000,0.000000}%
\pgfsetstrokecolor{textcolor}%
\pgfsetfillcolor{textcolor}%
\pgftext[x=3.659887in,y=0.223457in,,top]{\color{textcolor}\rmfamily\fontsize{10.000000}{12.000000}\selectfont \(\displaystyle {200}\)}%
\end{pgfscope}%
\begin{pgfscope}%
\pgfsetbuttcap%
\pgfsetroundjoin%
\definecolor{currentfill}{rgb}{0.000000,0.000000,0.000000}%
\pgfsetfillcolor{currentfill}%
\pgfsetlinewidth{0.803000pt}%
\definecolor{currentstroke}{rgb}{0.000000,0.000000,0.000000}%
\pgfsetstrokecolor{currentstroke}%
\pgfsetdash{}{0pt}%
\pgfsys@defobject{currentmarker}{\pgfqpoint{0.000000in}{-0.048611in}}{\pgfqpoint{0.000000in}{0.000000in}}{%
\pgfpathmoveto{\pgfqpoint{0.000000in}{0.000000in}}%
\pgfpathlineto{\pgfqpoint{0.000000in}{-0.048611in}}%
\pgfusepath{stroke,fill}%
}%
\begin{pgfscope}%
\pgfsys@transformshift{4.465082in}{0.320679in}%
\pgfsys@useobject{currentmarker}{}%
\end{pgfscope}%
\end{pgfscope}%
\begin{pgfscope}%
\definecolor{textcolor}{rgb}{0.000000,0.000000,0.000000}%
\pgfsetstrokecolor{textcolor}%
\pgfsetfillcolor{textcolor}%
\pgftext[x=4.465082in,y=0.223457in,,top]{\color{textcolor}\rmfamily\fontsize{10.000000}{12.000000}\selectfont \(\displaystyle {250}\)}%
\end{pgfscope}%
\begin{pgfscope}%
\pgfsetbuttcap%
\pgfsetroundjoin%
\definecolor{currentfill}{rgb}{0.000000,0.000000,0.000000}%
\pgfsetfillcolor{currentfill}%
\pgfsetlinewidth{0.803000pt}%
\definecolor{currentstroke}{rgb}{0.000000,0.000000,0.000000}%
\pgfsetstrokecolor{currentstroke}%
\pgfsetdash{}{0pt}%
\pgfsys@defobject{currentmarker}{\pgfqpoint{0.000000in}{-0.048611in}}{\pgfqpoint{0.000000in}{0.000000in}}{%
\pgfpathmoveto{\pgfqpoint{0.000000in}{0.000000in}}%
\pgfpathlineto{\pgfqpoint{0.000000in}{-0.048611in}}%
\pgfusepath{stroke,fill}%
}%
\begin{pgfscope}%
\pgfsys@transformshift{5.270276in}{0.320679in}%
\pgfsys@useobject{currentmarker}{}%
\end{pgfscope}%
\end{pgfscope}%
\begin{pgfscope}%
\definecolor{textcolor}{rgb}{0.000000,0.000000,0.000000}%
\pgfsetstrokecolor{textcolor}%
\pgfsetfillcolor{textcolor}%
\pgftext[x=5.270276in,y=0.223457in,,top]{\color{textcolor}\rmfamily\fontsize{10.000000}{12.000000}\selectfont \(\displaystyle {300}\)}%
\end{pgfscope}%
\begin{pgfscope}%
\pgfsetbuttcap%
\pgfsetroundjoin%
\definecolor{currentfill}{rgb}{0.000000,0.000000,0.000000}%
\pgfsetfillcolor{currentfill}%
\pgfsetlinewidth{0.803000pt}%
\definecolor{currentstroke}{rgb}{0.000000,0.000000,0.000000}%
\pgfsetstrokecolor{currentstroke}%
\pgfsetdash{}{0pt}%
\pgfsys@defobject{currentmarker}{\pgfqpoint{-0.048611in}{0.000000in}}{\pgfqpoint{-0.000000in}{0.000000in}}{%
\pgfpathmoveto{\pgfqpoint{-0.000000in}{0.000000in}}%
\pgfpathlineto{\pgfqpoint{-0.048611in}{0.000000in}}%
\pgfusepath{stroke,fill}%
}%
\begin{pgfscope}%
\pgfsys@transformshift{0.374692in}{0.346893in}%
\pgfsys@useobject{currentmarker}{}%
\end{pgfscope}%
\end{pgfscope}%
\begin{pgfscope}%
\definecolor{textcolor}{rgb}{0.000000,0.000000,0.000000}%
\pgfsetstrokecolor{textcolor}%
\pgfsetfillcolor{textcolor}%
\pgftext[x=0.100000in, y=0.298668in, left, base]{\color{textcolor}\rmfamily\fontsize{10.000000}{12.000000}\selectfont \(\displaystyle {0.0}\)}%
\end{pgfscope}%
\begin{pgfscope}%
\pgfsetbuttcap%
\pgfsetroundjoin%
\definecolor{currentfill}{rgb}{0.000000,0.000000,0.000000}%
\pgfsetfillcolor{currentfill}%
\pgfsetlinewidth{0.803000pt}%
\definecolor{currentstroke}{rgb}{0.000000,0.000000,0.000000}%
\pgfsetstrokecolor{currentstroke}%
\pgfsetdash{}{0pt}%
\pgfsys@defobject{currentmarker}{\pgfqpoint{-0.048611in}{0.000000in}}{\pgfqpoint{-0.000000in}{0.000000in}}{%
\pgfpathmoveto{\pgfqpoint{-0.000000in}{0.000000in}}%
\pgfpathlineto{\pgfqpoint{-0.048611in}{0.000000in}}%
\pgfusepath{stroke,fill}%
}%
\begin{pgfscope}%
\pgfsys@transformshift{0.374692in}{0.894283in}%
\pgfsys@useobject{currentmarker}{}%
\end{pgfscope}%
\end{pgfscope}%
\begin{pgfscope}%
\definecolor{textcolor}{rgb}{0.000000,0.000000,0.000000}%
\pgfsetstrokecolor{textcolor}%
\pgfsetfillcolor{textcolor}%
\pgftext[x=0.100000in, y=0.846058in, left, base]{\color{textcolor}\rmfamily\fontsize{10.000000}{12.000000}\selectfont \(\displaystyle {0.5}\)}%
\end{pgfscope}%
\begin{pgfscope}%
\pgfsetbuttcap%
\pgfsetroundjoin%
\definecolor{currentfill}{rgb}{0.000000,0.000000,0.000000}%
\pgfsetfillcolor{currentfill}%
\pgfsetlinewidth{0.803000pt}%
\definecolor{currentstroke}{rgb}{0.000000,0.000000,0.000000}%
\pgfsetstrokecolor{currentstroke}%
\pgfsetdash{}{0pt}%
\pgfsys@defobject{currentmarker}{\pgfqpoint{-0.048611in}{0.000000in}}{\pgfqpoint{-0.000000in}{0.000000in}}{%
\pgfpathmoveto{\pgfqpoint{-0.000000in}{0.000000in}}%
\pgfpathlineto{\pgfqpoint{-0.048611in}{0.000000in}}%
\pgfusepath{stroke,fill}%
}%
\begin{pgfscope}%
\pgfsys@transformshift{0.374692in}{1.441674in}%
\pgfsys@useobject{currentmarker}{}%
\end{pgfscope}%
\end{pgfscope}%
\begin{pgfscope}%
\definecolor{textcolor}{rgb}{0.000000,0.000000,0.000000}%
\pgfsetstrokecolor{textcolor}%
\pgfsetfillcolor{textcolor}%
\pgftext[x=0.100000in, y=1.393449in, left, base]{\color{textcolor}\rmfamily\fontsize{10.000000}{12.000000}\selectfont \(\displaystyle {1.0}\)}%
\end{pgfscope}%
\begin{pgfscope}%
\pgfsetbuttcap%
\pgfsetroundjoin%
\definecolor{currentfill}{rgb}{0.000000,0.000000,0.000000}%
\pgfsetfillcolor{currentfill}%
\pgfsetlinewidth{0.803000pt}%
\definecolor{currentstroke}{rgb}{0.000000,0.000000,0.000000}%
\pgfsetstrokecolor{currentstroke}%
\pgfsetdash{}{0pt}%
\pgfsys@defobject{currentmarker}{\pgfqpoint{-0.048611in}{0.000000in}}{\pgfqpoint{-0.000000in}{0.000000in}}{%
\pgfpathmoveto{\pgfqpoint{-0.000000in}{0.000000in}}%
\pgfpathlineto{\pgfqpoint{-0.048611in}{0.000000in}}%
\pgfusepath{stroke,fill}%
}%
\begin{pgfscope}%
\pgfsys@transformshift{0.374692in}{1.989064in}%
\pgfsys@useobject{currentmarker}{}%
\end{pgfscope}%
\end{pgfscope}%
\begin{pgfscope}%
\definecolor{textcolor}{rgb}{0.000000,0.000000,0.000000}%
\pgfsetstrokecolor{textcolor}%
\pgfsetfillcolor{textcolor}%
\pgftext[x=0.100000in, y=1.940839in, left, base]{\color{textcolor}\rmfamily\fontsize{10.000000}{12.000000}\selectfont \(\displaystyle {1.5}\)}%
\end{pgfscope}%
\begin{pgfscope}%
\pgfsetbuttcap%
\pgfsetroundjoin%
\definecolor{currentfill}{rgb}{0.000000,0.000000,0.000000}%
\pgfsetfillcolor{currentfill}%
\pgfsetlinewidth{0.803000pt}%
\definecolor{currentstroke}{rgb}{0.000000,0.000000,0.000000}%
\pgfsetstrokecolor{currentstroke}%
\pgfsetdash{}{0pt}%
\pgfsys@defobject{currentmarker}{\pgfqpoint{-0.048611in}{0.000000in}}{\pgfqpoint{-0.000000in}{0.000000in}}{%
\pgfpathmoveto{\pgfqpoint{-0.000000in}{0.000000in}}%
\pgfpathlineto{\pgfqpoint{-0.048611in}{0.000000in}}%
\pgfusepath{stroke,fill}%
}%
\begin{pgfscope}%
\pgfsys@transformshift{0.374692in}{2.536455in}%
\pgfsys@useobject{currentmarker}{}%
\end{pgfscope}%
\end{pgfscope}%
\begin{pgfscope}%
\definecolor{textcolor}{rgb}{0.000000,0.000000,0.000000}%
\pgfsetstrokecolor{textcolor}%
\pgfsetfillcolor{textcolor}%
\pgftext[x=0.100000in, y=2.488230in, left, base]{\color{textcolor}\rmfamily\fontsize{10.000000}{12.000000}\selectfont \(\displaystyle {2.0}\)}%
\end{pgfscope}%
\begin{pgfscope}%
\pgfsetbuttcap%
\pgfsetroundjoin%
\definecolor{currentfill}{rgb}{0.000000,0.000000,0.000000}%
\pgfsetfillcolor{currentfill}%
\pgfsetlinewidth{0.803000pt}%
\definecolor{currentstroke}{rgb}{0.000000,0.000000,0.000000}%
\pgfsetstrokecolor{currentstroke}%
\pgfsetdash{}{0pt}%
\pgfsys@defobject{currentmarker}{\pgfqpoint{-0.048611in}{0.000000in}}{\pgfqpoint{-0.000000in}{0.000000in}}{%
\pgfpathmoveto{\pgfqpoint{-0.000000in}{0.000000in}}%
\pgfpathlineto{\pgfqpoint{-0.048611in}{0.000000in}}%
\pgfusepath{stroke,fill}%
}%
\begin{pgfscope}%
\pgfsys@transformshift{0.374692in}{3.083845in}%
\pgfsys@useobject{currentmarker}{}%
\end{pgfscope}%
\end{pgfscope}%
\begin{pgfscope}%
\definecolor{textcolor}{rgb}{0.000000,0.000000,0.000000}%
\pgfsetstrokecolor{textcolor}%
\pgfsetfillcolor{textcolor}%
\pgftext[x=0.100000in, y=3.035620in, left, base]{\color{textcolor}\rmfamily\fontsize{10.000000}{12.000000}\selectfont \(\displaystyle {2.5}\)}%
\end{pgfscope}%
\begin{pgfscope}%
\pgfsetbuttcap%
\pgfsetroundjoin%
\definecolor{currentfill}{rgb}{0.000000,0.000000,0.000000}%
\pgfsetfillcolor{currentfill}%
\pgfsetlinewidth{0.803000pt}%
\definecolor{currentstroke}{rgb}{0.000000,0.000000,0.000000}%
\pgfsetstrokecolor{currentstroke}%
\pgfsetdash{}{0pt}%
\pgfsys@defobject{currentmarker}{\pgfqpoint{-0.048611in}{0.000000in}}{\pgfqpoint{-0.000000in}{0.000000in}}{%
\pgfpathmoveto{\pgfqpoint{-0.000000in}{0.000000in}}%
\pgfpathlineto{\pgfqpoint{-0.048611in}{0.000000in}}%
\pgfusepath{stroke,fill}%
}%
\begin{pgfscope}%
\pgfsys@transformshift{0.374692in}{3.631236in}%
\pgfsys@useobject{currentmarker}{}%
\end{pgfscope}%
\end{pgfscope}%
\begin{pgfscope}%
\definecolor{textcolor}{rgb}{0.000000,0.000000,0.000000}%
\pgfsetstrokecolor{textcolor}%
\pgfsetfillcolor{textcolor}%
\pgftext[x=0.100000in, y=3.583011in, left, base]{\color{textcolor}\rmfamily\fontsize{10.000000}{12.000000}\selectfont \(\displaystyle {3.0}\)}%
\end{pgfscope}%
\begin{pgfscope}%
\pgfpathrectangle{\pgfqpoint{0.374692in}{0.320679in}}{\pgfqpoint{4.960000in}{3.696000in}}%
\pgfusepath{clip}%
\pgfsetrectcap%
\pgfsetroundjoin%
\pgfsetlinewidth{1.505625pt}%
\definecolor{currentstroke}{rgb}{0.121569,0.466667,0.705882}%
\pgfsetstrokecolor{currentstroke}%
\pgfsetdash{}{0pt}%
\pgfpathmoveto{\pgfqpoint{0.600146in}{0.488679in}}%
\pgfpathlineto{\pgfqpoint{0.922224in}{0.710002in}}%
\pgfpathlineto{\pgfqpoint{1.244302in}{2.443480in}}%
\pgfpathlineto{\pgfqpoint{1.566380in}{3.539528in}}%
\pgfpathlineto{\pgfqpoint{1.888458in}{3.848679in}}%
\pgfpathlineto{\pgfqpoint{2.210536in}{3.158876in}}%
\pgfpathlineto{\pgfqpoint{2.532614in}{2.134621in}}%
\pgfpathlineto{\pgfqpoint{2.854692in}{1.694686in}}%
\pgfpathlineto{\pgfqpoint{3.176770in}{1.345973in}}%
\pgfpathlineto{\pgfqpoint{3.498848in}{1.161071in}}%
\pgfpathlineto{\pgfqpoint{3.820926in}{1.077521in}}%
\pgfpathlineto{\pgfqpoint{4.143004in}{0.844842in}}%
\pgfpathlineto{\pgfqpoint{4.465082in}{0.727138in}}%
\pgfpathlineto{\pgfqpoint{4.787159in}{0.722813in}}%
\pgfpathlineto{\pgfqpoint{5.109237in}{0.700810in}}%
\pgfusepath{stroke}%
\end{pgfscope}%
\begin{pgfscope}%
\pgfsetrectcap%
\pgfsetmiterjoin%
\pgfsetlinewidth{0.803000pt}%
\definecolor{currentstroke}{rgb}{0.000000,0.000000,0.000000}%
\pgfsetstrokecolor{currentstroke}%
\pgfsetdash{}{0pt}%
\pgfpathmoveto{\pgfqpoint{0.374692in}{0.320679in}}%
\pgfpathlineto{\pgfqpoint{0.374692in}{4.016679in}}%
\pgfusepath{stroke}%
\end{pgfscope}%
\begin{pgfscope}%
\pgfsetrectcap%
\pgfsetmiterjoin%
\pgfsetlinewidth{0.803000pt}%
\definecolor{currentstroke}{rgb}{0.000000,0.000000,0.000000}%
\pgfsetstrokecolor{currentstroke}%
\pgfsetdash{}{0pt}%
\pgfpathmoveto{\pgfqpoint{5.334692in}{0.320679in}}%
\pgfpathlineto{\pgfqpoint{5.334692in}{4.016679in}}%
\pgfusepath{stroke}%
\end{pgfscope}%
\begin{pgfscope}%
\pgfsetrectcap%
\pgfsetmiterjoin%
\pgfsetlinewidth{0.803000pt}%
\definecolor{currentstroke}{rgb}{0.000000,0.000000,0.000000}%
\pgfsetstrokecolor{currentstroke}%
\pgfsetdash{}{0pt}%
\pgfpathmoveto{\pgfqpoint{0.374692in}{0.320679in}}%
\pgfpathlineto{\pgfqpoint{5.334692in}{0.320679in}}%
\pgfusepath{stroke}%
\end{pgfscope}%
\begin{pgfscope}%
\pgfsetrectcap%
\pgfsetmiterjoin%
\pgfsetlinewidth{0.803000pt}%
\definecolor{currentstroke}{rgb}{0.000000,0.000000,0.000000}%
\pgfsetstrokecolor{currentstroke}%
\pgfsetdash{}{0pt}%
\pgfpathmoveto{\pgfqpoint{0.374692in}{4.016679in}}%
\pgfpathlineto{\pgfqpoint{5.334692in}{4.016679in}}%
\pgfusepath{stroke}%
\end{pgfscope}%
\end{pgfpicture}%
\makeatother%
\endgroup%